\documentclass{gtart} 

\input gtoutput
\volumenumber{4}\papernumber{1}\volumeyear{2000}
\pagenumbers{1}{83}\published{28 January 2000}
\proposed{Frances Kirwan}\seconded{Joan Birman, Robion Kirby}
\received{30 October 1999}\revised{27 January 2000}
\accepted{14 January 1999}

\usepackage{amssymb,amsmath,amscd}
\usepackage{graphicx}

\theoremstyle{plain}
\newtheorem{theorem}{Theorem}[section]
\newtheorem{lemma}[theorem]{Lemma}
\newtheorem{corollary}[theorem]{Corollary}
\newtheorem{proposition}[theorem]{Proposition}
\newtheorem{conjecture}[theorem]{Conjecture}

\theoremstyle{definition}
\newtheorem{definition}[theorem]{Definition}

\newtheorem{convention}[theorem]{Convention}

\theoremstyle{remark}
\newtheorem{remark}[theorem]{Remark}
\newtheorem*{claim}{Claim}
\newtheorem*{acknowledgement}{Acknowledgements}

\newcommand\C{{\mathcal C}}
\newcommand\f{{\mathcal F}}
\newcommand\calL{{\mathcal L}}
\newcommand\m{{\mathcal M}}
\newcommand\bfM{{\mathbf M}}
\newcommand\q{{\mathbf Q}}
\newcommand\z{{\mathbf Z}}

\begin{document}


\title{Claspers and finite type invariants of links} %
\authors{Kazuo Habiro} 
\address{Graduate School of Mathematical Sciences, University
of Tokyo\\ 3--8--1 Komaba Meguro-ku, Tokyo 153, Japan}
\asciiaddress{Graduate School of Mathematical Sciences, University
of Tokyo\\ 3-8-1 Komaba Meguro-ku, Tokyo 153, Japan}
\email{habiro@ms.u-tokyo.ac.jp} 

\begin{abstract}
We introduce the concept of ``claspers,'' which are surfaces in
$3$--manifolds with some additional structure on which surgery
operations can be performed.  Using claspers we define for each
positive integer $k$ an equivalence relation on links called
``$C_k$--equivalence,'' which is generated by surgery operations of a
certain kind called ``$C_k$--moves''.  We prove that two knots in the
$3$--sphere are $C_{k+1}$--equivalent if and only if they have equal
values of Vassiliev--Goussarov invariants of type $k$ with values in
any abelian groups.  This result gives a characterization in terms of
surgery operations of the informations that can be carried by
Vassiliev--Goussarov invariants.  In the last section we also describe
outlines of some applications of claspers to other fields in
$3$--dimensional topology.
\end{abstract}
\asciiabstract{We introduce the concept of `claspers,' which are
surfaces in 3-manifolds with some additional structure on which
surgery operations can be performed.  Using claspers we define for
each positive integer k an equivalence relation on links called
`C_k-equivalence,' which is generated by surgery operations of a
certain kind called `C_k-moves'.  We prove that two knots in the
3-sphere are C_{k+1}-equivalent if and only if they have equal
values of Vassiliev-Goussarov invariants of type k with values in
any abelian groups.  This result gives a characterization in terms of
surgery operations of the informations that can be carried by
Vassiliev--Goussarov invariants.  In the last section we also describe
outlines of some applications of claspers to other fields in
3-dimensional topology.}


\maketitlepage
%
%


\section{Introduction}

In the theory of finite type invariants of knots and links, also
called Vassiliev--Goussarov invariants \cite{Vassiliev:Cohomology}
\cite{Goussarov:ANewForm} \cite{Birman:NewPoints}
\cite{Birman-Lin:KnotPolynomials} \cite{Bar-Natan:OnTheVassiliev}
\cite{Kontsevich:VassilievsKnotInvariants}, we have a descending
filtration, called the Vassiliev--Goussarov filtration, on the free
abelian group generated by ambient isotopy classes of links, and
dually an ascending filtration on the group of invariants of links
with values in an abelian group.  Invariants which lie in the $k$th
subgroup in the filtration are characterized by the property that they
vanish on the $k+1$st subgroup of the Vassiliev--Goussarov filtration,
and called invariants of type $k$.

It is natural to ask when the difference of two links lies in the
$k+1$st subgroup of the Vassiliev--Goussarov filtration, ie, when the
two links are not distinguished by any invariant of type $k$.  If this
is the case, then the two links are said to be ``$V_k$--equivalent.''
T~Stanford proved in \cite{Stanford:BraidCommutators} that two links
are $V_k$--equivalent if one is obtained from the other by inserting a
pure braid commutator of class $k+1$.  
One of the main purposes of this paper is to prove a modified version
of the converse of this result:
\begin{theorem}
\label{main}
For two knots $\gamma$ and $\gamma'$ in $S^3$ and for $k\ge 0$, the
following conditions are equivalent.
\begin{enumerate}
\item $\gamma$ and $\gamma'$ are $V_k$--equivalent.
\item $\gamma$ and $\gamma'$ are related by an element of the $k+1$st
lower central series subgroup (ie, the subgroup generated by the
iterated commutators of class $k+1$) of the pure braid group of $n$
strands for some $n\ge 0$.  
\item $\gamma$ and $\gamma'$ are related by a finite sequence of
``simple $C_k$--moves'' and ambient isotopies.
\end{enumerate}
Here a ``simple $C_k$--move'' is a local operation on knots defined
using ``claspers''.  (Loosely speaking, a simple ${C_k}$--move on a link
is an operation which ``band-sums a $(k+1)$--component iterated Bing
double of the Hopf link.''  See Figure~\ref{exCk} for the case that $k=1$,
$2$ and $3$.) 
\end{theorem}

Theorem~\ref{main} is a part of Theorem~\ref{JkCk}.
M Goussarov independently proved a similar result.  Recently,
T~Stanford proved (after an earlier version
\cite{Habiro:ClaspersAndTheVassilievSkeinModules} of the present
paper, in which the equivalence of (1) and (3) of Theorem~\ref{main}
was proved, was circulated) that two knots in $S^3$ are
$V_k$--equivalent if and only if they are presented as two closed
braids differing only by an element of the $k+1$st lower central
series subgroup of the corresponding pure braid group
\cite{Stanford:VassilievInvariantsAndKnots}.  The equivalence of 1 and
2 in the above theorem can be derived also from this result of
Stanford.  His proof seems to be simpler than ours in some respects,
mostly due to the use of commutator calculus in groups, which is well
developed in literature.  However, we believe that it is worth
presenting the proof using claspers here because we think of our
technique, {\em calculus of claspers}, as a calculus of a new kind in
$3$--dimensional topology which plays a fundamental role in studying
finite type invariants of links and $3$--manifolds and, moreover, in
studying the category theoretic and algebraic structures in
$3$--dimensional topology.
\vspace{6pt}

Calculus of claspers is closely related to three well-known calculi:
Kirby's calculus of framed links \cite{Kirby:ACalculus}, the diagram
calculus of morphisms in braided categories
\cite{Majid:AlgebrasAndHopfAlgebras}, and the calculus of trivalent
graphs appearing in theories of finite type invariants of links and
$3$--manifolds \cite{Bar-Natan:OnTheVassiliev}
\cite{Garoufalidis-Ohtsuki:OnFiniteTypeIII}.  Let us briefly explain
these relationships here.

First, we may think of calculus of claspers as a variant of Kirby's
calculus of framed links \cite{Kirby:ACalculus}.  The Kirby calculus
reduces, to some extent, the study of closed oriented $3$--manifolds to
the study of framed links in $S^3$.  Claspers are topological objects
in $3$--manifolds on which we can perform surgery, like framed links.
In fact, surgery on a clasper is defined as surgery on an ``associated
framed link''.  Therefore we may think of calculus of claspers as
calculus of framed links of a special kind.\footnote{We can easily
derive from Kirby's theorem a set of operations on claspers that
generate the equivalence relation which says when two claspers yield
diffeomorphic results of surgeries.  But these moves seems to be not
so interesting.  An interesting version of ``Kirby type theorem''
would be equivalent to a presentation of the braided category
${\mathbf{Cob}}$ described just below.}

Second, we may think of calculus of claspers as a kind of diagram
calculus for a category ${\mathbf{Cob}}$ {\em embedded in a
$3$--manifold}.  Here ${\mathbf{Cob}}$ denotes the rigid braided strict
monoidal category of cobordisms of oriented connected surfaces with
connected boundary (see \cite{Crane-Yetter:OnAlgebraicStructures} or
\cite{Kerler:Genealogy}).  Recall that ${\mathbf{Cob}}$ is generated
as a braided category by the ``handle Hopf algebra,'' which is a
punctured torus as an object of ${\mathbf{Cob}}$.  Recall also that in
diagram calculus for braided category, an object is represented by a
vertical line or a parallel family of some vertical lines, and a
morphism by a vertex which have some input lines corresponding to the
domain and some output lines the codomain (see, eg,
\cite{Majid:FoundationsOfQuantum}).  If the braided category in
question is the cobordism category ${\mathbf{Cob}}$, then a diagram
represents a cobordism.  Speaking roughly and somewhat inaccurately, a
clasper is a flexible generalization of such a diagram embedded in a
$3$--manifold and we can perform {\em surgery} on it, which means
removing a regular neighborhood of it and gluing back the cobordism
represented by the diagram.  In this way, we may sometimes think of (a
part of) a clasper as a diagram in ${\mathbf{Cob}}$.  This enables us
to think of claspers {\em algebraically}.

Third, calculus of claspers is a kind of ``topological version'' of
the calculus of uni-trivalent graphs which appear in theories of
finite type invariants of links and $3$--manifolds
\cite{Bar-Natan:OnTheVassiliev}
\cite{Garoufalidis-Ohtsuki:OnFiniteTypeIII}.  Claspers of a special
kind, which we call ``(simple) graph claspers'' look very like
trivalent graphs, but they are embedded in a $3$--manifold and have
framings on edges.  We can think of a graph clasper as a ``topological
realization'' of a trivalent graph.  This aspect of calculus of
claspers is very important in that it gives an unifying view on finite
type invariants of both links and $3$--manifolds.  Moreover, we can
develop theories of clasper surgery equivalence relations on links and
$3$--manifolds.  We may think of this theory as more fundamental than
that of finite type invariants.

From the category theoretical point of view explained above, we may
think of the aspect of calculus of claspers related to trivalent
graphs as {\em commutator calculus in the braided category ${\mathbf{Cob}}$}.
This point of view clarifies that {\em the Lie algebraic structure of
trivalent graphs originates from the Hopf algebraic structure in the
category ${\mathbf{Cob}}$}.  This observation is just like that the Lie algebra
structure of the associated graded of the lower central series of a
group is explained in terms of the group structure.
\vspace{6pt}

The organization of the rest of this paper is as follows.
Sections~2--7 are devoted to definitions of claspers and theories of
$C_k$--equivalence relations and finite type invariants of
links.  Section~8 is devoted to giving a survey on other theories
stemming from calculus of claspers.

In section 2, we define the notion of claspers. A {\em basic
clasper\/} in an oriented $3$--manifold $M$ is a planar surface with
$3$ boundary components embedded in the interior of $M$ equipped with
a decomposition into two annuli and a band.  For a basic clasper $C$,
we associate a 2--component framed link $L_C$, and we define ``surgery
on a basic clasper $C$'' as surgery on the associated framed link
$L_C$.  Basic claspers serve as building blocks of claspers.  A {\em
clasper} in $M$ is a surface embedded in the interior of $M$
decomposed into some subsurfaces.  We associate to a clasper a union
of basic claspers in a certain way and we define surgery on the
clasper $G$ as surgery on associated basic claspers.  A {\em tame}
clasper is a clasper on which the surgery does not change the
$3$--manifold up to a canonical diffeomorphism.  We give some
moves on claspers and links which does not change the results of
surgeries (Proposition~\ref{moves}).

In section 3, we define {\em strict tree claspers}, which are tame
claspers of a special kind.  We define the notion of $C_k$--moves on
links as surgery on a strict tree clasper of degree $k$.  The
${C_k}$--equivalence is generated by $C_k$--moves and ambient isotopies.
The ${C_k}$--equivalence relation becomes finer as $k$ increases
(Proposition~\ref{leqkeq}).  In Theorem~\ref{equivalence}, we give
some necessary and sufficient conditions that two links are
$C_k$--equivalent.

In section 4, we define the notion of {\em homotopy} of claspers with
respect to a link ${\gamma_0}$ in a $3$--manifold $M$.  If two simple
strict forest claspers of degree $k$ (ie, a union of simple strict
tree claspers of degree $k$) are homotopic to each other, then they
yield ${C_{k+1}}$--equivalent results of surgeries
(Theorem~\ref{FactorThroughFhkg}).  Moreover, a certain abelian group
maps onto the set of ${C_{k+1}}$--equivalence classes of links which
are ${C_k}$--equivalent to a fixed link ${\gamma_0}$
(Theorem~\ref{FactorThroughFthkg}).  This abelian group is finitely
generated if $\pi_1M$ is finite.

In section 5, we define a monoid $\calL(\Sigma,n)$ of $n$--string links
in $\Sigma\times [0,1]$, where $\Sigma $ is a compact connected
oriented surface, and study the quotient $\calL(\Sigma,n)/{C_{k+1}}$
by the ${C_{k+1}}$--equivalence.  The monoid
$\calL(\Sigma,n)/{C_{k+1}}$ forms a residually solvable group, and 
the subgroup $\calL_1(\Sigma,n)/{C_{k+1}}$ of
$\calL(\Sigma,n)/{C_{k+1}}$ consisting of the ${C_{k+1}}$--equivalence
classes of homotopically trivial $n$--string links forms a group
(Theorem~\ref{group}).  These groups are
finitely generated if $\Sigma $ is a disk or a sphere
(Corollary~\ref{disksphere}).  The pure braid group $P(\Sigma,n)$ of
$n$--strands in ${\Sigma\times[0,1]}$ forms the unit subgroup of the
monoid $\calL(\Sigma,n)$ of $n$--string links in ${\Sigma\times[0,1]}$.
We show that the commutators of class $k$ of the subgroup
$P_1(\Sigma,n)$ of $P(\Sigma,n)$ consisting of homotopically trivial
pure braids are ${C_k}$--equivalent to ${1_n}$
(Proposition~\ref{pksnlksn}).  Using this result, we prove that two
links in a $3$--manifold are ${C_k}$--equivalent if and only if they are
``$P'_k$--equivalent'' (ie, related by an element of the $k$th lower
central series subgroup of a pure braid group in ${D^2\times[0,1]}$)
(Theorem~\ref{CkPnk}).  We give a definition of a graded Lie algebra
of string links.

In section 6, we study Vassiliev--Goussarov filtrations using claspers.
In 6.1, we recall the usual definition of Vassiliev--Goussarov
filtrations and finite type invariants using singular links.  In 6.2,
we redefine Vassiliev--Goussarov filtrations on links using {\em forest
schemes}, which are finite sets of disjoint strict tree claspers.  In 6.3, we
restrict our attention to Vassiliev--Goussarov filtrations on string
links, and in 6.4, to that on ``string knots'', ie, $1$--string links
in ${D^2\times[0,1]}$ up to ambient isotopy.  Clearly, there is a
natural one-to-one correspondence between the set of string knots and
that of knots in $S^3$.  We define an additive invariant ${\psi_k}$ of
type $k$ of string knots with values in the group of
${C_{k+1}}$--equivalence classes of string knots.  The invariant
${\psi_k}$ is universal among the additive invariants of type $k$ of
string knots (Theorem~\ref{UniversalAdditiveVassilievInvariant}).
Using this, we prove Theorem~\ref{JkCk}, which contains
Theorem~\ref{main}.

In section 7, we give some examples.  A simple $C_k$--move is a
$C_k$--move of a special kind and can be defined also as a band-sum
operation of a $(k+1)$--component iterated Bing double of the Hopf
link.  The Milnor $\bar\mu$ invariants of length $k+1$ of links in
$S^3$ are invariants of ${C_{k+1}}$--equivalence (Theorem~\ref{mubar}).  The
${C_k}$--equivalence relation is more closely related with the Milnor
$\bar\mu$ invariants than the $V_k$--equivalence relation.

In section 8, we give a survey of some other aspects of calculus of
claspers.  In 8.1, we explain the relationships between claspers and a
category of surface cobordisms.  In 8.2, we generalize the notion of
tree claspers to ``graph claspers'' and explain that graph claspers is
regarded as {\em topological realizations} of uni-trivalent graphs.
In 8.3, we give a definition of new filtrations on links and ``special
finite type invariants'' of links.  In 8.4, we apply claspers to the
theory of finite type invariants of $3$--manifolds.  In 8.5, we define
``groups of homology cobordisms of surfaces,'' which are extensions of
certain quotient of mapping class groups.  In 8.6, we relate claspers
to embedded gropes in $3$--manifolds.
\vspace{6pt}

We remark that, after almost finishing the present paper, the author
was informed that M Goussarov has given some constructions similar to
claspers.

\medskip

\begin{acknowledgement}
The author was partially supported by Research Fellowships of the
Japan Society for the Promotion of Science for Young Scientists.

This paper is a based on my Ph.D thesis
\cite{Habiro:ClaspersAndTheVassilievSkeinModules}, and I would like to
thank my advisor Yukio Matsumoto for helpful advice and continuous
encouragement.  I also thank Mikhail Goussarov, Thang Le, Hitoshi
Murakami and Tomotada Ohtsuki for useful comments and stimulating
conversations.
\end{acknowledgement}

\subsection{Preliminaries}
Throughout this paper all manifolds are smooth, compact, connected and
oriented unless otherwise stated.  Moreover, $3$--manifolds are always
oriented, and embeddings and diffeomorphism of $3$--manifolds are
orientation-preserving.

For a $3$--manifold $M$, a {\em pattern} $P=(\alpha ,i)$ on $M$ is the
pair of a compact, oriented $1$--manifold $\alpha $ and an embedding
$i\colon \partial \alpha \hookrightarrow\partial M$.  A {\em link}
$\gamma$ in $M$ of pattern $P$ is a proper embedding of $\alpha$ into
$M$ which restricts to $i$ on boundary.  Let $\gamma$ denote also the
image.  Two links $\gamma$ and $\gamma'$ in $M$ of the same pattern
$P$ are said to be {\em equivalent} (denoted $\gamma\cong\gamma'$) if
$\gamma$ and $\gamma'$ are related by an ambient isotopy relative to
endpoints.  Let $[\gamma]$ often denote the equivalence class of a
link $\gamma$.  (In literature, a `link' usually means a finite
disjoint union of embedded {\em circles}.  However, we will work with
the above extended definition of `links' in this paper.)  We simply
say that two links of pattern $P$ are {\em homotopic} to each other if
they are homotopic to each other relative to endpoints.

A {\em framed link} will mean a link consisting of only circle
components which are equipped with framings, ie, homotopy classes of
non-zero sections of the normal bundles.  In other words, a ``framed
link'' mean a ``usual framed link''.  Surgery on a framed link is
defined in the usual way.  The result from a $3$--manifold $M$ by
surgery on a framed link $L$ in $M$ is denoted by $M^L$.

For an equivalence relation $R$ on a set $S$ and an element $s$ of
$S$, let $[s]_R$ denote the element of the quotient set $S/R$
corresponding to $s$.  Similarly, for a normal subgroup $H$ of a group
$G$ and an element $g$ of $G$, let $[g]_H$ denote the coset $gH$ of
$g$ in the quotient group $G/H$.

For a group $G$, the $k$th lower central series subgroup $G_k$ of $G$
is defined by $G_1=G$ and $G_{k+1}=[G,G_k]$ ($k\ge1$), where
$[\cdot,\cdot]$ denotes commutator subgroup.

\section{Claspers and basic claspers}

In this section we introduce the notion of claspers and basic claspers
in $3$--manifolds.  A clasper is a kind of surface embedded in a
$3$--manifold on which one may perform surgery, like framed links.  A
clasper in a $3$--manifold $M$ is said to be ``tame'' if the result of
surgery yields a $3$--manifold which is diffeomorphic to $M$ in a
canonical way.  We may use a tame clasper to transform a link in $M$
into another.  At the end of this section we introduce some operations
on claspers and links which do not change the results of surgeries.

\subsection{Basic claspers}

\begin{definition}
A {\em basic clasper} $C=A_1\cup A_2\cup B$ in a $3$--manifold $M$ is
a non-oriented planar surface embedded in $M$ with three boundary
components equipped with a decomposition into two annuli $A_1$ and
$A_2$, and a band\footnote{A ``band'' will mean a disk parametrized by
$[0,1]\times [0,1]$ such that the two arcs in the boundary
corresponding to $\{i\} \times [0,1] \quad (i=0,1)$ are attached to
the boundary of other surfaces.} $B$ connecting $A_1$ and $A_2$.  We
call the two annuli $A_1$ and $A_2$ the {\em leaves} of $C$ and the
band $B$ the {\em edge} of $C$.  

Given a basic clasper $C=A_1\cup A_2\cup B$ in $M$, we associate to it
a $2$--component framed link $L_C=L_{C,1}\cup L_{C,2}$ in a small
regular neighborhood $N_C$ of $C$ in $M$ as follows.  Let $A_1'$ and
$A_2'$ be the two annuli in $N_C$ obtained from $A_1$ and $A_2$ by a
crossing change along the band $B$ as illustrated in
Figure~\ref{AssociatedFramedLink}a.  (Here the crossing must be just
as depicted, and it must not be in the opposite way.)  The framed link
$L_C$ is unique up to isotopy.  The framed link $L_C=L_{C,1}\cup
L_{C,2}$ is determined by $A'_1$ and $A'_2$ as depicted in
Figure~\ref{AssociatedFramedLink}b.  Observe that in the definition of
$L_C$, we use the orientation of $N_C$, but we do {\em not} need that
of the surface $C$.  Observe also that the order of $A_1$ and $A_2$ is
irrelevant.

\begin{figure}
\cl{\includegraphics{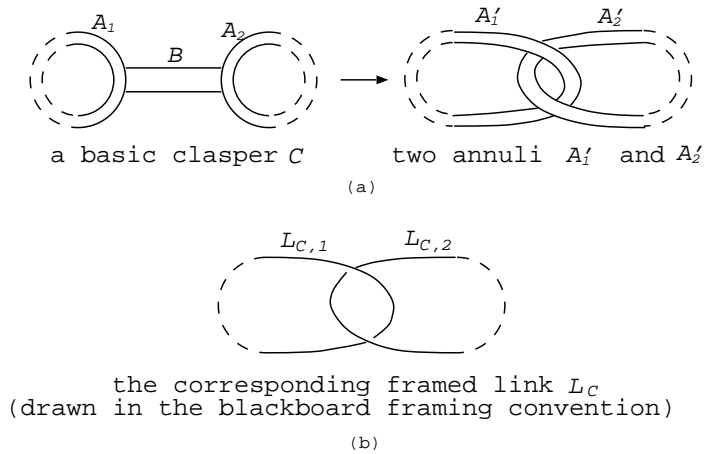}}
\caption{How to associate a framed link to a basic
clasper}
\label{AssociatedFramedLink}
\end{figure}

We define {\em surgery on the basic clasper} $C$ to be surgery on the
associated framed link $L_C$.  The $3$--manifold that we obtain from
$M$ by surgery on $C$ is denoted by $M^C$.  When a small regular
neighborhood $N_C$ of $C$ in $M$ is specified or clear from context,
we may identify $M^C$ with $({M\setminus{\operatorname{int}}}
N_C)\cup_{\partial N_C}{N_C}^C$ (via a diffeomorphism which is
identity outside $N_C$).
\end{definition}

The following Proposition is fundamental in that most of the properties
of claspers that will appear in what follows are derived from it.

\begin{figure}
\cl{\includegraphics{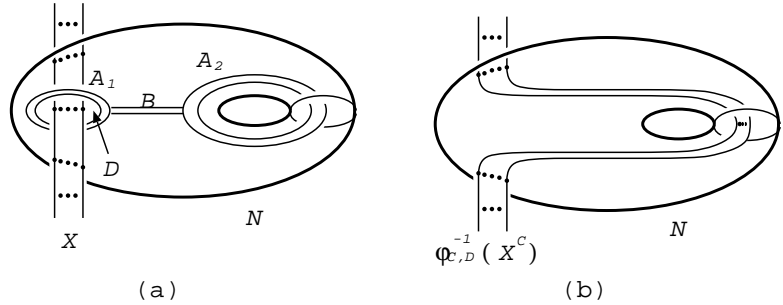}}
\caption{Effect of surgery on a ``disked'' basic clasper}
\label{slide1}
\end{figure}

\begin{proposition}
\label{slide}
{\rm(1)}\qua Let $C=A_1\cup A_2\cup B$ be a basic clasper in a $3$--manifold
$M$, and $D$ a disk embedded in $M$ such that $A_1$ is a collar
neighborhood of $\partial D$ in $D$ and such that $D\cap C=A_1$.  Let
$N$ be a small regular neighborhood of $C\cup D$ in $M$, which is a
solid torus.  Then there is a diffeomorphism $\varphi_{C,D}|_N\colon
N{\overset{\cong}{\longrightarrow}} N^C$ fixing $\partial N=\partial
N^C$ pointwise, which extends to a diffeomorphism $\varphi_{C,D}\colon
M{\overset{\cong}{\longrightarrow}} M^C$ restricting to the identity
on $M\setminus\operatorname{int}N$.

{\rm(2)}\qua In (1) assume that there is a parallel family of ``objects''
(eg, links, claspers, etc)  transversely intersecting the open disk
${D\setminus A_1}$ as depicted in Figure~\ref{slide1}a.  Then the
object $\varphi_{C,D}^{-1}(X^C)$ in $M$ looks as depicted in
Figure~\ref{slide1}b.
\end{proposition}

\begin{figure}
\cl{\includegraphics{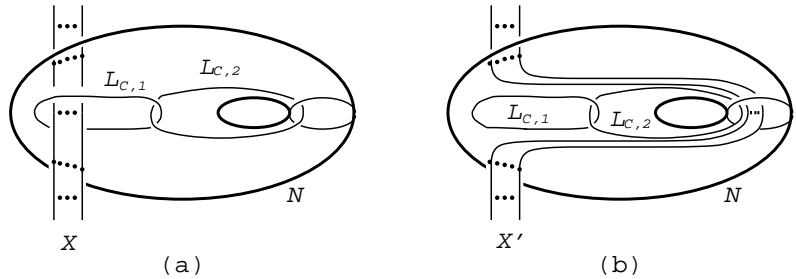}}
\caption{Proof of Lemma~\ref{slide} (2)}
\label{slide2}
\end{figure}

\begin{proof}
(1)\qua Let $L_C=L_{C,1}\cup L_{C,2}\subset N$ be the framed link
associated to $C$.  The component $L_{C,1}$ bounds a disk $D'$ in
${\operatorname{int}} N$ intersecting $L_{C,2}$ transversely once, and $L_{C,1}$ is of
framing zero.  Hence there is a diffeomorphism
$\varphi_{C,D}|_N\colon N\overset{\cong}{\longrightarrow} N^{L_C}(=N^C)$
restricting to the identity on $\partial N$.

(2)\qua The associated framed link $L_C$ looks as depicted in
Figure~\ref{slide2}a.  Before performing surgery on $L_C$, we slide
the object $X$ along the component $L_{C,2}$, obtaining an object $X'$
in $M$ depicted in Figure~\ref{slide2}b.  Since Dehn surgery on $L_C$
in this situation amounts to simply discarding $L_C$ (up to
diffeomorphism), the object $\varphi_{C,D}^{-1}(X^C)$ in $M$ looks as
depicted in Figure~\ref{slide1}b.  
\end{proof}

\begin{remark}
Let $C$ and $D$ be as given in Proposition~\ref{slide}(1).  The
isotopy class of the diffeomorphism $\varphi_{C,D}$ depends not only
on $C$ but also to the disk $D$: If the second homotopy group $\pi_2M$
of $M$ is not trivial, then, for two different bounding disks $D_1$
and $D_2$ for $L_1$ in $M$, the two diffeomorphisms $\varphi_{C,D_1}$
and $\varphi_{C,D_2}$ are not necessarily isotopic to each other.
Thus the data $D$ is necessary in the definition of the diffeomorphism
$\varphi_{C,D}$.  However, if $M$ is a $3$--ball or a $3$--sphere, then
$D$ is unique up to ambient isotopy, and hence $\varphi_{C,D}$ does
not depend on $D$ up to isotopy.
\end{remark}

\begin{figure}
\cl{\includegraphics{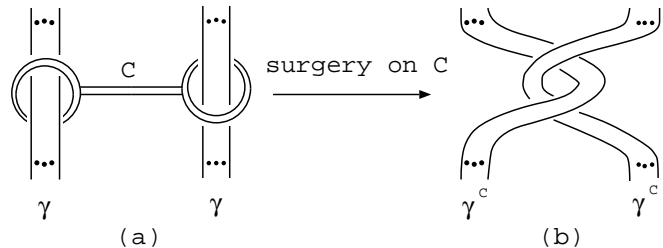}}
\caption{Surgery on a basic clasper clasps two parallel
families of strings}
\label{ClasperClasps}
\end{figure}

\begin{remark}
As a special case of Proposition~\ref{slide}, surgery on a basic
clasper $C$ linking with two parallel families of strings in a link
$\gamma$ as depicted in Figure~\ref{ClasperClasps}a amounts to
producing a ``clasp'' of the two parallel families as depicted in
Figure~\ref{ClasperClasps}b.  This fact explains the name ``clasper.''
\end{remark}

\subsection{Claspers}

\begin{definition}
A {\em clasper} $G=\mathbf{A}\cup\mathbf{B}$ for a link $\gamma$ in a
$3$--manifold $M$ is a non-oriented compact surface embedded in the
interior of $M$ and equipped with a decomposition into two subsurfaces
$\mathbf{A}$ and $\mathbf{B}$.  We call the connected components of
$\mathbf{A}$ the {\em constituents} of $G$, and that of $\mathbf{B}$ the
{\em edges} of $G$.  Each edge of $G$ is a band disjoint from $\gamma$
connecting two distinct constituents, or connecting one constituent
with itself.  An {\em end} of an edge $B$ of $G$ is one of the two
components of $B\cap\mathbf{A}$, which is an arc in $\partial B$.  There
are four kinds of constituents: {\em leaves, disk-leaves, nodes} and
{\em boxes}.  The leaves are annuli, while the disk-leaves, the nodes
and the boxes are disks.  The leaves, the nodes and the boxes are
disjoint from $\gamma $, but the disk-leaves may intersect $\gamma $
transversely.  Also, the constituents must satisfy the following
conditions.
\begin{enumerate}
\item Each node has three incident ends, where it may happen that two
of them are the two ends of one edge.
\item Each leaf (resp.\ disk-leaf) has just one incident end, and hence
has just one incident edge.
\item Each box $R$ of $G$ has three incident ends one of which is
distinguished from the other two.  We call the distinguished incident
end the {\em output end} of $R$, and the other two the {\em input
ends} of $R$.  (In Figures we draw a box $R$ as a rectangle as
depicted in Figure~\ref{BoxAndEdges} to distinguish the output end.)
The edge containing the output (resp.\ an input) end of $R$ is called
the output (resp.\ an input) edge of $R$.
(The two ends of an edge $B$ in a clasper may possibly incident to one
box $R$.  They may be either the two input ends of $R$, or one input
end and the output end of $R$.  In the latter case $B$ is called both
an input edge and the output edge of $R$.)
\end{enumerate}

\begin{figure}
\cl{\includegraphics{BoxAndE.eps}}
\caption{A box}
\label{BoxAndEdges}
\end{figure}

A {\em component} of a clasper $G$ is a connected component of the
underlying surface of $G$ together with the decomposition into
constituents and edges inherited from that of $G$.

Two constituents $P$ and $Q$ of $G$ are said to be {\em adjacent} to
each other if there is an edge $B$ incident both to $P$ and to $Q$.
If this is the case, then we also say that $P$ and $Q$ are {\em
connected} by $B$.

A disk-leaf of a clasper for a link $\gamma$ is called {\em trivial}
if it does not intersect $\gamma$, and {\em simple} if it intersects
$\gamma$ by just one point.

Given a clasper $G$, we obtain a clasper $C_G$ consisting of some
basic claspers in a small regular neighborhood $N_G$ of $G$ in $M$ by
replacing the nodes, the disk-leaves and the boxes of $G$ with some
leaves as illustrated in Figure~\ref{ReplaceNodesDiskleavesAndBoxes}.  The
number of basic claspers contained in $C_G$ is equal to the number of
edges in $G$.  We define {\em surgery on a clasper $G$} to be surgery
on the clasper $C_G$.  More precisely, we define the result $M^G$ from
$M$ of surgery on $G$ by
$$M^G=(M\setminus{\operatorname{int}} N_G)\underset{\partial N_G}{\cup
}{N_G}^{C_G}.$$ So, if a regular neighborhood $N_G$ is explicitly
specified, then we can identify $M\setminus\operatorname{int}N_G$ with
$M^G\setminus\operatorname{int}{N_G}^G$.

\begin{figure}
\cl{\includegraphics{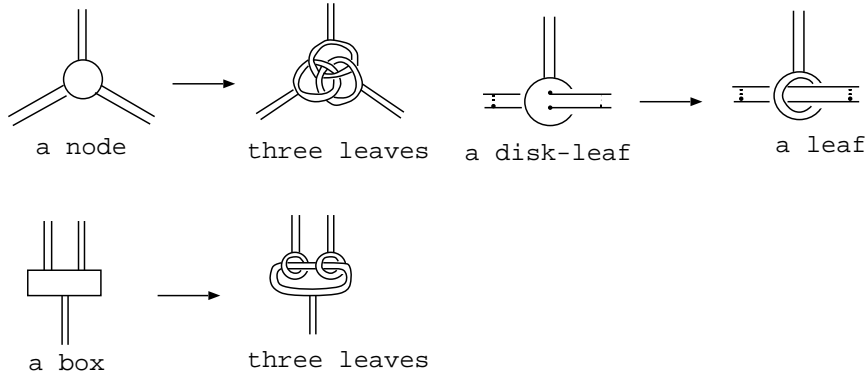}}
\caption{How to replace nodes,
disk-leaves and boxes with leaves}
\label{ReplaceNodesDiskleavesAndBoxes}
\end{figure}

\end{definition}

\begin{figure}
\cl{\includegraphics{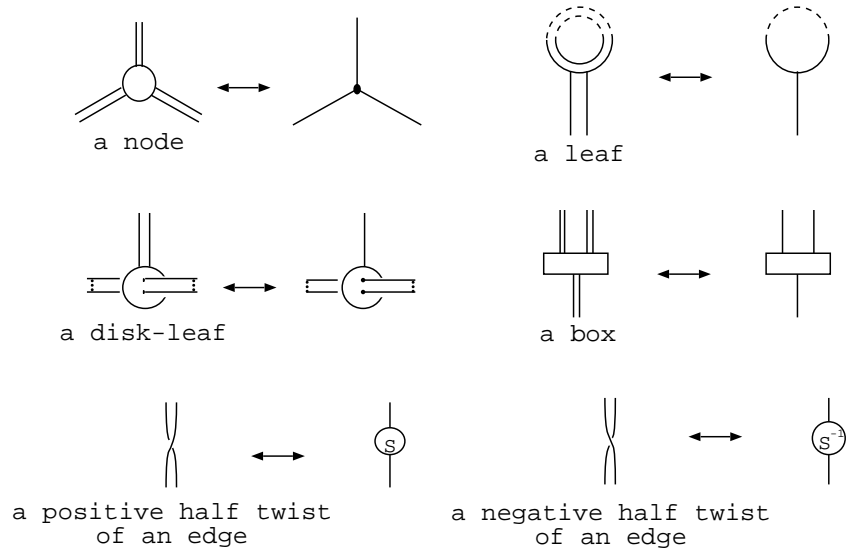}}
\caption{Convention in drawing claspers}
\label{BlackboardFraming}
\end{figure}

\begin{convention}
In Figures we usually draw claspers as illustrated in
Figure~\ref{BlackboardFraming}.  We follow the so-called
blackboard-framing convention to determine the full twists of leaves
and edges.  The last two rules in Figure~\ref{BlackboardFraming} show
how half twists of edges are depicted.
\end{convention}

\subsection{Tame claspers}
Let $V=V_1\cup\dots\cup V_n\quad (n\ge 0)$ be a disjoint union of
handlebodies in the interior of a $3$--manifold $M$, $\gamma$ a link in
$M$ transverse to $\partial V$, and $G\subset {\operatorname{int}} V$
a clasper for $ \gamma$.  We say that $G$ is {\em tame} in $V$ if
there is an orientation-preserving diffeomorphism
$\Phi_{(V,G)}|_V\colon V{\overset{\cong}{\longrightarrow}} V^G$ that
restricts to the identity on $\partial V$.  If this is the case, then
the diffeomorphism $\Phi_{(V,G)}|_V$ extends to the diffeomorphism
$$\Phi_{(V,G)}\colon M{\overset{\cong}{\longrightarrow}}
M^G(=M\setminus{\operatorname{int}} V\cup V^G)$$ restricting to the
identity outside $V$.  Observe that $\Phi_{(V,G)}$ is unique up to
isotopy relative to $M\setminus {\operatorname{int}} V$.  If there is
no fear of confusion, then let $\gamma^{(V,G)}$, or simply $\gamma^G$,
denote the link ${\Phi_{(V,G)}}^{-1}(\gamma^G)$ in $M$, and call it
the result from $\gamma$ of surgery on the pair ${(V,G)}$, or
often simply on $G$.  Observe that surgery on ${(V,G)}$ transforms a
link in $M$ into another link in $M$.

We simply say that $G$ is {\em tame} if $G$ is tame in a regular
neighborhood ${N_G}$ of $G$ in $M$.  If this is the case, we usually
let $\gamma^G$ denote the link $\gamma^{({N_G},G)}$.

If a clasper $G$ is tame in a disjoint union of handlebodies, $V$, and
if ${V'}\subset {\operatorname{int}} M$ is a disjoint union of
embedded handlebodies containing $V$, then $G$ is tame also in ${V'}$,
and the two diffeomorphisms $\Phi_{(V,G)},\Phi_{({V'},G)}\colon
M{\overset{\cong}{\longrightarrow}} M^G$ are isotopic relative to
$M\setminus{\operatorname{int}} {V'}$.  Especially, a tame clasper $G$
is tame in any disjoint union of handlebodies in ${\operatorname{int}}
M$ which contains $G$ in the interior.

\subsection{Some basic properties of claspers}

Let $(\gamma,G)$ and $(\gamma',G')$ be two pairs of links and tame
claspers in a $3$--manifold $M$.  By $(\gamma,G)\sim(\gamma ',G')$, or
simply by $G\sim G'$ if $\gamma=\gamma'$, we mean that the results of
surgeries $\gamma^G$ and ${\gamma'}^{G'}$ are equivalent.

Let ${(\gamma_A,G_A)}$ and ${(\gamma_B,G_B)}$ be two pairs of links
and claspers in $M$ and let $A$ and $B$ be two figures which depicts a
part of ${(\gamma_A,G_A)}$ and a part of ${(\gamma_B,G_B)}$,
respectively.  In such situations we usually assume that the
non-depicted parts of ${(\gamma_A,G_A)}$ and ${(\gamma_B,G_B)}$ are
equal.  We mean by `$A\sim B$' in figures that
${(\gamma_A,G_A)}\sim{(\gamma_B,G_B)}$.

\begin{proposition}
\label{moves}
Let $(\gamma,G)$ and $(\gamma',G')$ be two pairs of links and claspers
in $M$.  Suppose that $V$ is a union of handlebodies in $M$ in which
$G$ and $G'$ are tame.  Suppose that $(\gamma,G)$ and $(\gamma',G')$
are related by one of the moves 1--12 performed in $V$.  Then the
results of surgery $(\gamma\cap V)^G, ({\gamma'}\cap V)^{G'}$ are
equivalent in $V$, and hence $\gamma^G$ and ${\gamma'}^{G'}$ are
equivalent in $M$.
\end{proposition}

\begin{figure}
\cl{\includegraphics[width=4.6in]{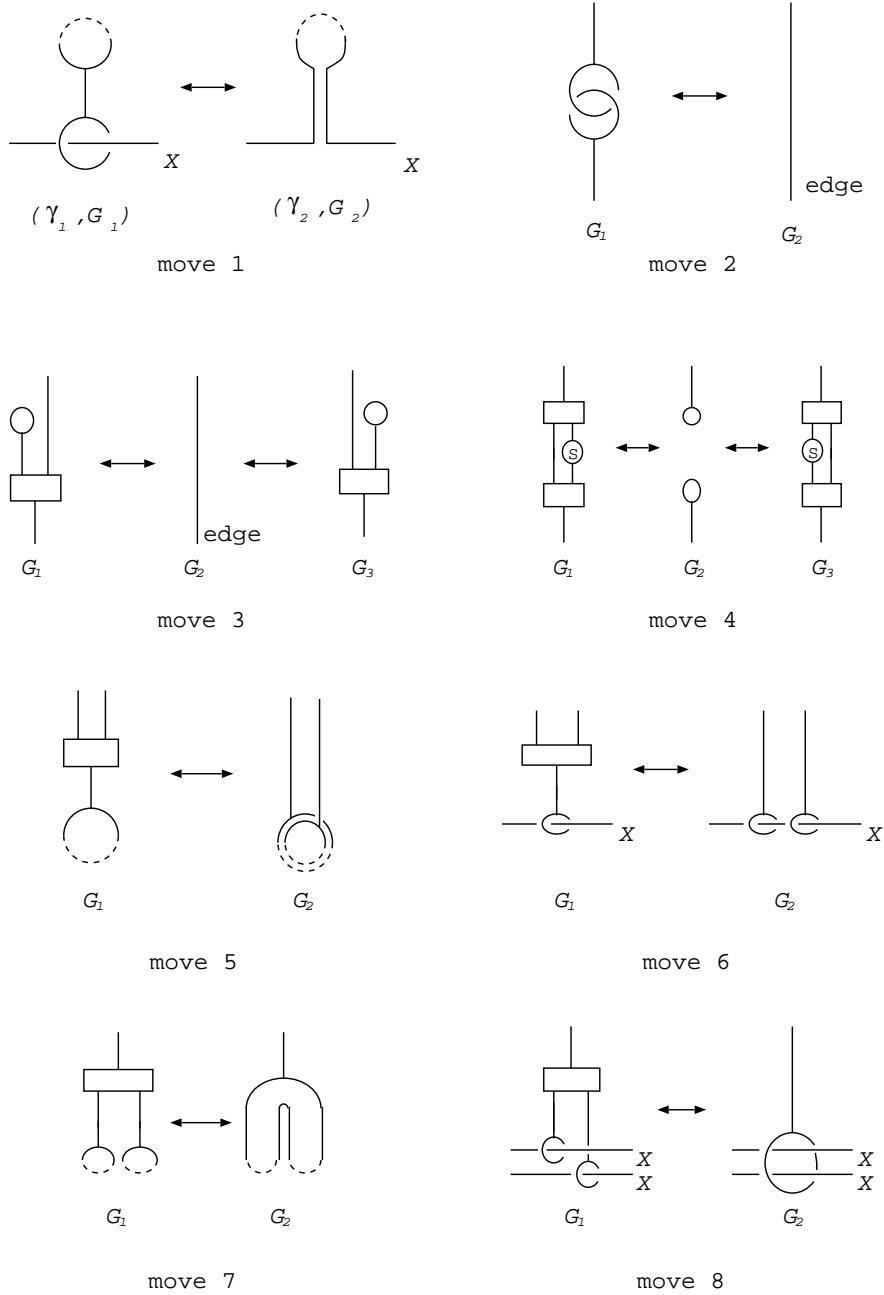}}
\caption{{Moves on claspers and links which do not change the
result of surgeries.  Here $\text{---}_X$ represents a
parallel family of strings of a link and/or edges and leaves of
claspers.}}
\label{moves1}
\end{figure}

\begin{figure}[ht!]
\cl{\includegraphics{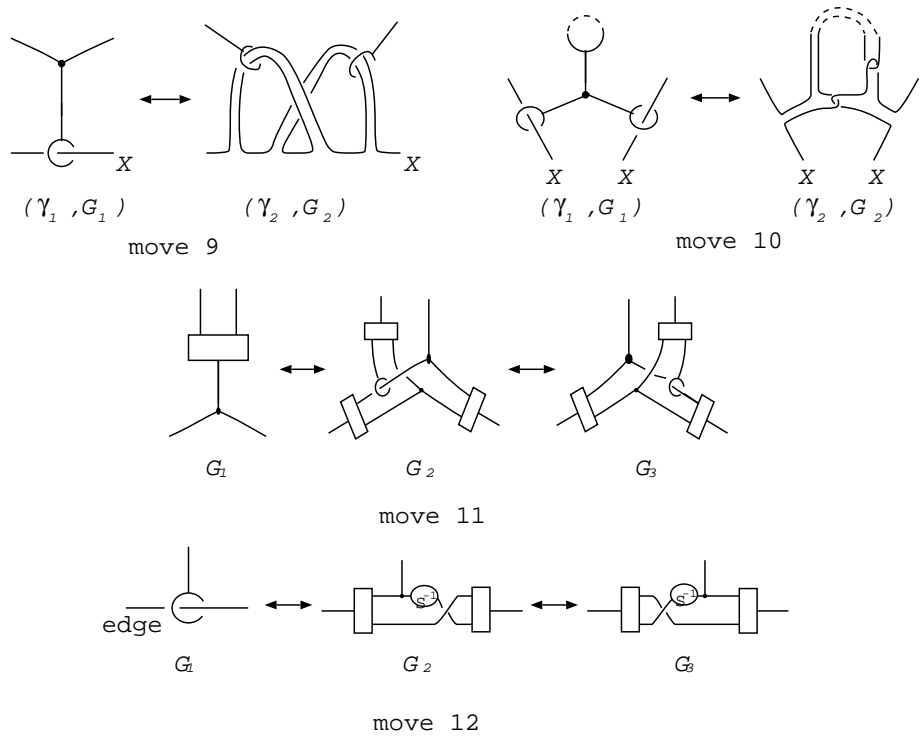}}
\caption{(Continued)}
\label{moves2}
\end{figure}

\begin{figure}
\cl{\includegraphics{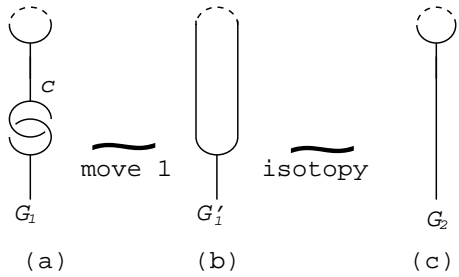}}
\caption{Proof of move~{2}}
\label{pfmove2}
\end{figure}

\begin{figure}
\cl{\includegraphics{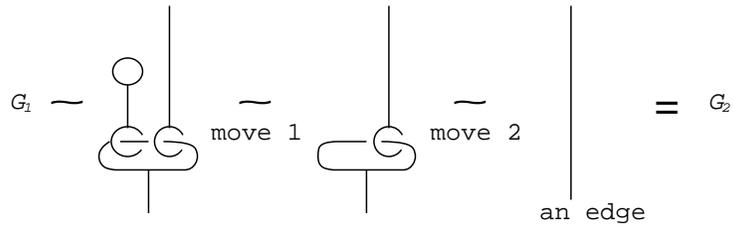}}
\caption{Proof of move~{3}}
\label{pfmove3}
\end{figure}

\begin{figure}[ht!]
\cl{\includegraphics{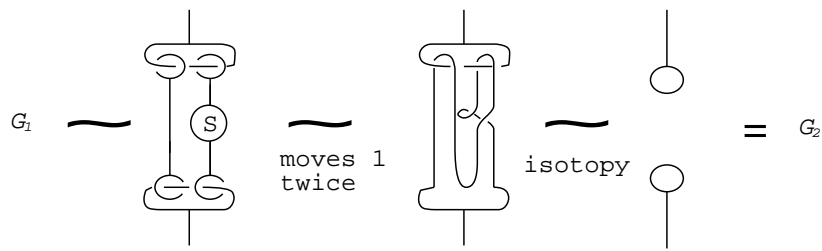}}
\caption{Proof of move~{4}}
\label{pfmove4}
\end{figure}

\begin{figure}[ht!]
\cl{\includegraphics{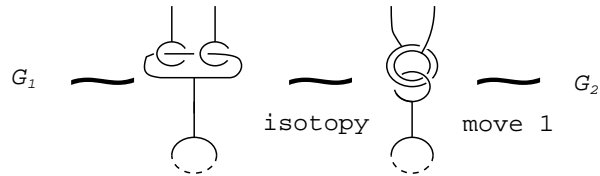}}
\caption{Proof of move~{5}}
\label{pfmove5}
\end{figure}

\begin{figure}[ht!]
\cl{\includegraphics{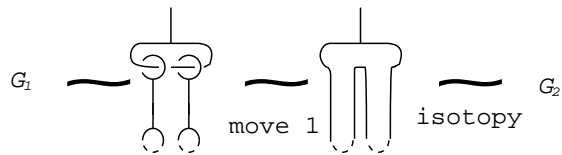}}
\caption{Proof of move~{7}}
\label{pfmove7}
\end{figure}

\begin{figure}[ht!]
\cl{\includegraphics{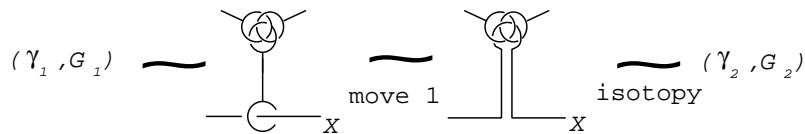}}
\caption{Proof of move~{9}}
\label{pfmove9}
\end{figure}

\begin{figure}[ht!]
\cl{\includegraphics{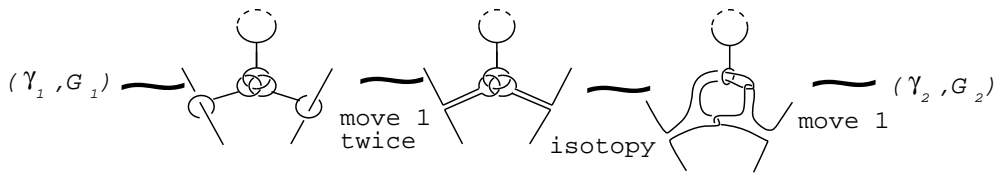}}
\caption{Proof of move~{10}}
\label{pfmove10}
\end{figure}

\begin{figure}[ht!]
\cl{\includegraphics{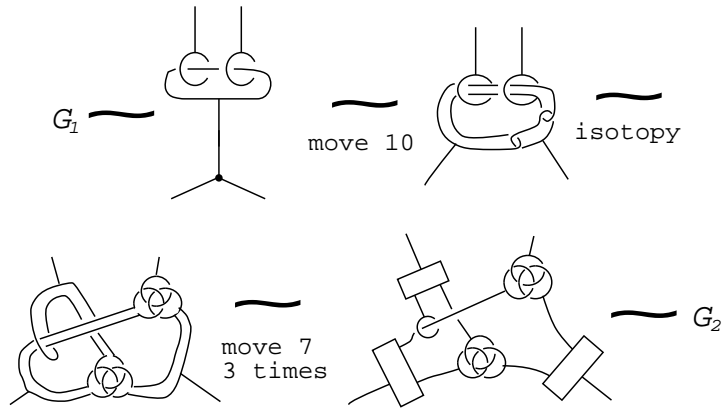}}
\caption{Proof of move~{11}}
\label{pfmove11}
\end{figure}

\begin{figure}[ht!]
\cl{\includegraphics{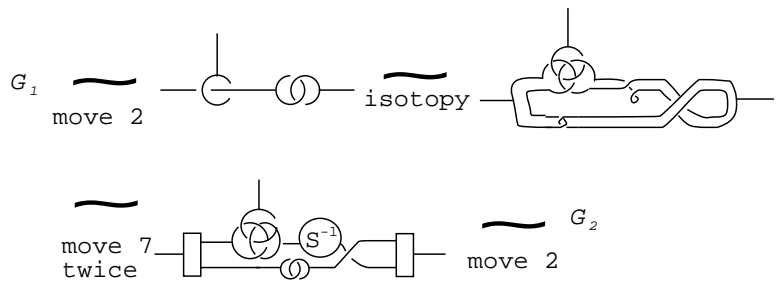}}
\caption{Proof of move~{12}}
\label{pfmove12}
\end{figure}

\begin{proof}
In this proof, let $(\gamma_i,G_i)\quad (i=1,2,3)$ denote the pair
of the link and the clasper depicted in the $i$th term in
each row in Figures \ref{moves1} and \ref{moves2}.

{\bf Move 1}\qua This is just Proposition~\ref{slide}.

{\bf Move 2}\qua We may assume that the edge depicted in the right
side and hence one of the edges in the left side are incident to
leaves since, if not, we can replace the incident constituent of the
edges with some leaves without changing the results of surgeries.
Thus we may assume that the clasper on the left side is as depicted in
Figure~\ref{pfmove2}a.  Surgery on the basic clasper $C$ yields a
clasper $G_1'$ depicted in Figure~\ref{pfmove2}b, which is ambient
isotopic to $G_2$ depicted in Figure~\ref{pfmove2}c.  Hence we have
$G_1\sim G_2$.

{\bf Move 3}\qua Figure~\ref{pfmove3} implies $G_1\sim G_2$.  The
proof of $G_3\sim G_2$ is similar.

{\bf Move 4}\qua See Figure~\ref{pfmove4}.

{\bf Move 5}\qua See Figure~\ref{pfmove5}.

{\bf Move 6}\qua Use move~{5}.

{\bf Move 7}\qua See Figure~\ref{pfmove7}.

{\bf Move 8}\qua Use move~{7}.

{\bf Move 9}\qua See Figure~\ref{pfmove9}.

{\bf Move 10}\qua See Figure~\ref{pfmove10}.

{\bf Move 11}\qua For $G_1\sim G_2$, see Figure~\ref{pfmove11}.  The
proof of $G_1\sim G_3$ is similar.

{\bf Move 12}\qua For $G_1\sim G_2$, see Figure~\ref{pfmove12}.  The
proof of $G_1\sim G_3$ is similar.
\end{proof}

\begin{remark}
Proposition~\ref{moves} can be modified as follows.  If two pairs
$(\gamma,G)$ and $(\gamma',G')$ are pairs of links and claspers in $M$
with $G$ and $G'$ not necessarily tame, and if they are related by one
of the moves in Proposition~\ref{moves} then the results of surgeries
$(M,\gamma)^G, (M,\gamma')^{G'}$ are related by a diffeomorphism
restricting to the identity on boundary.  This fact will not be used
in this paper but in future papers in which we will prove the results
announced in Section~8.
\end{remark}

\begin{remark}
In Figures \ref{moves1} and \ref{moves2}, there are no disk-leaves
depicted.  However, we often use these moves on claspers with
disk-leaves by freely replacing disk-leaves with leaves and vice
versa.
\end{remark}

\section{Tree claspers and the $C_k$--equivalence relations on links}
\label{section3}

\subsection{Definition of tree claspers}

\begin{figure}
\cl{\includegraphics{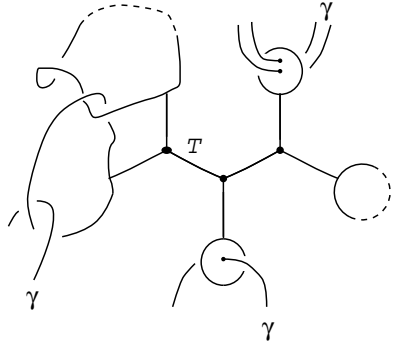}}
\caption{An example of a tree clasper $T$ for a link
$\gamma $ in a $3$--manifold $M$.  Leaves of $T$ may link with other
leaves, and may run through any part of the manifold $M$.}
\label{ExTreeClasper}
\end{figure}

\begin{definition}
A {\em tree clasper} $T$ for a link $\gamma$ in a $3$--manifold $M$ is
a connected clasper without box such that the union of the nodes
and the edges of $T$ is simply connected, and is hence ``tree-shaped.''
Figure~\ref{ExTreeClasper} shows an example of a tree clasper for a
link $\gamma$.

A tree clasper $T$ is {\em admissible} if $T$ has at least one
disk-leaf, and is {\em strict\/} if (moreover) $T$ has no leaves.
Observe that the underlying surface of a strict tree clasper is
diffeomorphic to the disk $D^2$.  A strict tree clasper $T$ is {\em
simple\/} if every disk-leaf of $T$ is simple.
\end{definition}

\begin{definition}
A {\em forest clasper} $T={T_1\cup\dots\cup T_p}$ ($p\ge 0$) for a
link $\gamma$ is a clasper $T$ consisting of $p$ tree claspers
$T_1,\dots,T_p$ for $\gamma$.  The forest clasper $T$ is {\em
admissible}, (resp.\ {\em strict}, {\em simple}) if every component of
$T$ is admissible (resp.\ strict, simple).
\end{definition}

\begin{proposition}
\label{admtc}
Every admissible tree clasper for a link in a $3$--manifold $M$ is
tame.  Especially, every strict tree clasper is tame.
\end{proposition}

\begin{proof}
Let $T$ be an admissible tree clasper for a link $\gamma$ in $M$,
${N_T}\subset {\operatorname{int}} M$ a small regular neighborhood of
$T$ in $M$, and $D$ a disk-leaf of $T$.  If there are other
disk-leaves of $T$, then we may safely replace them with leaves since
the tameness in $N_T$ of the new $T$ will imply that of the old $T$.
Assume  that $D$ is the only disk-leaf in $T$.  If $T$ has no node,
then $D$ is adjacent to a leaf $A$, and $T$ is tame in ${N_T}$ by
Proposition~\ref{slide}.  Hence we may assume that $T$ has at least
one node, and that the proposition holds for admissible tree claspers
which have less nodes than $T$ has.  Applying move~{9} to $D$ and the
adjacent node, we obtain two disjoint admissible tree claspers $T_1$
and $T_2$ in ${N_T}$ for $\gamma'\cong\gamma$ such that there is a
diffeomorphism
${{N_T}}^T{\overset{\cong}{\longrightarrow}}{{N_T}}^{T_1\cup T_2}$
fixing $\partial {N_T}$ pointwise.  Since $T_1$ and $T_2$ are tame,
there is a diffeomorphism
${{N_T}}^T{\overset{\cong}{\longrightarrow}}{{N_T}}^{T_1\cup
T_2}{\overset{\cong}{\longrightarrow}}{N_T}$ fixing $\partial {N_T}$
pointwise.  Hence $T$ is tame.  
\end{proof}

By Proposition~\ref{admtc}, an admissible tree clasper $T$ for a link $\gamma$ in
a $3$--manifold $M$ determines a link $\gamma^T$ in $M$.  Hence we may
think of surgery on an admissible tree clasper as an operation on
links in a {\em fixed} $3$--manifold $M$.

\begin{proposition}
\label{trivialleaf}
Let $T$ be an admissible tree clasper for a link $\gamma$ in $M$ with
at least one trivial disk-leaf.  Then $\gamma^T$ is equivalent to
$\gamma$.
\end{proposition}

\begin{proof}
There is a sequence of admissible forest claspers for $\gamma$,
$G_0=T,G_1,\dots ,$ $ G_p=\emptyset$ ($p\ge 0$) from $T$ to $\emptyset$
such that, for each $i=0,\dots ,p-1$, $G_{i+1}$ is obtained from $G_i$
by move~{1} or by move~{9}, where the ``object to be slided'' is empty.
Hence we have $\gamma^T=\gamma^{G_0}\cong\gamma^{G_p}=\gamma
^\emptyset=\gamma$.  
\end{proof}

\subsection{$C_k$--moves and $C_k$--equivalence}

\begin{definition}
The {\em degree}, $\deg T$, of a strict tree clasper $T$ for a link
$\gamma$ is the number of nodes of $T$ plus $1$.  The degree of a
strict forest clasper is the minimum of the degrees of its component
strict tree claspers.
\end{definition}

\begin{definition}
Let $M$ be a $3$--manifold and let $k\ge 1$ be an integer.  A {\em
(simple) ${C_k}$--move} on a link $\gamma$ in $M$ is a surgery on a
(simple) strict tree clasper of degree $k$.  More precisely, we say
that two links $\gamma$ and $\gamma'$ in $M$ are related by a (simple)
${C_k}$--move if there is a (simple) strict tree clasper $T$ for $\gamma$
of degree $k$ such that $\gamma^T$ is equivalent to $ \gamma'$.  We
write $\gamma{\underset{C_k}{\longrightarrow}}\gamma'$ ($\gamma {\underset{{sC_k}}{\longrightarrow}}\gamma'$) to mean that two
links $\gamma$ and $\gamma'$ are related by a (simple) ${C_k}$--move.

The {\em ${C_k}$--equivalence} (resp.\ {\em ${sC_k}$--equivalence}) is the
equivalence relation on links generated by the ${C_k}$--moves
(resp.\ simple ${C_k}$--moves) and ambient isotopies.  By $\gamma
{\underset{C_{k}}{\sim}}\gamma'$ (resp.\ $\gamma{\underset{{sC_k}}{\sim}}\gamma'$) we mean that $\gamma$
and $\gamma'$ are ${C_k}$--equivalent (resp.\ ${sC_k}$--equivalent).
\end{definition}

The following result means that the ${C_k}$--equivalence relation becomes
finer as $k$ increases.

\begin{proposition}
\label{leqkeq}
If $1\le k\le l$, then a $C_l$--move is achieved by a $C_k$--move, and
hence $C_l$--equivalence implies $C_k$--equivalence.
\end{proposition}

\begin{proof}
It suffices to show that, for each $k\ge1$ and for  a strict tree
clasper $T$ of degree $k+1$ for a link $\gamma$ in a $3$--manifold $M$,
there is a strict tree clasper $T'$ of degree $k$ for $\gamma$
such that $\gamma ^T\cong\gamma^{T'}$.  We choose a node $V$ of $T$
which is adjacent to at least two disk-leaves $D_1$ and $D_2$; see
Figure~\ref{leqkeq1}a.  
Applying move~{2} to the edge $B$ of $T$ that is incident to $V$
but neither to $D_1$ nor $D_2$, we obtain a clasper $T_1\cup T_2$ which is
tame in a small regular neighborhood $N_T$ of $T$ in $M$ consisting of
two admissible tree claspers $T_1$ and $T_2$ such that
$\gamma^{T_1\cup T_2}\cong \gamma^T$, see Figure~\ref{leqkeq1}b.  Here $T_1$
contains the node $V$ and the two disk-leaves $D_1$ and $D_2$.  By
move~{10} we obtain a link $\gamma^{T_2}$ such that $\gamma^{T_1\cup
T_2}=(\gamma ^{T_1})^{T_2}$, see Figure~\ref{leqkeq1}c.  Regarding the leaf
$A_2$ as a disk-leaf in the obvious way, we obtain a strict tree
clasper $T_2$ for $\gamma^{T_1}$ of degree $k$.  Observe that
$\gamma^{T_1}$ is equivalent to $\gamma$ and that $
(\gamma^{T_1})^{T_2}\cong \gamma ^T$.  Therefore there is a strict
tree clasper $T'$ for $\gamma$ of degree $k$ such that
$\gamma^{T'}\cong\gamma^T$.  
\end{proof}

\begin{figure}
\cl{\includegraphics{leqkeq1.eps}}
\nocolon\caption{}
\label{leqkeq1}
\end{figure}

\begin{definition}
Two links in $M$ are said to be {\em $C_\infty$--equivalent} if they
are ${C_k}$--equivalent for all $k\ge 1$.
\end{definition}

\begin{conjecture}
\label{conjCinfty}
Two links in a $3$--manifold $M$ are equivalent if
and only if they are $C_\infty$--equivalent.
\end{conjecture}

\subsection{Zip construction}
Here we give a technical construction which we call a {\em zip
construction} and which is crucial in what follows.

\begin{definition}
A {\em subtree} $T$ in a clasper $G$ is a union of some leaves,
disk-leaves, nodes and edges of $G$ such that
\begin{enumerate}
\item the total space of $T$ is connected,
\item $T\setminus(\text{leaves of $T$})$ is simply connected,
\item $T\cap\overline{C\setminus T}$ consists of  ends of some edges in $T$.

\end{enumerate}
We call each connected component of the intersection of $T$ and the
closure of $G\setminus T$ an {\em end of $T$}, and the edge containing
it an {\em end-edge} of $T$.  A subtree is said to be {\em strict} if
$T$ has no leaves.

An {\em output subtree} $T$ in $G$ is a subtree of $G$ with just one
end that is an output end of a box.
\end{definition}

\begin{definition}
A {\em marking} on a clasper $G$ is a set $\bfM$ of input ends of
boxes such that for each box $R$ of $G$, at most one input end of $R$
is an element of $\bfM$ and such that for each $e\in\bfM$, the box
$R\supset e$ is incident to an output subtree.
\end{definition}

\begin{definition}
Let $G$ be a clasper for a link $\gamma$ in $M$, and $\bfM$ a marking
on $G$.  A {\em zip construction} ${\operatorname{Zip}}(G,\bfM)$ is a
clasper for $\gamma$ contained in a small regular neighborhood $N_G$
of $G$  constructed as follows.  If $\bfM$ is empty, then we set
${\operatorname{Zip}}(G,\emptyset)=G$.  Otherwise we define
${\operatorname{Zip}}(G,\bfM)$ to be a clasper for $\gamma$ contained
in $N_G$ obtained from $(G,\bfM)$ by iterating the operations of the
following kind until the marking $\bfM$ becomes empty.
\begin{itemize}
\item We choose an element $e\in\bfM$ and let $R$ be the box
containing $e$, $T$ the output subtree, and $B$ the end-edge of $T$.
Let $G'$ be the clasper obtained from $G$ by applying move {5}, {6} or
{11} to $R$ according as the constituent incident to $B$ at the
opposite side of $R$ is a leaf, a disk-leaf or a node, respectively.
In the first two cases we set $\bfM'=\bfM\setminus\{e\}$, and in the
last case we set $\bfM'=(\bfM \setminus\{e\})\cup\{e_1,e_2\}$, where
$e_1$ and $e_2$ are ends in $G'$ determined as in Figure~\ref{e1e2}.
Then let $G'$ be the new $G$ and $\bfM'$ the new $\bfM$.
\end{itemize}
This procedure clearly terminates, and the result
${\operatorname{Zip}}(G,\bfM)$ does not depend on the choice of $e$ in
each step.  Observe that if there are more than one element in $\bfM$,
then ${\operatorname{Zip}}(G,\bfM)$ is obtained from $G$ by separately
applying the above construction to each element of $\bfM$; eg,
${\operatorname{Zip}}(G,\{e,e'\})={\operatorname{Zip}}({\operatorname{Zip}}(G,\{e\}),\{e'\})$.

\begin{figure}
\cl{\includegraphics*{e1e2.eps}}
\nocolon\caption{}
\label{e1e2}
\end{figure}

The clasper ${\operatorname{Zip}}(G,\bfM)$ is unique up to isotopy in
$N_G$.  We call it the {\em zip construction} for $(G,\bfM)$.  By
construction, $G$ and ${\operatorname{Zip}}(G,\bfM)$ have
diffeomorphic results of surgeries.  Hence, if $G$ is tame, then
$\operatorname{Zip}(G,\bfM)$ is tame in $N_G$ and that the results of
surgeries on $G$ and $\operatorname{Zip}(G,\bfM)$ are equivalent.

If $\bfM$ is a singleton set $\{e\}$, then we set
${\operatorname{Zip}}(G,e)={\operatorname{Zip}}(G,\{e\})$ and call it
the zip construction for $(G,e)$.
\end{definition}

Figure~\ref{zip} shows an example of zip construction.  The name ``zip
construction'' comes from the fact that the procedure of obtaining a
zip construction looks like ``opening a zip-fastener.''

\begin{figure}
\cl{\includegraphics*{zip.eps}}
\nocolon\caption{}
\label{zip}
\end{figure}

\begin{definition}
An {\em input subtree} $T$ of $G$ is a subtree of $G$ each of whose
ends is an input end of a box.  An input subtree $T$ is said to be
{\em good} if the following conditions hold.
\begin{enumerate}
\item $T$ is strict.
\item The ends of $T$ form a marking of $G$.  
\item For each box $R $ incident to $T$, the output subtree of $R$ is
strict. 
\end{enumerate}
Each strict output subtree in the condition  3 above is said to be
{\em adjacent} to $T$.  
\end{definition}

\begin{definition}
The {\em degree} of a strict subtree $T$ of a clasper $G$ is half the
number of disk-leaves and nodes, which is a half-integer.  The {\em
$e$--degree} (`e' for `essential') of a good input subtree  $T$ of $G$ is
defined to be the sum ${\operatorname{deg}} T + {\operatorname{deg}}
T_1 +\cdots + {\operatorname{deg}} T_m$, where $T_1,\dots,T_m$ ($m\ge0$) 
are the adjacent strict output subtrees of $T$.  The $e$--degree is
always a positive integer.  We say that $T$ is {\em $e$--simple} if $T$ 
and the $T_1,\dots,T_m$ are all simple.
\end{definition}

\begin{definition}
Let $G$ be a clasper and let $X$ be a union of constituents and edges
of $G$.  Assume that the incident edges of the leaves, disk-leaves and
nodes in $X$ are in $X$, that the incident constituents of the edges
are in $X$, and that for each box $R$ in $X$, the output edge of $R$
is in $X$ and at least one of the input edges is in $X$.  Thus $X$ may 
fail to be a clasper only at some {\em one-input boxes}, see
Figure~\ref{smoothing}a.  Let $X\thinspace\tilde{}$ denote the clasper 
obtained from $X$ by ``smoothing'' the one-input boxes, see
Figure~\ref{smoothing}b.  We call $X\thinspace\tilde{}$ the {\em
smoothing} of $X$.

Let $Y$ be a union of constituents and edges of a clasper $G$ such
that the closure of $G\setminus Y$ can be smoothed as above.  Then the 
smoothing $(\overline{G\setminus Y})\thinspace\tilde{}$ is denoted by $G\ominus Y$.
\end{definition}

\begin{figure}
\cl{\includegraphics{smoothin.eps}}
\nocolon\caption{}
\label{smoothing}
\end{figure}

\begin{lemma}
\label{zipnew}
Let $G$ be a tame clasper for a link $\gamma$ in a $3$--manifold $M$,
and $T$ a good input subtree of $G$ of $e$--degree $k\ge1$.  Then
$\gamma^G$ is obtained from $\gamma^{G\ominus T}$ by a $C_k$--move.
If, moreover, $T$ is $e$--simple, then $\gamma^G$ is obtained from
$\gamma^{G\ominus T}$ by a simple $C_k$--move.
\end{lemma}

\begin{figure}
\cl{\includegraphics{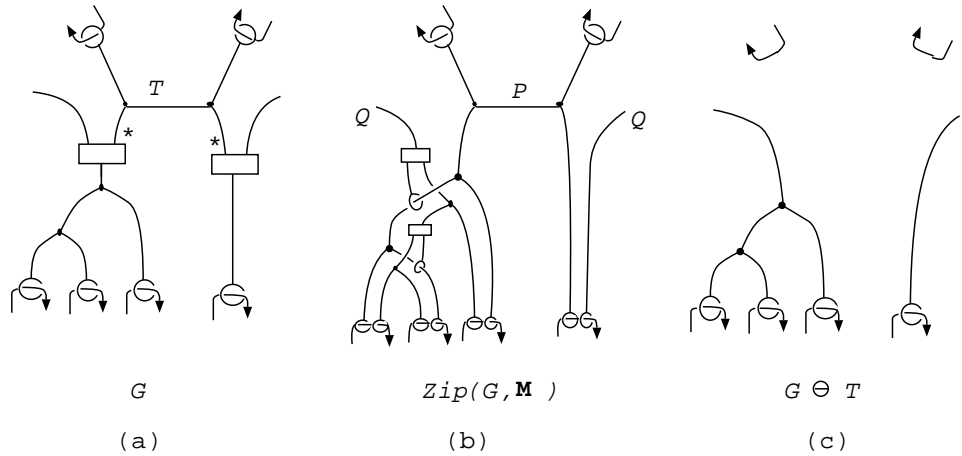}}
\caption{In (a), two asterisks are placed near the two ends of the 
good input subtree $T$ in $G$.}
\label{zipnew1}
\end{figure}

\begin{proof}
Let $\bfM$ denote the set of ends of $T$.  Then
$\operatorname{Zip}(G,\bfM)$ is a disjoint union of a strict tree
clasper $P$ of degree $k$ and a clasper $Q$, see Figure~\ref{zipnew1}a
and b.  We have $\gamma^Q\cong\gamma^{G\ominus T}$, see
Figure~\ref{zipnew1}c.  Hence $\gamma^{G\ominus
T}\cong\gamma^Q\overset{P}{\underset{C_k}{\longrightarrow}}\gamma^{P\cup
Q}\cong\gamma^T$.

If $T$ and the output trees adjacent to
$T$ are simple, then so is $P$.  Hence $\gamma^T$ is obtained from
$\gamma^{G\ominus T}$ by one simple $C_k$--move.
\end{proof}

\subsection{$C_k$--equivalence and simultaneous application of $C_k$--moves}

The rest of this section is devoted to proving the following theorem.

\begin{theorem}
\label{equivalence}
Let $\gamma$ and $\gamma'$ be two links in a $3$--manifold $M$ and let
$k\ge 1$ be an integer.  Then the following conditions are equivalent.
\begin{enumerate}
\item $\gamma$ and $\gamma'$ are ${C_k}$--equivalent.
\item $\gamma$ and $\gamma'$ are  ${sC_k}$--equivalent.
\item $\gamma'$ is obtained from $\gamma$ by surgery on a strict
forest clasper $T={T_1\cup\dots\cup T_l}$ ($l\ge0$) consisting of
strict tree claspers $T_1,\dots,T_l$ of degree $k$.
\item $\gamma'$ is obtained from $\gamma$ by surgery on a simple strict
forest clasper $T={T_1\cup\dots\cup T_l}$ ($l\ge0$) consisting of
simple strict tree claspers $T_1,\dots,T_l$ of degree $k$.
\end{enumerate}	
\end{theorem}

\begin{remark}
By Proposition~\ref{leqkeq}, we may allow in the conditions 3 and 4 above
(simple) strict forest claspers of degree $k$ possibly containing
components of degree $\ge k$.
\end{remark}		

\begin{proof}[Proof of 2$\Rightarrow$1, 4$\Rightarrow$3, 3$\Rightarrow$1
and 4$\Rightarrow$2 of Theorem~\ref{equivalence}] The implications
2$\Rightarrow$1 and 4$\Rightarrow$3 are clear.  The implications
3$\Rightarrow$1 and 4$\Rightarrow$2 come from the following
observation: If $T=T_1\cup\dots\cup T_l$ ($l\ge 0$) is a (simple) strict
forest clasper for $\gamma$ of degree $k$, then there is a sequence of
(simple) ${C_k}$--moves
$$\gamma{\overset{T_1}{\underset{(s){C_k}}{\longrightarrow}}}
\gamma^{T_1}{\overset{T_2}{\underset{(s){C_k}}{\longrightarrow}}}
\gamma^{T_1\cup T_2}{\overset{T_3}{\underset{(s){C_k}}
{\longrightarrow}}}\dots {\overset{T_l}{\underset{(s){C_k}}
{\longrightarrow}}}\gamma^{{T_1\cup\dots\cup T_l}}$$ from $\gamma $ to
$\gamma^{{T_1\cup\dots\cup T_l}}$.
\end{proof}

In the following we first prove 1$\Rightarrow$2 by showing that a
${C_k}$--move can be achieved by a finite sequence of simple
${C_k}$--moves, and then prove 2$\Rightarrow$4 by showing that a
sequence of simple ${C_k}$--moves and inverses of simple ${C_k}$--moves
can be achieved by a surgery on a simple strict forest clasper of
degree $k$.

\begin{proof}[Proof of 1$\Rightarrow$2 of Theorem~\ref{equivalence}]
It suffices to prove the following claim.

\begin{claim}
If a link $\gamma'$ is obtained from a link $\gamma $ by surgery on a
strict tree clasper $T$ for $\gamma$ of degree $k$, then there is a
sequence of simple ${C_k}$--moves from $\gamma$ to $\gamma'$.
\end{claim}

Before proving the claim, we make some definitions which is used only
in this proof and the next remark: For a disk-leaf $D$ in a strict
tree clasper $T$ for a link $\gamma$, let $n(D)$ denote the number of
intersection points of $D$ with $\gamma$.  We also set
$n(T)=\prod_{D}n(D)$, where $D$ runs over all disk-leaves of $T$.

The proof of the claim is by induction on $n=n(T)$.  If $n=0$, then
$\gamma'$ is equivalent to $\gamma$ by Proposition~\ref{trivialleaf}.
If $n=1$, then $T$ is simple, and therefore $\gamma$ and $\gamma'$ are
related by one simple ${C_k}$--move.  Let $n\ge 2$ and suppose that the
claim holds for strict tree claspers with smaller $n$.  Then there is
at least one disk-leaf $D$ of $T$ with $n(D)\ge 2$.  Applying move~{8}
to $D$, we obtain a clasper $G_1$ which is tame in a small regular
neighborhood $N_T$ of $T$ in $M$ consisting of a box $R$, a strict
output subtree
$T'$, two input edges $B_1$ and $B_2$ of $R$, and two disk-leaves
$D_1$ and $D_2$ incident to $B_1$ and $B_2$, respectively.  Here we
have $n(D_1)=n(D)-1$ and $n(D_2)=1$, see Figure~\ref{1to2-1}a and b.
The union $D_1\cup B_1$ is a good input subtree of $e$--degree $k$.  We
consider the zip construction $\operatorname{Zip}(G_1,\{B_1\cap
R\})=P\cup Q$, where $P$ is a strict tree clasper of degree $k$ with
$n(P)=(n(D)-1)n(T)/n(D)<n(T)$.  By the induction hypothesis, there is
a sequence of simple $C_k$--moves from $\gamma^Q$ to $\gamma^{P\cup
Q}\cong\gamma^{G_1}\cong\gamma^T$.  We have $\gamma^Q\cong\gamma^{Q'}$
by move~{3}, where $Q'=G_1\ominus (D_1\cup B_1)$ is a strict tree
clasper of degree $k$ with $n(Q')=n(P)/n(D)<n(P)$.  By the induction
hypothesis, there is a sequence of simple $C_k$--moves from $\gamma$ to
$\gamma^{Q'}\cong\gamma^Q$. This completes the proof of the claim and
hence that of 1$\Rightarrow$2.
\end{proof}

\begin{figure}
\cl{\includegraphics{1to2-1.eps}}
\nocolon\caption{}
\label{1to2-1}
\end{figure}

\begin{remark}
It is clear from the above proof that surgery on a strict tree
clasper $T$ of degree $k$ is achieved by a sequence of $n(T)$ 
simple ${C_k}$--moves.
\end{remark}

Before proving 2$\Rightarrow$4 of Theorem~\ref{equivalence}, we need some
definitions and lemmas.

In the following, a {\em tangle} $\gamma$ will mean a link in a
$3$--ball $B^3$ consisting of only some arcs.  A tangle $\gamma$ is
called {\em trivial} if the pair $(B^3,\gamma)$ is diffeomorphic to
the pair $(D^2\times [-1,1],{\gamma_0})$ with ${\gamma_0}\subset
D^2\times \{0\}$ after smoothing the corners.

For later convenience, the following lemma is stated more strongly
than actually needed here.

\begin{lemma}
\label{trivialtangle}
Let $\gamma$ be a trivial tangle in $B^3$, and let $T$ be a simple
strict tree clasper for $\gamma$ of degree $k\ge 1$.  Suppose that
there is a properly embedded disk $D\subset B^3$ such that $T\subset
D$ and such that each component of $\gamma$ transversely intersects
$D$ at a point in a disk-leaf of $T$, see Figure~\ref{trivialtangle1}a for
example.  Then the tangle $\gamma^T$ is trivial.  Moreover, $\gamma
^T$ is of the form depicted in Figure~\ref{trivialtangle1}b, where $\beta $
is a pure braid of $2k+2$ strands such that
\begin{enumerate}
\item $\beta $ is contained in the $k$th lower central series subgroup
$P(2k+2)_k$ of the pure braid group $P(2k+2)$,
\item for each $i=1,\dots ,k+1$, the result from $\beta $ of removing
the $(2i-1)$st and the $2i$th strands is a trivial pure braid of $2k$
strands, where we number the strings from left to right,
\item the first strand of $\beta$ is trivial and not linked with each
others, ie, $\beta$ has a projection with no crossings on the first
strand  (by the condition 2, $\beta\setminus(\text{the 2nd strand})$
is trivial).
\end{enumerate}
\end{lemma}

\begin{proof}
The proof is by induction on $k$.  If $k=1$, then the lemma holds
since $T$ and $\gamma^T$ look as depicted in Figure~\ref{trivialtangle1}c.

Let $k\ge 2$ and suppose that the lemma holds for tree claspers with
degree $\le k-1$.  Applying move~{2} to $T$ in an appropriate way, we
obtain an admissible forest clasper $T_0\cup T_1$ such that $T_0$ has
just one node, see Figure~\ref{trivialtangle1}d.  (By an appropriate rotation
of $B^3$, we may assume that $T_0$ intersects the first and second
strings of $\gamma $.)  By assumption, there is a $2k$--strand pure
braid $\beta_1$ such that
\begin{enumerate}
\item $\gamma^{T_1}$ and ${T_0}^{T_1}$ look as depicted in
Figure~\ref{trivialtangle1}e (here the framing of the (only) leaf of
${T_0}^{T_1}$ is zero),
\item $\beta_1$ is contained in $P(2k)_{k-1}$,
\item $\beta_1\setminus (\text{$2i-1$st and $2i$th strands})$ is
trivial for $i=1,\dots ,k$,
\item the first strand of $\beta_1$ is trivial and not linked with the
others  (hence $\beta_1\setminus (\text{2nd strand})$ is trivial).
\end{enumerate}
By move~{10}, the result of surgery $(\gamma
^{T_1})^{{T_0}^{T_1}}\cong\gamma^{T_0\cup T_1}\cong\gamma^T$ looks as
depicted in Figure~\ref{trivialtangle1}f, where the $(2k+2)$--strand
pure braid $\beta '_1$ is obtained from $\beta_1$ by duplicating the
first and second strands.  By the condition 4 above, ${\beta '_
1}^{-1}\setminus(\text{the 3rd and 4th strands})$ is trivial, and
hence $(\gamma ^{T_1})^{{T_0}^{T_1}}\cong\gamma^T$ is equivalent to
the tangle $\gamma'$ depicted in Figure~\ref{trivialtangle1}g.  It is
easy to see that $\beta={\beta_1'}^{-1}
{\beta_2}^{-1}\beta_1'\beta_2\in P(2k+2)$ satisfies the condition 1, 2
and 3 of Lemma~\ref{trivialtangle}.
\end{proof}

\begin{figure}
\cl{\includegraphics{trivt1.eps}}
\nocolon\caption{}
\label{trivialtangle1}
\end{figure}

By Lemma~\ref{trivialtangle}, a  ${C_k}$--move is an operation which
replaces a trivial tangle in a link into another trivial tangle.  It
is well known that a sequence of such operations can be achieved by a
set of simultaneous  operations of such kind as in
Lemma~\ref{trivialtanglemoves}. 

\begin{lemma}
\label{trivialtanglemoves}
Let ${\gamma_0},\gamma_1,\dots ,\gamma_p$ ($p\ge0$) be a sequence of
links of the same pattern in a $3$--manifold $M$.  Suppose that, for
each $i=0,\dots, p-1$, there is a $3$--ball $B_i$ in the interior of
$M$ such that the two links $\gamma_i$ and $\gamma_{i+1}$ coincide
outside $B_i$ and such that the tangles $(B_i,\gamma_i\cap B_i)$ and
$(B_i,\gamma_{i+1}\cap B_i)$ are trivial and of the same pattern.
Then there are disjoint $3$--balls $B'_0,\dots ,B'_{p-1}$ in the
interior of $M$ and diffeomorphisms $\varphi_i\colon
B_i{\overset{\cong}{\longrightarrow}} B'_i$ ($i=0,\dots, p-1$) such
that the following conditions hold.
\begin{enumerate}
\item For each $i=0,\dots,p-1$, we have
$\varphi_i(B_i\cap\gamma_i)=B'_i\cap {\gamma_0}$.
\item The link $\gamma_p$ is equivalent to the link
\begin{equation}
\label{eq:gammap}
\gamma_0\setminus (\gamma_0\cap (B'_1\cup\dots\cup B'_p))\cup
{\bigcup_{i=0}^{p-1}\varphi_i(B_i\cap \gamma_{i+1})}.
\end{equation}
\end{enumerate}
\end{lemma}

\begin{proof}
The proof is by induction on $p$.  If $p=0$, the result obviously
holds.  Let $p\ge 1$ and suppose that the lemma holds for smaller
$p$.  Thus there are disjoint $3$--balls $B'_0,\ldots,B'_{p-2}$ in
$\operatorname{int}M$ and diffeomorphisms $\varphi_{i}\colon
B_i\overset{\cong}{\longrightarrow}B_i'$ ($i=0,\ldots,p-2$) such that
$\varphi_i(B_i\cap\gamma_i)=B'_i\cap\gamma_0$, and such that
$\gamma_{p-1}$ is equivalent to the link
\begin{equation}
\label{eq:gammap-1}
\gamma_0\setminus(\gamma_0\cap(B'_1\cup\cdots\cup B'_{p-1}))\cup
\bigcup_{i=0}^{p-2}\varphi_i(B_i\cap \gamma_{i+1}).
\end{equation}
We may safely assume that $\gamma_{p-1}$ is equal to
\eqref{eq:gammap-1}.  There is a $3$--ball $B_{p-1}$ in
$\operatorname{int}M$ such that $(B_{p-1},\gamma_{p-1}\cap B_{p-1})$
and $(B_{p-1},\gamma_p\cap B_{p-1})$ are trivial tangles with the same
pattern.  Since $(B_{p-1},\gamma_{p-1})$ is a trivial tangle, there is
an isotopy $f_t\colon M\overset{\cong}{\longrightarrow}M$ fixing
$\partial M$ pointwise and fixing $\gamma_{p-1}$ as a set, such that
$f_0=\operatorname{id}_M$ and $f_1(B_{p-1})$ is disjoint from
$B'_1\cup\cdots\cup B'_{p-1}$.  We set $B'_{p-1}=f_1(B_{p-1})$ and
$\varphi_{p-1}=f_1|_{B_{p-1}}$.  Then $B'_0,\ldots,B'_{p-1}$ and
$\varphi_0,\ldots,\varphi_{p-1}$ clearly satisfies the conditions 1
the lemma.  The condition 2 follows since the link $f_1(\gamma_p)$,
which is obviously equivalent to $\gamma_p$, is equal to
$\gamma_{p-1}\setminus (\gamma_{p-1}\cap B'_{p-1}) \cup
\varphi_{p-1}(\gamma_p\cap B_{p-1})$ and hence to the link
\eqref{eq:gammap}.
\end{proof}

Using Lemmas~\ref{trivialtangle} and \ref{trivialtanglemoves}, it is
easy to verify the following.

\begin{proposition}
\label{simultaneous}
Let ${\gamma_0},\dots ,\gamma_p$ ($p\ge 0$) be a sequence of links in a
$3$--manifold $M$.  Suppose that, for each $i=0,\dots ,p-1$, the links
$\gamma_i$ and $\gamma_{i+1}$ are related by a (simple) $C_{k_i}$--move
($k_i\ge 1$).  Then there is a (simple) strict forest clasper
$T=T_0\cup\dots\cup T_{p-1}$ such that $\deg T_i=k_i$ for
$i=0,\dots,p-1$ and such that ${{\gamma_0}}^{T_1\cup\dots\cup T_{p-1}}$ is
equivalent to $\gamma_p$.
\end{proposition}

The relation on links defined by (simple) ${C_k}$--moves is symmetric as
follows.

\begin{proposition}
\label{inverse}
If a (simple) ${C_k}$--move on a link $\gamma$ in a $3$--manifold $M$
yields a link $\gamma'$ in $M$, then a (simple) ${C_k}$--move on
$\gamma'$ can yield $\gamma$.
\end{proposition}

\begin{proof}
Assume that there is a (simple) strict tree clasper $T$ for $\gamma$
of degree $k$ such that $\gamma^T\cong\gamma'$.  It suffices to show
that there is a (simple) strict tree clasper $T'$ for $\gamma^T$ of
degree $k$ disjoint from $T$ such that $(\gamma^T)^{T'}\cong\gamma$.

We choose an edge $B$ of $T$ and replace $B$ with two edges and two
trivial disk-leaves, obtaining a strict forest clasper $T_1\cup T_2$,
see Figure~\ref{inverse1}a and b.  By Proposition~\ref{trivialleaf},
we have $\gamma \cong\gamma^{T_1\cup T_2}$.  By move~{4}, we have
$\gamma^{T_1\cup T_2}\cong\gamma^{G_1}$, where $G_1$ is as depicted in
Figure~\ref{inverse1}c.  Observe that the edge $B_1$ is an
($e$--simple) good input subtree of $G_1$ of $e$--degree $k$.  By
Lemma~\ref{zipnew}, $\gamma^{G_1}$ is obtained from
$\gamma^{G_1\ominus B_1}$ by a (simple) $C_k$--move.  Clearly, we have
$\gamma^{G_1\ominus B}\cong \gamma^{T}$.  Hence $\gamma$ is obtained
from $\gamma^{T}$ by one (simple) $C_k$--move.
\end{proof}

\begin{figure}
\cl{\includegraphics{inverse1.eps}}
\nocolon\caption{}
\label{inverse1}
\end{figure}

\begin{proof}[Proof of 2$\Rightarrow$4 of Theorem~\ref{equivalence}]
Suppose that a link $\gamma$ in $M$ is ${sC_k}$--equiv\-alent to a link
$\gamma'$ in $M$.  Then there is a sequence from $\gamma$ to $\gamma '
$ of simple ${C_k}$--moves and inverses of simple ${C_k}$--moves.  By
Proposition~\ref{inverse}, the inverse simple ${C_k}$--moves are
replaced with direct simple ${C_k}$--moves.  By
Proposition~\ref{simultaneous}, such a sequence can be achieved by a
surgery on a simple strict forest clasper consisting of simple strict
tree claspers of degree $k$.
\end{proof}

\section{Structure of the set of ${C_{k+1}}$--equivalence classes of links}
\label{structureofsets}

\subsection{Set of ${C_{k+1}}$--equivalence classes of links}
It is natural and important to ask when two links of the same pattern
are ${C_k}$--equivalent.  This question decomposes inductively to the
question of when two mutually ${C_k}$--equivalent links are
${C_{k+1}}$--equivalent.  Thus the problem reduces to {\em classifying
the ${C_{k+1}}$--equivalence classes of links which are ${C_k}$--equivalent to
a fixed link ${\gamma_0}$}.  For a link $\gamma$ which is ${C_k}$--equivalent
to ${\gamma_0}$, Theorem~\ref{equivalence} enables us to measure ``how much they are
different'' by a simple strict forest clasper for ${\gamma_0}$ of degree $k$.
Hence we wish to know when two such forest claspers give
${C_{k+1}}$--equivalent results of surgeries.  

Let $M$ be a $3$--manifold, and ${\gamma_0}$ a link in $M$ of pattern
$P$.  In the following, ${\gamma_0}$ will serve as a kind of ``base
point'' or ``origin'' in the set of links which are of pattern $P$.
Let ${\calL(M,\gamma_0)}$ denote the set of equivalence classes of
links in $M$ which are of pattern $P$.  Though we have
$\calL(M,\gamma_0)=\calL(M,\gamma_0')$ for any link $\gamma_0'$ of
pattern $P$, we denote it by $\calL(M,\gamma_0)$ and not by
$\calL(M,P)$ to remember that $\gamma_0$ is the ``base point.''  We
usually write `$({\gamma_0})$' for `$(M,{\gamma_0})$' if `$M$' is
clear from context.  For each $k\ge 1$, let ${\calL_k({\gamma_0})}$
denote the subset of ${\calL(\gamma_0)}$ consisting of equivalence
classes of links which are ${C_k}$--equivalent to ${\gamma_0}$.  Then
we have the following descending family of subsets of
$\calL({\gamma_0})$
\begin{equation}
\label{eq-descending}
{\calL(\gamma_0)}\supset\calL_1({\gamma_0})\supset\calL_2({\gamma_0})\supset\cdots
\supset\calL_{\infty}({\gamma_0})\overset{\text{def}}{=}
\bigcap_{k\ge1}{\calL_k({\gamma_0})}\ni[{\gamma_0}],
\end{equation}
where $[{\gamma_0}]$ denotes the equivalence class of ${\gamma_0}$.
Conjecture~\ref{conjCinfty} is equivalent to that $\calL_{\infty}({\gamma_0})=\{[{\gamma_0}]\}$ for
any link ${\gamma_0}$ in a $3$--manifold $M$.

In order to study the descending family \eqref{eq-descending}, it is natural to
consider ${{\bar\calL}_k({\gamma_0})}=$\break
${\calL_k({\gamma_0})}/{C_{k+1}}$, the set of ${C_{k+1}}$--equivalence classes of
links in $M$ which are ${C_k}$--equivalent to the link ${\gamma_0}$.

\begin{remark}
Before proceeding to study ${{\bar\calL}_k({\gamma_0})}$, we comment
on the structure of the set ${\calL(\gamma_0)}/C_1$.  Since a simple
$C_1$--move is just a crossing change of strings, the set
${\calL(\gamma_0)}/C_1$ is identified with the set of homotopy classes
(relative to endpoints) of links that are of the same pattern as
${\gamma_0}$.  Therefore elements of ${\calL(\gamma_0)}/C_1$ are
described by the homotopy classes of the components of links.  There is
not any natural group (or monoid) structure on the set
${\calL(\gamma_0)}/C_1$ in general, but there {\em is} in the case of
string links as we will see later.
\end{remark}

\begin{definition}
Two claspers for a link ${\gamma_0}$ in $M$ are {\em isotopic with
respect to ${\gamma_0}$} if they are related by an isotopy of $M$
which preserves the set $\gamma_0$.  
Two claspers $G$ and $G'$ for a link ${\gamma_0}$ are {\em homotopic
with respect to ${\gamma_0}$} if there is a homotopy $f_t\co G\to M$ 
($t\in[0,1]$) such that 
\begin{enumerate}
\item $f_0$ is the identity map of $G$,
\item $f_1$ maps $G$ onto $G'$, respecting the decompositions into
constituents,
\item for every $t\in[0,1]$ and for every disk-leaf $D$ of $G$,
$f_t(D)$ intersects $\gamma_0$ transversely at just one point in
$f_t(\operatorname{int}D)$,
\item for each pair of two disk-leaves $D$ and $D'$ contained in one
component of $G$, the points $f_t(D)\cap\gamma_0$ and
$f_t(D')\cap\gamma_0$ are disjoint for all $t\in[0,1]$.
\end{enumerate}
\end{definition}

For $k\ge 1$, let ${\f_k({\gamma_0})}$ denote the set of simple strict
forest claspers of degree $k$ for ${\gamma_0}$.  We define a map
\begin{displaymath}
{\sigma_k}\colon {\f_k({\gamma_0})}\to{{\bar\calL}_k({\gamma_0})}
\end{displaymath}
by ${\sigma_k}({T_1\cup\dots\cup
T_p})=[{{\gamma_0}}^{{T_1\cup\dots\cup T_p}}]_{C_{k+1}}$.  Let
${\f^h_k({\gamma_0})}$ denote the quotient of ${\f_k({\gamma_0})}$ by
homotopy with respect to $\gamma_0$.

\begin{theorem}
\label{FactorThroughFhkg}
For a link ${\gamma_0}$ in a $3$--manifold $M$ and for $k\ge 1$, the map
${\sigma_k}\colon {\f_k({\gamma_0})}\to{{\bar\calL}_k({\gamma_0})}$ factors through ${\f^h_k({\gamma_0})}$.
\end{theorem}

To prove Theorem~\ref{FactorThroughFhkg}, we need some results.  The
following three Propositions are used in the proof of
Theorem~\ref{FactorThroughFhkg} and also in later sections.

\begin{proposition}
\label{slidediskleaf}
Let $T_1\cup T'_1$ be a strict forest clasper for a link $\gamma$ in a
$3$--manifold $M$ with $\deg T_1=k\ge 1$ and $\deg T_1'=k'\ge 1$.  Let
$T_2\cup T_2'$ be a strict forest clasper obtained from $T_1\cup T_1'$
by sliding a disk-leaf of $T_1$ over that of $T_1'$ along a component
of $\gamma$ as depicted in Figure~\ref{slidediskleaf1}.  Then the two links
$\gamma^{T_1\cup T_1'}$ and $\gamma^{T_2\cup T_2'}$ in $M$ are related
by one $C_{k+k'}$--move.  If, moreover, $T_1$ and $T'_1$ (and hence
$T_2$ and $T'_2$) are simple, then $\gamma^{T_1\cup T_1'}$ and
$\gamma^{T_2\cup T_2'}$ are related by one simple $C_{k+k'}$--move.
\end{proposition}

\begin{figure}
\cl{\includegraphics{slidedl1.eps}}
\nocolon\caption{}
\label{slidediskleaf1}
\end{figure}

\begin{proof}
There is a sequence of claspers for $\gamma$ from $T_2\cup T_2'$ to
a clasper $G$ as depicted in Figure~\ref{slidediskleaf2}a--f, preserving the
result of surgery, as follows.  First we obtain b from a by replacing a
simple disk-leaf of $T$ with a leaf and then isotoping it.  Then we
obtain c from b by move~{7} and by replacing a leaf with a simple
disk-leaf.  We obtain d from c by ambient isotopy, e from d by
move~{12}, and f from e by move~{6}.  Let $T$ be the good input
subtree of the clasper $G$.  The $e$--degree of $T$ is equal to
$k_1+k_2$.  By Lemma~\ref{zipnew}, $\gamma^G\cong \gamma^{T_2\cup
T_2'}$ is obtained from $\gamma^{G\ominus T}\cong \gamma^{T_1\cup
T_1'}$ by one $C_{k_1+k_2}$--move.  

If $T_1$ and $T_1'$ are simple, then the input subtree $T$ is
$e$--simple and hence $\gamma^{T_1\cup T_1'}$ is obtained from
$\gamma^{T_2\cup T_2'}$ by one simple $C_{k_1+k_2}$--move.  
\end{proof}

\begin{figure}[ht!]
\cl{\includegraphics{slidedl2.eps}}
\nocolon\caption{}
\label{slidediskleaf2}
\end{figure}

\begin{proposition}
\label{passstring}
Let $T_1$ and $T_2$ be two strict tree claspers for a link $\gamma$ of
degree $k$ in a $3$--manifold $M$ differing from each other only by a
crossing change of an edge with a component of $\gamma$.  Then
$\gamma^{T_1}$ and $\gamma^{T_2}$ are related by one ${C_{k+1}}$--move.  If,
moreover, $T_1$ and hence $T_2$ are simple, then $\gamma^{T_1}$ and
$\gamma^{T_2}$ are related by one simple ${C_{k+1}}$--move.
\end{proposition}

\begin{proof}
We may assume that $(T_1,\gamma)$ and $(T_2,\gamma)$ coincide outside
a $3$--ball in which they look as depicted in Figure~\ref{passstring1}a and b,
respectively.  There is a sequence of  claspers for $\gamma$,
preserving the results of surgery, from $T_2$ to a  clasper
$G$ as depicted in Figure~\ref{passstring1}b--d.  Here we obtain c from
b by move~{1}, and d from c by move~{12}.  Let $T$ be the good input
subtree of $G$ of $e$--degree $k+1$ as in d.  By Lemma~\ref{zipnew},
$\gamma^{G}\cong \gamma^{T_2}$ is obtained from $\gamma^{G\ominus
T}\cong \gamma^{T_1}$ by a $C_{k+1}$--move.  If $T_1$, and hence $T_2$, 
are simple, then this $C_{k+1}$--move is simple.
\end{proof}

\begin{figure}[ht!]
\cl{\includegraphics{pstring1.eps}}
\nocolon\caption{}
\label{passstring1}
\end{figure}

\begin{proposition}
\label{passedge}
Let $T_1\cup T_1'$ be a strict forest clasper for a link $\gamma$ in
$M$ with $\deg T_1=k\ge 1$ and $\deg T_1'=k'\ge 1$.  Let $T_2\cup
T_2'$ be a forest clasper for $\gamma$ obtained from $T_1\cup T_1'$ by
passing an edge of $T_1$ across that of $T_2$.  Then $\gamma^{T_1\cup
T_1'}$ and $\gamma^{T_2\cup T_2'}$ are related by one
$C_{k+k'+1}$--move.  If, moreover, $T_1$ and $T_1'$ and hence $T_2$ and
$T_2'$ are simple, then $\gamma^{T_1\cup T_1'}$ and $\gamma^{T_2\cup
T_2'}$ are related by one simple $C_{k+k'+1}$--move.
\end{proposition}

\begin{proof}
We may assume that $(T_1\cup T_1',\gamma)$ and
$(T_2\cup{T_2'},\gamma)$ coincide outside a $3$--ball in which they
look as depicted in Figure~\ref{passedge1}a and b, respectively.
(Here the $3$--ball do not intersect $\gamma$.)  We obtain from
$T_2\cup T_2'$ a clasper $G$ depicted in Figure~\ref{passedge1}d as
follows.  First we obtain c from b by move~{1}, and d from c by
move~{12} twice.  Note that the input subtree $T$ in $G$ is good and
of $e$--degree $k+k'+1$.  The rest of the proof proceeds similarly to
that of Proposition~\ref{passstring}.
\end{proof}

\begin{figure}[ht!]
\cl{\includegraphics{pedge1.eps}}
\nocolon\caption{}
\label{passedge1}
\end{figure}

\begin{proof}[Proof of Theorem~\ref{FactorThroughFhkg}]
Suppose that $T={T_1\cup\dots\cup T_p}$ and $T'=T'_1\cup\dots\cup
T'_{p'}$ ($p,p'\ge0$) are two simple strict forest claspers for
${\gamma_0}$ of degree $k$ which are homotopic to each other with
respect to ${\gamma_0}$.  We must show that
${\sigma_k}(T)={\sigma_k}(T')$, ie,
${{\gamma_0}}^{T}{\underset{C_{k+1}}{\sim}}{{\gamma_0}}^{T'}$.  By
assumption, we have $p=p'$ and, after reordering $T'_1,\dots,T'_p$ if
necessary, there is a finite sequence $G_0=T,G_1,\dots ,G_q=T'$ ($q\ge
0$) from $T$ to $T'$ of simple strict forest claspers for ${\gamma_0}$
of degree $k$ such that, for each $i=0,\dots ,q-1$, the two
consecutive simple strict forest claspers $G_i$ and $G_{i+1}$ are
related either by
\begin{enumerate}
\item isotopy with respect to ${\gamma_0}$,
\item passing an edge of a component across an edge of another
component,
\item sliding a disk-leaf of a component over a disk-leaf of another
component,
\item passing an edge of a component across the link ${\gamma_0}$,
\item passing an edge of a component across another edge of the same
component,
\item full-twisting an edge of a component.
\end{enumerate}
In each case we must show that
${{\gamma_0}}^{G_i}{\underset{C_{k+1}}{\sim}}{{\gamma_0}}^{G_{i+1}}$.
The case 1 is clear.  The cases 2, 3 and 4 comes from
Propositions~\ref{passedge}, \ref{slidediskleaf} and \ref{passstring},
respectively.  The case 5 reduces to the case 4 since passing an edge
of a component across another edge of the same component is achieved
by a finite sequence of passing an edge across ${\gamma_0}$ and
isotopy with respect to ${\gamma_0}$.  The case 6 reduces to the cases
4 and 5 since full-twisting an edge is achieved by a finite sequence
of isotopy with respect to ${\gamma_0}$, passing an edge across
another, and full-twisting an edge incident to a disk-leaf, which is
achieved by passing an edge across ${\gamma_0}$ and isotopies with
respect to ${\gamma_0}$.
\end{proof}

There is a natural monoid structure on ${\f^h_k({\gamma_0})}$ with
multiplication induced by union and with unit the empty forest
clasper.  There is a natural 1--1 correspondence between the monoid
${\f^h_k({\gamma_0})}$ and the free commutative monoid generated by
the homotopy classes with respect to ${\gamma_0}$ of simple strict
tree claspers for ${\gamma_0}$ of degree $k$. If $\pi_1M$ is finite,
then the commutative monoid ${\f^h_k({\gamma_0})}$ is finitely
generated.

Let ${\tilde\f^h_k}({\gamma_0})$ denote the (free) abelian group
obtained from the free commutative monoid ${\f^h_k({\gamma_0})}$ by
imposing the relation $[S]_h+[S']_h=0$, where $S$ and $S'$ are two
simple strict tree claspers of degree $k$ for ${\gamma_0}$ related to
each other by one half twist of an edge, and $[\thinspace \cdot
\thinspace]_h$ denotes homotopy class with respect to ${\gamma_0}$.
If $\pi_1M$ is finite, then the abelian group
${\tilde\f^h_k}({\gamma_0})$ is finitely generated.

\begin{theorem}
\label{FactorThroughFthkg}
For a link ${\gamma_0}$ in a $3$--manifold $M$ and for $k\ge 1$, the
map ${\sigma_k}\colon {\f_k({\gamma_0})}\to {{\bar\calL}_k
({\gamma_0})}$ factors through the abelian group ${\tilde\f^h_k}
({\gamma_0})$.
\end{theorem}

\begin{proof}
We have only to prove the following claim.

\begin{claim}
Let $T={T_1\cup\dots\cup T_p}$ ($p\ge 0$) be a simple strict forest
clasper for ${\gamma_0}$ in $M$ of degree $k$ and let $S$ and $S'$ be
two disjoint simple strict tree claspers for ${\gamma_0}$ of degree
$k$ which are disjoint from $T$.  Suppose that $S$ and $S'$ are
related by one half-twist of an edge and homotopy with respect to
${\gamma_0}$ in $M$.  Then the two links ${\gamma_0}^{T}$ and
${\gamma_0}^{T\cup S\cup S'}$ are ${C_{k+1}}$--equivalent.
\end{claim}

Since, by Theorem~\ref{FactorThroughFhkg}, homotopy with respect to
${\gamma_0}$ preserves the ${C_{k+1}}$--equivalence class of the result
of surgery on forest claspers of degree $k$, we may safely assume that
the $S'$ is contained in the interior of a small regular neighborhood
$N$ of $S$ in $M$.  Moreover, we may assume that $S'$ is obtained from
$S$ by a {\em positive} half twist on an edge $B$, since, if not, we
may exchange the role of $S$ and $S'$.  Let ${\gamma_N}$ denote the
link ${\gamma_0}\cap{N}$ in $N$.

Let $G=G_1\cup G_2$ be the simple strict forest clasper consisting of
two strict tree claspers $G_1$ and $G_2$ of degree $k_1$ and $k_2$,
respectively, ($k_1+k_2=k+1$) such that $G$ is obtained from $S$ by
inserting two trivial disk-leaves into the edge $B$.  By
Proposition~\ref{trivialleaf}, ${\gamma_N}^G$ is equivalent to
${\gamma_N}$.  Let $G'$ be the clasper in $N$ obtained from $G$ by
applying move~{4}.  We have ${\gamma_N}^{G'}\cong{\gamma_N}^{G}$.
Let $B$ be the edge in $G'$ that is incident to the two boxes and is
half twisted, like the edge $B_1$ in Figure~\ref{inverse1}.  Let
$\bfM$ denote the set of the two ends of $B$.  The zip construction
$\operatorname{Zip}(G',\bfM)$ consists of two components $P$ and $Q$,
satisfying the following properties.
\begin{enumerate}
\item ${\gamma_N}^{P\cup Q}\cong{\gamma_N}$.
\item $Q$ is a connected admissible clasper with ${\gamma_N}^{Q}\cong {\gamma_N}^S$.
\item $P$ is a simple strict tree clasper in $N$ for ${\gamma_N}$ of
degree $k$ such that $P$ is homotopic with respect to ${\gamma_N}$ to
$S'$ in $N$.
\end{enumerate}

Let $N_1$ be a small regular neighborhood of $N$ in $M$ which is
disjoint from $T$ and let $\gamma_1={\gamma_0}\cap N_1$.  Let $P'$ be
a simple strict tree clasper for $\gamma_1$ in $N_1\setminus N_0$ of
degree $k$ which is isotopic to $S'$, and hence to $P$, with respect
to $\gamma_1$ in $N_1$.  We have ${\gamma_1}^{P'}\cong
{\gamma_1}^{P}\cong {\gamma_1}^{S'}$.  By the construction of $P\cup
Q$, it follows that $P$ is homotopic to $P'$ with respect to
${\gamma_1}^{Q}$ in $N_1$, and hence that
$({\gamma_1}^{Q})^{P}{\underset{C_{k+1}}{\sim}}({\gamma_1}^{Q})^{P'}$.
Then we have
\begin{displaymath}
{\gamma_1}\cong{{\gamma_1}}^G\cong{\gamma_1}^{G'}\cong{\gamma_1}^{P\cup
Q}\cong({\gamma_1}^{Q})^{P}{\underset{C_{k+1}}{\sim}}({\gamma_1}^{Q})^{P'}
\cong{\gamma_1}^{Q\cup P'}\cong {\gamma_1}^{S\cup S'}
\end{displaymath}
This implies that
${\gamma_0}^T{\underset{C_{k+1}}{\sim}}{\gamma_0}^{T\cup S\cup S'}$.
This completes the proof of the claim and hence that of
Theorem~\ref{FactorThroughFthkg}.
\end{proof}

\begin{remark}
\label{remark:nk}
By Theorem~\ref{FactorThroughFthkg}, there is a surjection ${\nu_k}\colon
{\tilde\f^h_k}({\gamma_0})\to{{\bar\calL}_k({\gamma_0})}$ satisfying
${\sigma_k}={\nu_k}\circ\operatorname{proj}$, where
$\operatorname{proj}\colon {\f_k({\gamma_0})}\to{\tilde\f^h_k}({\gamma_0})$ is the
projection.
\end{remark}


\section{Groups and Lie algebras of string links}
In this section we study groups of string links in the product of a
connected oriented surface $\Sigma $ and the unit interval $[0,1]$ modulo the
${C_{k+1}}$--equivalence relation, and we also study the associated graded
Lie algebras.

In the following we fix a connected oriented surface $\Sigma$ and
distinct points $x_1,\dots,x_n$ in the interior of $\Sigma$, where $n\ge 0$.

\subsection{Definition of string links}
String links are introduced in \cite{Habegger-Lin:TheClassification}
to study link-homotopy classification of links in $S^3$.  We here
generalize this notion to string links in $\Sigma\times[0,1]$.  This
generalization is  natural and almost obvious.

\begin{definition}
An {\em $n$--string link} $\gamma=\gamma_1\cup\dots\cup\gamma_n$ in ${\Sigma\times[0,1]}$ is a
link in $\Sigma\times[0,1]$ consisting of $n$ disjoint oriented arcs
$\gamma_1,\dots,\gamma_n$, such that, for each $i=1,\dots ,n$,
$\gamma_i$ runs from $(x_i,0)$ to $(x_i,1)$.
\end{definition}

\begin{definition}
An {\em $n$--string pure braid} $\gamma$ in ${\Sigma\times[0,1]}$ is an
$n$--string link $\gamma$ in ${\Sigma\times[0,1]}$ such that, for each
$t\in[0,1]$, the surface $\Sigma\times\{t\}$ transversely intersects
$\gamma$ with just $n$ points.
\end{definition}

\begin{figure}
\cl{\includegraphics{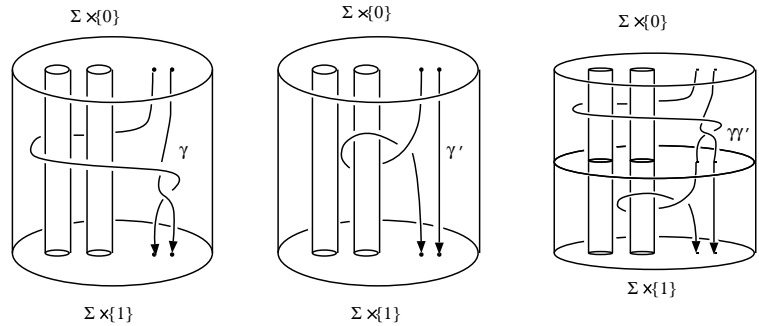}}
\caption{Example of composition of two $2$--string
links}
\label{Example:composition}
\end{figure}

Composition of $n$--string links is defined as follows.  Let
$\gamma={\gamma}_1\cup\dots\cup {\gamma}_n$ and $\gamma'
={\gamma'}_1\cup\dots\cup {\gamma'}_n$ be two string links in
${\Sigma\times[0,1]}$.  Then the {\em composition} $\gamma\gamma'
=(\gamma\gamma')_1\cup\dots\cup(\gamma\gamma')_n$ of $\gamma$ and $\gamma'$ is a
string link in ${\Sigma\times[0,1]}$ defined by
\begin{displaymath}
(\gamma\gamma')_i=h_0(\gamma_i)\cup h_1(\gamma'_i)
\end{displaymath}
for $i=1,\dots ,n$, where $h_0, h_1\colon {\Sigma\times[0,1]}\hookrightarrow{\Sigma\times[0,1]}$ are
embeddings defined by 
\begin{equation}
\label{eq:h0h1}
h_0(x,t)=(x,{\frac{1}{2}} t),\quad \mbox{and}\quad
h_1(x,t)=(x,{\frac{1}{2}}+{\frac{1}{2}} t)
\end{equation}
for $x\in\Sigma$ and $t\in[0,1]$.
For example, see Figure~\ref{Example:composition}.

The {\em trivial} $n$--string link ${1_n}$ in ${\Sigma\times[0,1]}$
consists of $n$ arcs $(1_n)_i=\{x_i\}\times [0,1],\quad i=1,\dots ,n$.
It is clear that the set ${\calL(\Sigma ,n)}$ of equivalence classes
of $n$--string links in $\Sigma\times[0,1]$ forms a monoid with multiplication induced by the
composition operation and with unit the equivalence class $[{1_n}]$ of
${1_n}$.  Here two string links are said to be {\em equivalent\/} if
they are equivalent as two links in ${\Sigma\times[0,1]}$, ie,
ambient isotopic to each other relative to endpoints.  The subset
${P{(\Sigma,n)}}$ of ${\calL(\Sigma ,n)}$ consisting of the
equivalence classes of $n$--string pure braids in ${\Sigma\times[0,1]}$
forms the unit subgroup of the monoid ${\calL(\Sigma ,n)}$, ie, the
(maximal) subgroup in ${\calL(\Sigma ,n)}$ consisting of all the
invertible elements.

\subsection{String links modulo ${C_k}$--equivalence}
For $k\ge 1$, let ${\calL_k{(\Sigma,n)}}$ denote the submonoid of ${\calL(\Sigma ,n)}$ consisting of
the equivalence classes of $n$--string links which are ${C_k}$--equivalent
to the trivial $n$--string link ${1_n}$.  That is, ${\calL_k{(\Sigma,n)}}=\calL_k({\Sigma\times[0,1]},{1_n})$.
There is a descending filtration of monoids
\begin{equation}
\label{eq:filtlsn}
{\calL(\Sigma ,n)}\supset\calL_1{(\Sigma,n)}\supset \calL_2{(\Sigma,n)}\supset \cdots.
\end{equation}
Observe that $\calL_1{(\Sigma,n)}$ is just the set of equivalence
classes of homotopically trivial $n$--string links in
${\Sigma\times[0,1]}$, where $\gamma$ is said to be {\em homotopically
trivial} if it is homotopic to ${1_n}$.   If $\Sigma$ is a disk
$D^2$ or a sphere $S^2$, then we have ${\calL(\Sigma ,
n)}=\calL_1{(\Sigma,n)}$.

Let $l\ge k$.  Let ${\calL_k{(\Sigma,n)}}/{C_l}$ denote the quotient
of ${\calL_k{(\Sigma,n)}}$ by the ${C_l}$--equiv\-al\-ence.  Also let
${\calL(\Sigma ,n)}/C_l$ denote the quotient of ${\calL(\Sigma ,n)}$ by
${C_l}$--equiv\-al\-ence.  In the obvious way, the set
${\calL_k{(\Sigma,n)}}/{C_l}$ is identified with the set of
${C_l}$--equivalence classes of $n$--string links that are
${C_k}$--equivalent to ${1_n}$.  It is easy to see that the monoid
structure on ${\calL_k{(\Sigma,n)}}$ induces that of
${\calL_k{(\Sigma,n)}}/{C_l}$.  There is a filtration on
${\calL(\Sigma ,n)}/C_l$ of finite length
\begin{equation}
\label{eq:descendinglsncl}
{\calL(\Sigma ,n)}/C_l\supset\calL_1(\Sigma ,n)/{C_l}\supset \calL_2(\Sigma ,n)/{C_l}\supset
\dots \supset \calL_l(\Sigma ,n)/{C_l}=\{1\}.
\end{equation}

Since $C_1$--equivalence is just homotopy (relative to endpoints), we
have the following.

\begin{proposition}
The monoid ${\calL(\Sigma ,n)}/C_1$ is isomorphic to the direct product
$(\pi_1\Sigma)^n$ of $n$ copies of the fundamental group
$\pi_1\Sigma $ of $\Sigma $.  Hence $\calL(\Sigma,n)/C_1$ is finitely
generated and residually nilpotent.
\end{proposition}

The following is the main result of this section.

\begin{theorem}
\label{group}
Let $\Sigma $ be a connected oriented surface.  Let $n\ge0$ and $1\le
k\le l$.  Then we have the following.
\begin{enumerate}
\item The monoid ${\calL_k{(\Sigma,n)}}/{C_l}$ is a nilpotent group.
\item The monoid ${\calL(\Sigma ,n)}/C_l$ is a residually solvable
group.  More precisely, ${\calL(\Sigma ,n)}/C_l$ is an extension of
the residually nilpotent group ${\calL(\Sigma ,n)}/C_1$ by the
nilpotent group ${\calL_1{(\Sigma,n)}}/C_l$.
\item If $\Sigma$ is a disk or a sphere, then the groups
${\calL(\Sigma ,n)}/C_l={\calL_1{(\Sigma,n)}/C_l}$ and
${\calL_k{(\Sigma,n)}}/{C_l}$ are finitely generated.
\item If $\Sigma$ is a disk or a sphere and if $n=1$, then
${\calL(\Sigma ,n)}/C_l={\calL_1{(\Sigma,n)}/C_l}$ and
${\calL_k{(\Sigma,n)}}/{C_l}$ are abelian.
\item We have
\begin{displaymath}
[{\calL_k{(\Sigma,n)}}/{C_l},\calL_{k'}(\Sigma,n)/{C_l}]
\subset\calL_{k+k'}(\Sigma,n)/{C_l},
\end{displaymath}
for $k,k'\ge 1$ with $k+k'\le l$, where $[-,-]$ denotes the commutator
subgroup.  Especially, ${\calL_k{(\Sigma,n)}}/{C_l}$ is abelian if
$1\le k\le l\le 2k$. 
\item The subgroup ${\calL_k{(\Sigma,n)}}/{C_l}$ of ${\calL(\Sigma ,
n)}/C_l$ is normal in ${\calL(\Sigma ,n)}/C_l$ (and hence in
$\calL_{k'}{(\Sigma,n)}/C_l$ with $1\le k'\le k$).  The quotient group
\[
({\calL(\Sigma ,n)}/C_l)/({\calL_k{(\Sigma,n)}}/{C_l})
\]
 is naturally
isomorphic to ${\calL(\Sigma ,n)}/{C_k}$. Similarly,
\[
(\calL_{k'}{(\Sigma,n)}/{C_l}) /({\calL_k{(\Sigma,n)}}/{C_l})
\cong\calL_{k'}{(\Sigma,n)}/{C_k}.
\]
\end{enumerate}
\end{theorem}

To prove Theorem~\ref{group}, we consider the submonoid
${{\bar\calL}_k{(\Sigma,n)}}
\overset{\operatorname{def}}{=}{\calL_k{(\Sigma,n)}}/{C_{k+1}}$
$(={{\bar\calL}_k}({\Sigma\times[0,1]},{1_n}))$ of ${\calL(\Sigma ,
n)}/{C_{k+1}}$.  We set
${\tilde\f^h_k}{(\Sigma,n)}={\tilde\f^h_k}({\Sigma\times[0,1]},{1_n})$.
By Remark~\ref{remark:nk}, there is a natural surjective map of sets
${\nu_k}\colon
{\tilde\f^h_k}{(\Sigma,n)}\to{{\bar\calL}_k{(\Sigma,n)}}$.  We have
the following lemma.

\begin{lemma}
\label{nkhomomorphism}
The map ${\nu_k}\colon {\tilde\f^h_k}{(\Sigma,n)}\to{{\bar\calL}_k{(\Sigma,n)}}$ is a surjective homomorphism of
monoids.  Hence ${{\bar\calL}_k{(\Sigma,n)}}$ is an abelian group.
\end{lemma}

\begin{proof}
First, we have
${\nu_k}(0)=[{1_n}^\emptyset]_{C_{k+1}}=[{1_n}]_{C_{k+1}}
=1_{{\bar\calL}_k{(\Sigma,n)}}$.  Second, for two elements $a$ and $b$
in ${\tilde\f^h_k}{(\Sigma,n)}$, we choose two forest claspers $T^a=
T^a_1\cup\dots\cup T^a_p$ and $T^b=T^b_1\cup\dots\cup T^b_q$ of degree
$k$ for ${1_n}$ representing $a$ and $b$, respectively.  We may assume
that $T^a\subset [0,{\frac{1}{2}}]$ and $T^b\subset \Sigma \times
[{\frac{1}{2}},1]$ since, by Theorem~\ref{FactorThroughFhkg}, homotopy
with respect to ${1_n}$ preserves the ${C_{k+1}}$--equivalence class of
results of surgeries.  Hence the forest clasper $T^a\cup T^b$
represents the element $a+b$.  We have

${\nu_k}(a+b)=[{1_n}^{T^a\cup
T^b}]_{C_{k+1}}=[{1_n}^{T^a}{1_n}^{T^b}]_{C_{k+1}}
=[{1_n}^{T^a}]_{C_{k+1}}[{1_n}^{T^b}]_{C_{k+1}}$\hfill\break
\hbox{}\hfill$={\nu_k}(a){\nu_k}(b).$

Hence $\nu_k$ is a surjective homomorphism of monoids.  Since
$\tilde{\f}^h_k (\Sigma,n)$ is a group, so is $\bar{\calL}_k
(\Sigma,n)$.
\end{proof}

The following is clear from Lemma~\ref{nkhomomorphism}.

\begin{corollary}
\label{disksphere}
If $\Sigma $ is a disk or a sphere, then ${{\bar\calL}_k{(\Sigma,n)}}$ is finitely
generated.
\end{corollary}

\begin{proof}[Proof of 1, 2, 3 and 4 of Theorem~\ref{group}]
We first prove that ${\calL_{k}{(\Sigma,n)}/C_l}$ is a group.
The proof is by a descending induction on $k$.  If $k=l$, then there is
nothing to prove.  Let $1\le k<l$ and suppose that
${\calL_{k+1}{(\Sigma,n)}/C_l}$ is a nilpotent group.  Then we have a
short exact sequence of monoids
\begin{displaymath}
1\to{\calL_{k+1}{(\Sigma,n)}/C_l}\to\calL_{k}{(\Sigma,n)}/C_l\to
{{\bar\calL}_k{(\Sigma,n)}}\to1,
\end{displaymath}	
where ${\calL_{k+1}{(\Sigma,n)}/C_l}$ and ${{\bar\calL}_k{(\Sigma,n)}}$ are  groups.  Hence ${\calL_k{(\Sigma,n)}}/{C_l}$
is also a  group.  The nilpotency is proved using the property (5) of
the theorem proved below.  
This completes the proof of 1.

The statement 2 holds since there is a  short exact sequence of monoids
\begin{displaymath}
1\to{\calL_1{(\Sigma,n)}/C_l}\to{\calL(\Sigma ,n)}/C_l\to{\calL(\Sigma ,n)}/C_1\to1.
\end{displaymath}

If $\Sigma $ is a disk or a sphere, then the group
${\calL_k{(\Sigma,n)}}/{C_l}$ is an iterated extension of finitely
generated abelian groups $\calL_k{(\Sigma,n)}/C_{k+1},\dots,
\calL_{l-1}{(\Sigma,n)}/{C_l}$.  Hence the statement 3 holds.

If $\Sigma $ is a disk or a sphere and if $n=1$, then the monoid
$\calL_1{(\Sigma,n)}$ is commutative.  Hence the statement 4 holds.
\end{proof}

Before proving the rest of Theorem~\ref{group}, we prove some results.

\begin{proposition}
\label{decompose}
Let $1\le k\le l$ and let $\gamma$ and $\gamma'$ be two $n$--string
links in ${\Sigma\times[0,1]}$ which are ${C_k}$--equivalent to each
other.  Then $\gamma'$ is ${C_l}$--equivalent to an $n$--string link
\begin{displaymath}
\gamma''=\gamma{1_n}^{T_1}\dots {1_n}^{T_p},\quad p\ge 0,
\end{displaymath}
where $T_1,\dots,T_p$ are simple strict tree claspers for ${1_n}$ such that
\begin{displaymath}
k\le \deg T_1\le \dots \le \deg T_p\le l-1.
\end{displaymath}
\end{proposition}

\begin{proof}
The proof is by induction on $l$.  If $l=k$, then there is nothing to
prove.  Let $l>k$ and suppose that $\gamma'$ is $C_l$--equivalent to
the $n$--string link $\gamma''$ given as above.  We must show that
$\gamma' $ is $C_{l+1}$--equivalent to $\gamma''{1_n}^{T_{p+1}}\dots
{1_n}^{T_{p+q}}$ ($q\ge 0$), where $T_{p+1},\dots ,T_{p+q}$ are simple
strict tree claspers for ${1_n}$ of degree $l$.  Since $\gamma''$ is
$C_{l}$--equivalent to $\gamma'$, by Theorem~\ref{equivalence} there is
a simple strict forest clasper $T'=T'_{p+1}\cup \dots \cup T'_{p+q}$
($q\ge 0$) for $\gamma''$ consisting of simple strict tree claspers of
degree $l$ such that ${\gamma''}^{T'}\cong\gamma'$.  By a homotopy
with respect to $\gamma''$ followed by an ambient isotopy fixing
endpoints, we obtain from $T'$ a simple strict forest clasper
$T''=T''_{p+1}\cup\dots\cup T''_{p+q}$ for the composition
$\gamma''{1_n}$ such that
\begin{enumerate}
\item for each $i=1,\dots, q$, $T''_i$ is contained in
$\Sigma\times[{\frac{1}{2}},1]$,
\item for each distinct $i,j\in\{1,\dots q\}$, we have $p(T''_i)\cap
p(T''_j)=\emptyset$, where $p\colon {\Sigma\times[0,1]}\to[0,1]$ is the projection.
\end{enumerate}
We have
\begin{displaymath}
(\gamma''{1_n})^{T''}\cong\gamma''{1_n}^{T''}\cong\gamma''{1_n}^{T''_{p+1}}\dots
{1_n}^{T''_{p+q}}
\end{displaymath}
(after renumbering if necessary).  By Theorem~\ref{FactorThroughFhkg},
$(\gamma''{1_n})^{T''}$ is $C_{l+1}$--equivalent to $\gamma'$.  That is,
the simple strict tree claspers $T''_{p+1},\dots,T''_{p+q}$ satisfies
the required condition.  
\end{proof}

\begin{proposition}
\label{commute}
Let $\gamma$ and $\gamma'$ be two $n$--string links in
$\Sigma\times[0,1]$ which are $C_k$--equivalent and
$C_{k'}$--equivalent, respectively, to ${1_n}$, where $k,k'\ge 1$.
Then the two compositions $\gamma\gamma'$ and $\gamma'\gamma$ are
$C_{k+k'}$--equivalent to each other.
\end{proposition}

\begin{proof}
By Proposition~\ref{decompose}, there is a simple strict forest
clasper $T=T_1\cup\dots\cup T_p$ for ${1_n}$ of degree $k$ with
${1_n}^T\cong\gamma$ and there is a simple strict forest clasper
$T'={T'}_1\cup\dots\cup{T'}_{p'}$ for $1_n$ of degree $k'$ with
${1_n}^{T'}\cong\gamma '$.  There is a sequence of claspers consisting
of simple strict tree claspers of degree $k$ or $k'$ for ${1_n}$ from
$T\cdot T'$ to $T'\cdot T$ (here we define $T\cdot T'=h_0(T)\cup
h_1(T')$ with $h_0$ and $h_1$ defined by \eqref{eq:h0h1}) such that
each consecutive two claspers are related by either one of the
following operations:
\begin{enumerate}
\item ambient isotopy fixing endpoints,
\item sliding a disk-leaf of a simple strict tree clasper of degree
$k$ with a disk-leaf of another simple strict tree clasper of degree
$k'$,
\item passing an edge of a simple strict tree clasper of degree $k$
across an edge of a simple strict tree clasper of degree $k'$.
\end{enumerate}
By Propositions~\ref{slidediskleaf} and \ref{passedge}, the result of
surgery does not change up to $C_{k+k'}$--equivalence under an
operation of the above type.  Therefore $\gamma\gamma'$ and
$\gamma'\gamma$ are $C_{k+k'}$--equivalent.
\end{proof}

\begin{proof}[Proof of 5 of Theorem~\ref{group}]
By Proposition~\ref{commute}, an element $a$ of ${\calL_k{(\Sigma,n)}}/{C_l}$ and an element $b$ of
$\calL_{k'}{(\Sigma,n)}/C_l$ commute up to $C_{k+k'}$--equivalence.  Hence the
commutator $[a,b]=a^{-1}b^{-1}ab$ is $C_{k+k'}$--equivalent to ${1_n}$.
This means that $[a,b]$ is contained in $\calL_{k+k'}{(\Sigma,n)}/C_l$.  
\end{proof}

\begin{proof}[Proof of 6 of Theorem~\ref{group}]
The subgroup ${\calL_k{(\Sigma,n)}}/{C_l}$ is normal in the subgroup
${\calL_1{(\Sigma,n)}/C_l}$ since for
$a\in{\calL_k{(\Sigma,n)}}/{C_l}$ and $b\in{\calL_1{(\Sigma,n)}/C_l}$,
we have
$b^{-1}ab=a[a,b]\in({\calL_k{(\Sigma,n)}}/{C_l})({\calL_{k+1}{(\Sigma,n)}/C_l})={\calL_k{(\Sigma,n)}}/{C_l}$.
From this fact and the fact that every $n$--string link $\gamma$ is
$C_1$--equivalent to a pure braid in ${\Sigma\times[0,1]}$, we have
only to show that ${\calL_k{(\Sigma,n)}}/{C_l}$ is closed under
conjugation of every element in ${\calL(\Sigma ,n)}/C_l$ which is
represented by a pure braid.  Let $a=[\alpha]_{C_l}\in{\calL(\Sigma ,
n)}/C_l$ be an element represented by a pure braid $\alpha$ and let
$b=[{1_n}^T]_{C_l}$ be an element of ${\calL_k{(\Sigma,n)}}/{C_l}$,
where $T$ is a simple strict forest clasper of degree $k$ for ${1_n}$.
Then the pair $(\alpha ^{-1}{1_n}\alpha, T)$ is ambient isotopic
relative to endpoints to a pair $({1_n}, T')$, where $\alpha^{-1}$ is
the inverse pure braid of $\alpha$, and $T'$ is a simple strict forest
clasper of degree $k$ for ${1_n}$.  Hence we have
$a^{-1}ba=[{1_n}^{T'}]_{C_l}\in{\calL_k{(\Sigma,n)}}/{C_l}$.
\end{proof}

\begin{remark}
Clearly, we can extend the pure braid group action on the subgroup
${\calL_k{(\Sigma,n)}}/{C_l}$ which appears in the proof of 6 of
Theorem~\ref{group} to a mapping class group action.  It is also clear
that the filtration \eqref{eq:filtlsn} is invariant under this mapping
class group action.
\end{remark}

\subsection{Lower central series of pure braid groups and groups of
string links}

Let ${P_1{(\Sigma,n)}}$ denote the subgroup of the pure braid group ${P{(\Sigma,n)}}$
consisting of equivalence classes of pure braids which are
$C_1$--equivalent to the trivial string link ${1_n}$ (ie, homotopically
trivial), ie, ${P_1{(\Sigma,n)}}={P{(\Sigma,n)}}\cap{\calL_1{(\Sigma,n)}}$.  Let
\begin{equation}
\label{eq:LHSpurebraids}
{P_1{(\Sigma,n)}}\supset P_2{(\Sigma,n)}\supset P_3{(\Sigma,n)}\supset \dots 
\end{equation}
be the lower central series of ${P_1{(\Sigma,n)}}$, which is defined by
\[
P_{k}{(\Sigma,n)}=[P_{k-1}{(\Sigma,n)},P_1{(\Sigma,n)}]
\] 
for $k\ge2$.  The following comes from Theorem~\ref{group}.

\begin{proposition}
\label{pksnlksn}
For each $k\ge 1$, we have ${P_k{(\Sigma,n)}}\subset {\calL_k{(\Sigma,n)}}$.  In other words,
every commutator of class $k$ of homotopically trivial pure braids in
${\Sigma\times[0,1]}$ is ${C_k}$--equivalent to ${1_n}$.
\end{proposition}

Now recall the definition of Stanford's equivalence relation of links
using lower central series subgroup of the (usual) pure braid group
$P(D^2,n)$ \cite{Stanford:BraidCommutators}.  

\begin{definition}
Let $M$ be a $3$--manifold.  We say that two links $\gamma$ and
$\gamma'$ are related by an element $b=[\beta]\in P(D^2,n)$ if there
is an embedding $i\colon D^2\times [0,1]\to M$ such that
$i{^{-1}}(\gamma)={1_n}$ and $i{^{-1}}(\gamma')=\beta$ as non-oriented
string links or, equivalently, as sets.  We say that two links
$\gamma$ and $\gamma'$ are {\em $P'_k$--equivalent} if $\gamma$ and
$\gamma'$ are related by an element of the $k$th lower central series
subgroup ${P_k}(D^2,n)$ of $P(D^2,n)$ for some $n\ge 0$.
\end{definition}

We can verify that $P'_k$--equivalence is actually an equivalence relation
using the fact that a pure braid in ${D^2\times[0,1]}$ is a trivial
tangle and also the fact that an element in ${P_k}(D^2,n)$ and an
element in ${P_k}(D^2,n')$ ($n,n'\ge 0$) placed `side by side' form an
element of ${P_k}(D^2,n+n')$.

The following theorem is a characterization of ${C_k}$--equivalence in
terms of pure braid commutators.

\begin{theorem}
\label{CkPnk}
Let $k\ge 0$ and let $\gamma$ and $\gamma'$ be two links in a
$3$--manifold $M$.  Then $\gamma$ and $\gamma'$ are ${C_k}$--equivalent if
and only if they are $P'_k$--equivalent.
\end{theorem}

\begin{proof}
That ${C_k}$--equivalence implies $P'_k$--equivalence follows from
Lemma~\ref{trivialtangle}.  That $P'_k$--equivalence implies ${C_k}$--equivalence
follows from Proposition~\ref{pksnlksn}.  
\end{proof}

\begin{remark}
\label{remarkPk}
We will prove in a future paper that the variant of $P'_k$--equiv\-al\-ence
which uses oriented pure braids in $D^2\times[0,1]$ instead of
non-oriented ones, which we call ``$P_k$--equivalence'', is equal to
the $P'_k$--equivalence and hence to the $C_k$--equivalence.  For knots
in $S^3$, this is derived from a recent result of
Stanford \cite{Stanford:VassilievInvariantsAndKnots}.
\end{remark}

\begin{remark}
One can redefine the notion of $P_k$--equiv\-al\-ence and
$P'_k$--equiv\-al\-ence using the lower central series of ${P_1{(\Sigma,n)}}$ for
connected oriented surface $\Sigma $.  However, it may be more
interesting to use the lower central series of ${P{(\Sigma,n)}}$, instead.
Equivalence relations thus obtained are equivalent to equivalence
relations defined using ``admissible graph claspers,'' see Section
8. 3.
\end{remark}


\subsection{Graded Lie algebras of string links}

Let ${{\hat\calL}{(\Sigma,n)}}=\varprojlim_l{\calL(\Sigma ,n)}/C_l$ and
${{\hat\calL}_k{(\Sigma,n)}}=\varprojlim_l{\calL_k{(\Sigma,n)}}/C_l$ 
($k\ge 1$) be projective limits of groups.  There is a descending
filtration of groups
\begin{displaymath}
{{\hat\calL}{(\Sigma,n)}}\supset {\hat\calL}_1{(\Sigma,n)}\supset {\hat\calL}_2{(\Sigma,n)}\supset \dots
\end{displaymath}
By construction we have
$\bigcap_{k=1}^\infty{{\hat\calL}_k{(\Sigma,n)}}=\{1\}$.  The natural
map ${\calL(\Sigma ,n)}\to{{\hat\calL}{(\Sigma,n)}}$ is injective if
and only if Conjecture~\ref{conjCinfty} holds for $n$--string links in
${\Sigma\times[0,1]}$.  If this is the case, we may think of the group
${{\hat\calL}{(\Sigma,n)}}$ as a {\em completion} of the monoid
${\calL(\Sigma ,n)}$.  However, at present, we can only say here that
${{\hat\calL}{(\Sigma,n)}}$ is a completion of the monoid
${\calL(\Sigma ,n)}/C_\infty$ of $C_\infty$--equivalence classes of
$n$--string links in ${\Sigma\times[0,1]}$.

By Theorem~\ref{group}, we
have$[{{\hat\calL}_k{(\Sigma,n)}},{\hat\calL}_{k'}{(\Sigma,n)}]\subset
{\hat\calL}_{k+k'}{(\Sigma,n)}$ for $k,k'\ge 1$.  Hence the filtration
\begin{displaymath}
{\hat\calL}_1{(\Sigma,n)}\supset {\hat\calL}_2{(\Sigma,n)}\supset \dots 
\end{displaymath}
yields the associated graded Lie algebra
$\bigoplus_{k=1}^\infty{{\hat\calL}_k{(\Sigma,n)}}/{\hat\calL}_{k+1}{(\Sigma,n)}$ with Lie bracket\nl
\hbox{}\hglue1cm$[\cdot,\cdot]\colon {{\hat\calL}_k{(\Sigma,n)}}/{\hat\calL}_{k+1}{(\Sigma,n)}\times {\hat\calL}_{k'}{(\Sigma,n)}/{\hat\calL}_{k'+1}{(\Sigma,n)}$\nl
\hbox{}\hglue7cm$\to{\hat\calL}_{k+k'}{(\Sigma,n)}/{\hat\calL}_{k+k'+1}{(\Sigma,n)}$

($k,k'\ge 1$) which maps the pair of the coset of $a$ and the coset of $b$
into the coset of the commutator $a{^{-1}}b{^{-1}}ab$.

Observe that the quotient group
${{\hat\calL}_k{(\Sigma,n)}}/{\hat\calL}_{k+1}{(\Sigma,n)}$ is
naturally isomorphic to ${{\bar\calL}_k{(\Sigma,n)}}$.  Therefore the
above graded Lie algebra structure defines that on the graded abelian
group
${\bar\calL}(\Sigma,n)=\oplus_{k=1}^\infty{{\bar\calL}_k{(\Sigma,n)}}$.
The Lie bracket
\begin{displaymath}
[\cdot,\cdot]\colon {{\bar\calL}_k{(\Sigma,n)}}\times {{\bar\calL}_{k'}(\Sigma ,
n)}\to{{\bar\calL}_{k+k'}(\Sigma ,n)}
\end{displaymath}
is given by $[[\gamma]_{C_{k+1}},[\gamma']_{C_{k'+1}}]
=[\bar{\gamma}\bar{\gamma'}\gamma\gamma']_{C_{k+k'+1}}$ for two $n$--string
links $\gamma{\underset{C_{k}}{\sim}}{1_n}$ and $\gamma'\underset{C_{k'}}{\sim}{1_n}$, where
$\bar{\gamma}$ (resp.\ $\bar{\gamma'}$) is an $n$--string link that is
inverse to $\gamma$ (resp.\ $\gamma'$) up to ${C_{k+1}}$
(resp.\ $C_{k'+1}$)--equivalence.

There is a natural action of ${\calL(\Sigma ,
n)}/C_1\cong(\pi_1\Sigma)^n$ on the graded Lie algebra
${\bar\calL}{(\Sigma,n)}$ via conjugation.

The lower central series \eqref{eq:LHSpurebraids} yields the
associated graded Lie algebra ${\bar
P}{(\Sigma,n)}=\bigoplus_{k=1}^\infty{\bar P}_k{(\Sigma,n)}$, where
${\bar P}_k{(\Sigma,n)}={P_k{(\Sigma,n)}}/P_{k+1}{(\Sigma,n)}$.  There
is an obvious homomorphism of graded Lie algebras
\begin{equation}
\label{eq:mappbsnlbsn}
i_*\colon {\bar P}{(\Sigma,n)}\to{\bar\calL}{(\Sigma,n)}.
\end{equation}

\begin{remark}
The map $i_*$ is far from surjective if $n\ge1$, as will be clear in
later sections.  If $\Sigma=D^2$, then the map $i_*$ is injective.  We
can prove this injectivity using results in the next section as
follows.  Suppose that a pure braid $\beta $ in $P_k(D^2,n)$ satisfies
$\beta{\underset{C_{k+1}}{\sim}}{1_n}$.  We must show that $\beta\in P_{k+1}(D^2,n)$.  By
Theorem~\ref{CkJk}, $\beta$ is $V_k$--equivalent to ${1_n}$, ie, $\beta$ is not
distinguished from ${1_n}$ by any invariants of type $k$.  By a theorem
of T~Kohno \cite{Kohno:VassilievInvariantsAndDeRham}, we have
$\beta\in P_{k+1}(D^2,n)$.  This completes the proof of injectivity.
\end{remark}

\begin{conjecture}
\label{istarinj}
The homomorphism $i_*\colon \bar{P}(\Sigma,n)\to\bar\calL (\Sigma,n)$
of graded Lie algebras is injective.
\end{conjecture}

\section{Vassiliev--Goussarov filtrations}


\subsection{Usual definition of the Vassiliev--Goussarov filtration}

First we recall the usual definition of the Vassiliev--Goussarov
filtration using singular links.  For details, see
\cite{Birman-Lin:KnotPolynomials} and \cite{Bar-Natan:OnTheVassiliev}.

\begin{definition}
A {\em singular link} $\gamma$ in a $3$--manifold $M$ of pattern
$P=(\alpha,i\colon$\break $\partial\alpha \hookrightarrow\partial M)$ is a
proper immersion of the $1$--manifold $\alpha$ into $M$ restricting to
$i$ on boundary such that the singularity set consists of finitely
many transverse double points.  The image of $\gamma$ is also called a
singular link of pattern $P$, and denoted by $\gamma$.  A {\em
component} of $\gamma$ is the image of a connected component of
$\alpha$ by $\gamma$.  It may happen that two distinct components of
$\gamma$ are contained in the same connected component of $\gamma$.

Two singular links $\gamma$ and $\gamma'$ of  pattern $P$ are
said to be {\em equivalent} if they are ambient isotopic relative to
endpoints.
\end{definition}

In the following, we fix a link ${\gamma_0}$ in $M$ of pattern
$P=(\alpha,i)$.

As in Section~\ref{structureofsets}, let ${\calL(M,\gamma_0)}$ denote
the set of equivalence classes of links in $M$ which are of the same
pattern with ${\gamma_0}$, and, for each $k\ge 0$, let
${\calL_k(M,{\gamma_0})}$ denote the subset of ${\calL(M,\gamma_0)}$
consisting of the equivalence classes of links which are
${C_k}$--equivalent to ${\gamma_0}$.  If $M$ is clear from the context, we
usually let ${\calL_k({\gamma_0})}$ denote ${\calL_k(M,{\gamma_0})}$.
Since $C_1$--equivalence is just homotopy (relative to endpoints),
${\calL_1({\gamma_0})}$ is the set of equivalence classes of links in
$M$ that  are homotopic to ${\gamma_0}$.

For each $k\ge 0$, let ${{{\mathcal
S}\calL}_k(M,{\gamma_0})}={{{\mathcal S}\calL}_k({\gamma_0})}$ denote
the set of equivalence classes of singular links in $M$ 
equipped with just $k$ double points, and  homotopic to
${\gamma_0}$.

For each $k\ge 0$, we construct a $\z$--linear map $e\colon
{\z{{{\mathcal S}\calL}_k({\gamma_0})}}\to\z{\calL_1({\gamma_0})}$ 
as follows.  Let $\gamma$ be a singular link with $k$ double points
which is homotopic to ${\gamma_0}$.  Let $p_1,\dots,p_k$ be
the double points of $\gamma$ and let ${\epsilon }_1,\dots,{\epsilon}_k\in\{+,-\}$ be
signs.  Let $\gamma_{({\epsilon }_1,\dots,{\epsilon}_k)}$ denote the link in $M$
obtained from $\gamma$ by replacing each double point $p_i$ with a
crossing of sign~$\epsilon_i$.  Then we set
$$e([\gamma])=\sum_{\epsilon_1,\dots,\epsilon_k\in\{+,-\}}\epsilon_
1\cdots\epsilon_k[\gamma_{(\epsilon_1,\dots,\epsilon_k)}],$$ where
$[\medspace\cdot\medspace]$ denotes equivalence class.

Let ${J_k({\gamma_0})}$ denote the subgroup of $\z{\calL_1({\gamma_0})}$ generated by the set
$X_k({\gamma_0})$ consisting of the elements $e([\gamma])$, where
$[\gamma]\in{\mathcal S}\calL_k(\gamma_0)$.  It is easy to see that the ${J_k({\gamma_0})}$'s form
a descending filtration on $\z{\calL_1({\gamma_0})}$
\begin{equation}
\label{eq:filtjk}
\z{\calL_1({\gamma_0})}=J_0({\gamma_0})\supset J_1({\gamma_0})\supset
J_2({\gamma_0})\supset \cdots,
\end{equation}
which we call the {\em Vassiliev--Goussarov filtration} on
$\z{\calL_1({\gamma_0})}$.  Later, we will redefine
${J_k({\gamma_0})}$ using claspers.

\begin{remark}
If two links ${\gamma_0}$ and ${\gamma_0}'$ are homotopic to each other, then we
have $\calL({\gamma_0} )=\calL({\gamma_0}')$ and the filtration
\eqref{eq:filtjk} is equal to the filtration
\begin{equation}
\label{eq:filtjk'}
\z\calL_1({\gamma_0}')=J_0({\gamma_0}')\supset J_1({\gamma_0}')\supset J_2({\gamma_0}')\supset \cdots.
\end{equation}
\end{remark}

\begin{remark}
We may consider a similar filtration on the abelian group
$\z{\calL(\gamma_0)}$ using singular links which are of the same
pattern with ${\gamma_0}$.  However, this filtration is the direct sum
of the filtrations on the $\z\calL_1(\gamma)$'s, where $\gamma$ runs
over a set of representatives of homotopy classes of links of pattern $P$.  Hence it suffices to study filtrations on
$\z{\calL_1({\gamma_0})}$ to study that on ${\z{\calL(\gamma_0)}}$.
\end{remark}

\begin{definition}
Let $A$ be an abelian group and $k\ge0$ an integer.  An $A$--valued
{\em invariant of type $k$} on $\calL_1(\gamma_0)$ is a homomorphism
of $\z{\calL_1({\gamma_0})}$ into $A$ which vanishes on
${J_{k+1}({\gamma_0})}$.  Thus the group of $A$--valued type $k$
invariants on ${\calL_1({\gamma_0})}$ is isomorphic to
$\operatorname{Hom}_{\z}(\z{\calL_1({\gamma_0})}/{J_{k+1}({\gamma_0})},A)$.
\end{definition}

For $k\ge0$, two links $\gamma$ and $\gamma'$ in $M$ are said to be
{\em $V_k$--equivalent}, if $[\gamma]-[\gamma']\in{J_{k+1}({\gamma_0})}$, or
equivalently, if $\gamma$ and $\gamma'$ are not distinguished by any
invariants of type $k$ with values in any abelian group.


\subsection{Definition of Vassiliev--Goussarov filtrations using claspers}

In the following, we first introduce a notion of ``schemes'' on a link
in $3$--manifolds.  Then using them we define some filtrations on
$\z{\calL_1({\gamma_0})}$ which will turn out to be equal to the
Vassiliev--Goussarov filtration defined above.

We fix a link ${\gamma_0}$ in a $3$--manifold $M$ in the following.

\begin{definition}
Let $l\ge 0$ and let $\gamma$ be a link in $M$ which is homotopic to
${\gamma_0}$.  A {\em scheme of size $l$}, $S=\{S_1,\dots,S_l\}$, for
a link $\gamma$ in a $3$--manifold $M$ is a set of $l$ disjoint
claspers for $\gamma$.  If $S_i$ is a tame clasper for every
$i=1,\dots,l$, then for each subset $S'\subset S$, the result of
surgery $\gamma^{\bigcup S'}$ is a link in $M$.  
Define an element $[\gamma,S]= [\gamma; {S_1,\dots,S_l}]$
of $\z{\calL({\gamma_0})}$ by
$$[\gamma,S]=\sum_{S'\subset
S}(-1)^{l-\operatorname{size}(S')}[\gamma^{\bigcup S'}],$$ where $S'$
runs over all $2^l$ subsets of $S$, and $[\gamma ^{\bigcup S'}]$
denotes the equivalence class of the result of surgery,
$\gamma^{\bigcup S'}$.
\end{definition}

We can easily check the following properties of the bracket notation.
\begin{itemize}
\item $[\gamma,\emptyset]=[\gamma;]=[\gamma]$,
\item $[\gamma;S_{\sigma (1)},\dots ,S_{\sigma (l)}]=[\gamma
;{S_1,\dots,S_l}]$, where $\sigma $ is a permutation on the set $\{1,\dots
,l\}$,
\item $[\gamma;{S_1,\dots,S_l}]=[\gamma^{S_1};S_2,\dots,S_l] -
[\gamma;S_2,\dots ,S_l]$,
\item $[\gamma;S_{1,1}\cup\dots\cup
S_{1,p},S_2,\dots,S_l]=\sum_{j=1}^p[\gamma^{S_{1,1}\cup\dots\cup
S_{1,j-1}};S_{1,j},S_2,\dots,S_l]$,\nl where $S_{1,1}\cup\dots\cup
S_{1,p}$ is a clasper consisting of $p$ disjoint tame claspers
$S_{1,1},\dots,S_{1,p}$.
\end{itemize}

\begin{definition}
A {\em forest scheme} $S=\{S_1,\dots,S_l\}$ is a scheme consisting
of strict tree claspers $S_1,\dots,S_l $.  We say that $S$ is {\em
simple} if the elements $S_1,\dots,S_l $ of $S$ are all simple.  The
{\em degree} $\deg S$ of a forest scheme $S$ is the sum of the degrees
of its elements.
\end{definition}	

For $k,l$ with $1\le l\le k$, let ${{J_{k,l}}({\gamma_0})}$
(resp.\ ${J^S_{k,l}({\gamma_0})}$) denote the subgroup of $\z{\calL_1({\gamma_0})}$
generated by the elements $[\gamma,S]$, where $\gamma$ is a link in
$M$ homotopic to ${\gamma_0}$ and $S$ is a forest scheme (resp.\ simple forest
scheme) of size $l$ for $\gamma$ of degree $k$.  Clearly, we have
${J^S_{k,l}({\gamma_0})}\subset {{J_{k,l}}({\gamma_0})}$.

The following theorem describes some of the inclusions of various subgroups
of $\z{\calL_1({\gamma_0})}$.  Especially, we can redefine the subgroup
$J_k({\gamma_0})$, which is previously defined using singular links, in terms
forest schemes.

\begin{theorem}
\label{redefinition}
Let ${\gamma_0}$ be a link in a $3$--manifold $M$.  Then we have the
following.
\begin{enumerate}
\item If $k\ge 1$, then we have ${J_k({\gamma_0})}=J^S_{k,k}({\gamma_0})=J_{k,k}({\gamma_0})$.
\item If $1\le l\le k$, then we have
${J^S_{k,l}({\gamma_0})}={{J_{k,l}}({\gamma_0})}$.
\item If $1\le l\le l'\le k$, then we have
$J_{k,l}({\gamma_0})\subset J_{k,l'}({\gamma_0})$.
\item If $1\le l\le k\le k'$, then we have
$J_{k',l}({\gamma_0})\subset J_{k,l}({\gamma_0})$.
\end{enumerate}
Hence we may redefine ${J_k({\gamma_0})}$ as the submodule of $\z{\calL_1({\gamma_0})}$
generated by elements $[\gamma,S]$, where $\gamma$ is a link in $M$
homotopic to ${\gamma_0}$ and $S$ is a (simple) forest scheme of degree $k$
(or of degree$\ge k$).
\end{theorem}

\begin{corollary}
\label{CkJk}
For $k\ge 0$, if two links $\gamma$ and $\gamma'$ in a $3$--manifold
$M$ are ${C_{k+1}}$--equivalent, then $\gamma$ and $\gamma'$ are
$V_k$--equivalent.
\end{corollary}

\begin{proof}
We have only to show that if $\gamma' $ is obtained from $\gamma$ by
surgery on a simple strict tree clasper of degree $k+1$, then we have
$[\gamma]-[\gamma']\in{J_{k+1}({\gamma_0})}$.  By definition, it is clear
that $[\gamma]-[\gamma']\in J^S_{k+1,1}({\gamma_0})$.  By
Theorem~\ref{redefinition}, we have $J^S_{k+1,1}({\gamma_0})\subset
{J_{k+1}({\gamma_0})}$.  Hence we have $[\gamma]-[\gamma']
\in{J_{k+1}({\gamma_0})}$.  
\end{proof}

Now we will prove Theorem~\ref{redefinition}.  The following is the
first half of the claim 1 of Theorem~\ref{redefinition}.

\begin{proposition}
\label{TwoDefsCoincide}
If $k\ge 1$, then we have ${J_k({\gamma_0})}=J^S_{k,k}({\gamma_0})$.
\end{proposition}

\begin{proof}
If $k=1$, then the result follows from Figure~\ref{proof:jkgjskkg}.
The general case follows from this case.
\end{proof}

\begin{figure}
\cl{\includegraphics{jkgjskkg.eps}}
\nocolon\caption{}
\label{proof:jkgjskkg}
\end{figure}

\begin{proof}[Proof of 2 of Theorem~\ref{redefinition}]
It suffices to show that ${{J_{k,l}}({\gamma_0})}\subset
{J^S_{k,l}({\gamma_0})}$.  Let $\gamma$ be a link in $M$ homotopic to
${\gamma_0}$ and let $S=\{{S_1,\dots,S_l}\}$ be a forest scheme of
size $l$ for $\gamma$ of degree $k$.  For each $i=1,\dots ,l$, let
$N_i$ be a small regular neighborhood of $S_i$ in $M$ such that
$\partial N_i$ is transverse to $\gamma$.  By
Theorem~\ref{equivalence}, there are finitely many disjoint simple
strict tree claspers $T_{i,1},\dots , T_{i,p_i}$ ($p_i\ge 0$) for
$\gamma\cap N_i$ in $N_i$ of degree $\deg S_i$ such that $(\gamma \cap
N_i)^{T_{i,1}\cup \dots \cup T_{i,p_i}}\cong(\gamma\cap N_i)^{S_i}$,
and hence $\gamma^{T_{i,1}\cup \dots \cup T_{i,p_i}\cup G}\cong
\gamma^{S_i\cup G}$ for any tame clasper $G$ for $\gamma$ which is
disjoint from $N_i$.  Therefore we have $[\gamma;
{S_1,\dots,S_l}]=[\gamma;T_{1,1}\cup \dots \cup T_{1,p_1},\dots ,
T_{l,1}\cup \dots \cup T_{l,p_l}]$
$=\sum_{j_1=1}^{p_1}\cdots\sum_{j_l=1}^{p_l}[\gamma^{\bigcup_{1\le
i\le l}\bigcup_{1\le j\le j_i}T_{i,j}};T_{1,j_1},\dots , T_{l,j_l}]$
$\in{J^S_{k,l}({\gamma_0})}$.  
\end{proof}

\begin{proof}[Proof of 1 of Theorem~\ref{redefinition}]
The first half is Proposition~\ref{TwoDefsCoincide}.  The rest comes
from the claim 2.
\end{proof}

\begin{proof}[Proof of 3 of Theorem~\ref{redefinition}]
For $k,l$ with $1\le l\le k-1$, we have only to prove that
${{J_{k,l}}({\gamma_0})}\subset{J_{k,l+1}({\gamma_0})}$.

By the claim 2, it suffices to show that, for a simple forest scheme
$S=\{S_1,\dots,S_l \}$ of size $l$ of degree $k$ for a link $\gamma$
in $M$ with $[\gamma]\in{\calL_1({\gamma_0})}$, we have
$[\gamma,S]\in J_{k,{l+1}}({\gamma_0})$.  By hypothesis, there is an element, say
$S_1$, of $S$ of degree $\ge 2$.  Let $V$ be a node of $S_1$ adjacent
to at least one disk-leaf.  Applying move~{9} and ambient isotopy to
the simple strict tree clasper $S_1$ and the node $V$, we obtain two
disjoint strict tree claspers $S_{1,1}$ and $S_{1,2}$ for $\gamma$ in
a small regular neighborhood $N$ of $S_1$ in $M$ such that
\begin{enumerate}
\item $\deg S_{1,1}+\deg S_{1,2}=\deg S_1$,
\item ${\gamma_N}\cong{\gamma_N}^{S_{1,1}}\cong{\gamma_N}^{S_{1,2}}$,
\item ${\gamma_N}^{S_{1,1}\cup S_{1,2}}\cong{\gamma_N}^{S_1}$,
\end{enumerate}
where ${\gamma_N}$ denotes the link $\gamma\cap N$ in $N$.  Hence we
have $[{\gamma_N};S_1]=$\nl
$[{\gamma_N};S_{1,1}\cup S_{1,2}]
=[{\gamma_N}^{S_{1,1}\cup S_{1,2}}]-[{\gamma_N}]
=[{\gamma_N}^{S_{1,1}\cup
S_{1,2}}]-[{\gamma_N}^{S_{1,1}}]-[{\gamma_N}^{S_{1,2}}]+[{\gamma_N}]$\nl
$=[{\gamma_N};S_{1,1},S_{1,2}]$.\par  Therefore $[\gamma,S]=[\gamma ;
S_1,S_2,\dots ,S_l] =[\gamma; S_{1,1},S_{1,2},S_2,\dots ,S_l]
\in{J_{k,l+1}({\gamma_0})}$.  
\end{proof}

\begin{proof}[Proof of 4 of Theorem~\ref{redefinition}]
We have to prove that if $1\le l\le k$, then
${J_{k+1,l}({\gamma_0})}\subset{{J_{k,l}}({\gamma_0})}$.

It suffices to show that, for a link $\gamma$ in $M$ and a simple
forest scheme $S=\{S_1,\dots,S_l\}$ of size $l$ for $\gamma$ of degree
$k+1$, we have $[\gamma,S]\in{{J_{k,l}}({\gamma_0})}$.  There is at least one
element, say $S_1$, of $S$ with degree $\ge2$.  By
Proposition~\ref{leqkeq} and Theorem~\ref{equivalence}, there is a
strict tree clasper $S_1'$ of degree $\deg{S_1}-1$ for $\gamma$
contained in a small regular neighborhood $N$ of $S_1$ in $M$ such
that ${\gamma_N}^{S_1}\cong{\gamma_N}^{S_1'}$, where
${\gamma_N}=\gamma\cap{N}$.  Hence we have $[\gamma,S]
=[\gamma;S_1,\dots,S_l]
=[\gamma;S_1',S_2,\dots,S_l]\in{{J_{k,l}}({\gamma_0})}$.  
\end{proof}	

\subsection{Vassiliev--Goussarov filtrations on string links}
In this section we study the Vassiliev--Goussarov filtration on
$n$--string links in ${\Sigma\times[0,1]}$, where $n\ge 0$ and $\Sigma$
is a connected oriented surface.  For $k\ge 0$, we set
${J_k{(\Sigma,n)}}={J_k}({\Sigma\times[0,1]},{1_n})$.  This defines a
descending sequence of two-sided ideals of the monoid ring
$\z{\calL_1{(\Sigma,n)}}$
\begin{displaymath}
\z{\calL_1{(\Sigma,n)}}=J_0{(\Sigma,n)}\supset J_1{(\Sigma,n)}\supset \cdots.
\end{displaymath}
We also set
${J_{k,l}{(\Sigma,n)}}=J_{k,l}({\Sigma\times[0,1]},{1_n})$.
As an alternative, we may define ${J_k{(\Sigma,n)}}$
(resp.\ ${J_{k,l}{(\Sigma,n)}}$) to be the subgroup of
$\z{\calL_1{(\Sigma,n)}}$ generated by the set of the elements
$[\gamma, S]$, where $\gamma$ is a homotopically trivial $n$--string
link and $S$ is a forest scheme (resp.\ a forest scheme of size $l$)
for $\gamma$ of degree $k$.  (We may assume that $S$ is simple.  See
Theorem~\ref{redefinition}.)

As we can easily see, ${J_{k,l}{(\Sigma,n)}}$ (and hence
${J_k{(\Sigma,n)}}$) is a two-sided ideal in the monoid ring
$\z{\calL_1{(\Sigma,n)}}$.  Moreover we have
$${J_{k,l}{(\Sigma,n)}}J_{k',l'}{(\Sigma,n)} \subset
J_{k+k',l+l'}{(\Sigma,n)}$$ and, especially,
${J_k{(\Sigma,n)}}J_{k'}{(\Sigma,n)}\subset
J_{k+k'}{(\Sigma,n)}$. Observe that $J_1{(\Sigma,n)}$ is the
augmentation ideal of the monoid ring $\z{\calL_1{(\Sigma,n)}}$, ie,
\[J_1{(\Sigma,n)}=\operatorname{ker}(\epsilon \colon
\z{\calL_1{(\Sigma,n)}}\to\z),\] where $\epsilon $ is given by
$\epsilon ([\gamma])=1$ for every $[\gamma]\in\calL_1 (\Sigma,n)$.

For two schemes $S=\{S_1,\dots,S_l \}$ and $S'=\{S'_1,\dots,S'_{l'}\}$
for $n$--string links $\gamma$ and $\gamma'$, respectively, let $S\cdot
S'$ denote the scheme of size $l+l'$ for the composition
$\gamma\gamma'$ defined by $S\cdot S'=h_0(S)\cup h_1(S')$, where
$h_0,h_1\colon \Sigma\times[0,1]\hookrightarrow \Sigma\times[0,1]$ are
as in \eqref{eq:h0h1}.

For $k,l$ with $1\le l\le k$, we set
\begin{displaymath}
{J^{(l)}_k{(\Sigma,n)}} =\sum_{\substack{k_1,\dots,k_l\ge1\\
k_1+\dots+k_l=k}} J_{k_1,1}{(\Sigma,n)}\cdots J_{k_l,1}{(\Sigma,n)}.
\end{displaymath}

\begin{proposition}
\label{newJ1J1}
If $k\ge 1$, then
\begin{equation}
\label{eqJ1J1}
{J_k{(\Sigma,n)}}\subset \sum_{l=1}^k {J^{(l)}_k{(\Sigma,n)}}.
\end{equation}
\end{proposition}

\begin{proof}
The proof is by induction on $k$.  If $k=1$, then \eqref{eqJ1J1}
clearly holds.  Let $k\ge 2$ and assume that \eqref{eqJ1J1} holds for
smaller $k$.  In this proof, we set $N_k=\sum_{l=1}^k
{J^{(l)}_k{(\Sigma,n)}}$.  It suffices to prove the following claim. 

\begin{claim}
Let $S=\{S_1,\dots,S_m\}$ be a simple strict forest scheme of size $m$
of degree $k$ for a homotopically trivial $n$--string link $\gamma$.
Then we have $[\gamma,S]\in N_k$.
\end{claim}

The proof of the claim is by induction on $m$.  If $m=1$, then
$[\gamma, S]\in J_{k,1}{(\Sigma,n)}\subset {N_k}$.  Let $m\ge 2$ and
assume that the claim holds for smaller $m$.

We first prove that we may assume $\gamma={1_n}$.  By
Theorem~\ref{group}, there is an $n$--string link ${\bar\gamma}$ which
is inverse to $\gamma$ up to ${C_k}$--equivalence.  Hence there is a
simple strict forest clasper $T= T_1\cup\dots\cup T_q$ for
${\bar\gamma}\gamma$ of degree $k$ such that $T$ is disjoint from
$\emptyset\cdot S$ and such that $({\bar\gamma}\gamma)^T\cong{1_n}$.
Since ${\bar\gamma}$ is $C_1$--equivalent to ${1_n}$, there is, by
Theorem~\ref{equivalence}, a simple strict forest clasper $C=
C_1\cup\dots\cup C_p$ of size $p$ for ${1_n}$ of degree $1$ such that
${1_n}^C\cong{\bar\gamma}$.  We have $[\gamma,
S]=[{1_n}][\gamma,S]=[{\bar\gamma}][\gamma,
S]-([{\bar\gamma}]-[{1_n}])[\gamma,
S]=[{\bar\gamma}\gamma,\emptyset\cdot S]-[{1_n};C][\gamma,
S]=[({\bar\gamma}\gamma) ^T;\emptyset\cdot
S]-[{\bar\gamma}\gamma,\{T\}\cup(\emptyset\cdot
S)]-[{1_n};C][\gamma,S]$.  It is clear that
$[{\bar\gamma}\gamma,\{T\}\cup(\emptyset\cdot S)]\in
J_{k,1}{(\Sigma,n)}$.  We have also
$[{1_n};C][\gamma,S]\in{J_1{(\Sigma,n)}}{J_k{(\Sigma,n)}}\subset
{J_1{(\Sigma,n)}}J_{k-1}{(\Sigma,n)}$.  By the induction hypothesis,
we have $J_{k-1}{(\Sigma,n)}\subset{N_{k-1}}$.  Hence
$[{1_n};C][\gamma,S]\in J_{1,1}{(\Sigma,n)}{N_{k-1}}\subset {N_k}$.
Therefore we have only to show that $[({\bar\gamma}
\gamma)^T;\emptyset\cdot S]\in{N_k}$.  This means that we have to
prove the claim only in the case that $\gamma={1_n}$.

Assume that $\gamma={1_n}$.  By
ambient isotopy, we may assume that $S_1,\dots,S_m\subset \Sigma\times
[{\frac{1}{2}},1]$.  There is a sequence of simple strict tree
claspers $S_{1,0},\dots , S_{1,r}$ ($r\ge 0$) for ${1_n}$ of degree
$k_1=\deg S_1$ satisfying the following conditions.
\begin{enumerate}
\item $S_{1,0}=S_1,\quad S_{1,r}\subset [0,{\frac{1}{2}}]$.
\item $S_{1,i}$ is disjoint from ${S_2,\dots,S_m}$ for $i=0,\dots ,r$.
\item For each $i=0,\dots ,r-1$, $S_i$ and $S_{i+1}$ are related by
one of the following operations:
\begin{enumerate}
\item ambient isotopy fixing $S_2\cup \dots \cup S_m$ pointwise and
${1_n}$ as a set,
\item sliding a disk-leaf of $S_{1,i}$ over a disk-leaf of some $S_j$
($2\le j\le m$),
\item passing an edge of $S_{1,i}$ across an edge of some $S_j$
($2\le{j}\le m$).
\end{enumerate}
\end{enumerate} 

We set
$d_i=[{1_n};S_{1,i+1},{S_2,\dots,S_m}]-[{1_n};S_{1,i},{S_2,\dots,S_m}]$
for $i=0,1,\dots ,r-1$.  We must show that $d_0,\dots,d_{r-1}$ and
$[{1_n};S_{1,r},{S_2,\dots,S_m}]$ are contained in ${N_k}$.  Since
$S_{1,r}\subset [0,{\frac{1}{2}}]$ and ${S_2,\dots,S_m}\subset
[{\frac{1}{2}},1]$, we have
$[{1_n};S_{1,r},{S_2,\dots,S_m}]=[{1_n};S_{1,r}][{1_n};{S_2,\dots,S_m}]\in
J_{k_1,1}{(\Sigma,n)} J_{k-k_1}{(\Sigma,n)}$.  
By hypothesis, $J_{k-k_1}{(\Sigma,n)}$ is contained in $N_{k-k_1}$. 
 Hence $[{1_n};S_{1,r},{S_2,\dots,S_m}]\in{N_k}$.  Now it suffices to show
that $d_i\in{N_k}$ in each of the three cases a, b and c above.

{\bf Case (a)}\qua  We clearly have $d_i=0$.

{\bf Case (b)}\qua  Suppose that $S_{1,i+1}$ is obtained from $S_{1,i}$
by sliding a disk-leaf $D_1$ of $S_{1,i}$ over a disk-leaf $D_2$ of
$S_j$ ($2\le j\le l$).  We may assume $j=2$ without loss of
generality.  Let $c \subset {1_n}$ be the segment in ${1_n}$ bounded
by $D_1\cap {1_n},D_2\cap {1_n}$, along which the slide occurs, and
let $N$ be a small regular neighborhood of $S_{1,i}\cup S_2\cup c$,
with $\partial N$ transverse to ${1_n}$. We may assume that
$S_{1,i+1}\subset {\operatorname{int}} N$.  By
Proposition~\ref{slidediskleaf} and Theorem~\ref{equivalence}, there
is a simple strict tree clasper $T$ for ${\gamma_N}={1_n}\cap N$ of
degree $k_1+k_2$ disjoint from $S_{1,i}\cup S_2$ such that
${\gamma_N}^{S_{1,i+1}\cup S_2}\cong{\gamma_N}^{S_{1,i}\cup S_2\cup
T}$.  Hence $[{\gamma_N};S_{1,i+1},S_2]-[{\gamma_N};S_{1,i},S_2]
=([{\gamma_N}^{S_{1,i+1}\cup
S_2}]-[{\gamma_N}^{S_{1,i+1}}]-[{\gamma_N}^{S_2}]+[{\gamma_N}]) -
([{\gamma_N}^{S_{1,i}\cup
S_2}]-[{\gamma_N}^{S_{1,i}}]-[{\gamma_N}^{S_2}]+[{\gamma_N}])
=[{\gamma_N}^{S_{1,i}\cup S_2\cup T}]-[{\gamma_N}^{S_{1,i}\cup {S_2}}]
=[{\gamma_N}^{S_{1,i}\cup S_2};T]$.  Therefore we have
$d_i=[{1_n}^{S_{1,i}\cup S_2};T,{S_3,\dots,S_m}]\in{N_k}$ by the
induction hypothesis.

{\bf Case (c)}\qua  Suppose that $S_{1,i+1}$ is obtained from $S_{1,i}$
by a crossing change of an edge of $S_{1,i}$ and an edge of $S_j$
($2\le j \le m$).  We may assume $j=2$ without loss of generality.
Let $B$ be a small $3$--ball in
${\operatorname{int}}({\Sigma\times[0,1]})$ in which the crossing
change occurs, and let $N$ be a small regular neighborhood of
$B\cup{}S_{1,i}$ in ${\Sigma\times[0,1]}$ with $\partial N$ transverse
to ${1_n}$.  We may assume that $S_{1,i+1}\subset {\operatorname{int}}
N$.  By Proposition~\ref{passedge} and Theorem~\ref{equivalence},
there is a simple strict tree clasper $T$ for ${\gamma_N}={1_n}\cap N$
of degree $k_1+k_2+1$ disjoint from $S_{1,i}$ and $S_2$ such that
${\gamma_N}^{S_{1,i+1}\cup S_2}\cong{\gamma_N}^{S_{1,i}\cup S_2\cup
T_0}$.  Modifying $T_0$ by move~{10}, we obtain a simple strict tree
clasper $T$ in a small regular neighborhood of $T_0$ for ${\gamma_N}$
of degree $k_1+k_2$ such that ${\gamma_N}^{S_{1,i+1}\cup
S_2}\cong{\gamma_N}^{S_{1,i}\cup S_2\cup T}$.  Then we can check
$d_i\in{N_k}$ as in the case b.
\end{proof}

\begin{theorem}
\label{CinftyJinfty}
If $n\ge 0$ and if $\Sigma$ is a connected oriented surface, then we
have
\begin{equation}
\label{eqCinftyJinfty}
\bigcap_{k=0}^\infty{J_k{(\Sigma,n)}}=\bigcap_{k=1}^\infty J_{k,1}{(\Sigma,n)}
=J_{\infty,1}(\Sigma,n).
\end{equation}
Here $J_{\infty,1}(\Sigma,n)$ is the subgroup of $\z\calL_1(\Sigma,n)$
generated by the elements of the form $[\gamma]-[\gamma']$, where
$\gamma$ and $\gamma'$ are $C_\infty$--equivalent.
\end{theorem}

\begin{proof}
It is clear that $\bigcap_{k=1}^\infty J_{k,1}{(\Sigma,n)}\subset
\bigcap_{k=0}^\infty{J_k{(\Sigma,n)}}$.  We will prove the reverse
inclusion.  We must show that, for each $k\ge 1$, we have
$\bigcap_{k=0}^\infty{J_k{(\Sigma,n)}}\subset J_{k,1}{(\Sigma,n)}$.

Let $i\colon \z{\calL_1{(\Sigma,n)}}/J_{k,1}{(\Sigma,n)}
{\overset{\cong}{\longrightarrow}}\z({\calL_1{(\Sigma,n)}}/{C_k})$ be
the obvious isomorphism of rings.  Let $J$ denote the augmentation
ideal of the group ring $\z({\calL_1{(\Sigma,n)}}/{C_k})$.  For each
$l\ge 1$, we have
$i^{-1}(J^l)=({J_1{(\Sigma,n)}}^l+J_{k,1}{(\Sigma,n)})/J_{k,1}{(\Sigma,n)}$.
Since $J$ is the augmentation ideal of the group ring of the {\em
nilpotent} group ${\calL_1{(\Sigma,n)}}/{C_k}$, we have
$\bigcap_{l=0}^\infty J^l=\{0\}$.  Hence we have $\bigcap_{l=0}^\infty
(({J_1{(\Sigma,n)}}^l+J_{k,1}{(\Sigma,n)})/J_{k,1}{(\Sigma,n)})=\{0\}$.
This implies that
$\bigcap_{l=0}^\infty({J_1{(\Sigma,n)}}^l+J_{k,1}(\Sigma,n))\subset
J_{k,1}{(\Sigma,n)}$.  
Hence we have $\bigcap_{l=1}^\infty
J_l(\Sigma,n)\subset\bigcap_{l=1}^\infty J_{kl}(\Sigma,n)\subset
\bigcap_{l=1}^\infty \sum_{i=1}^{kl} J_{kl}^{(i)}(\Sigma,n) \subset
\bigcap_{l=1}^\infty (J_1 (\Sigma,n)^l+ J_{k,1}(\Sigma,n)) \subset
J_{k,1}(\Sigma,n)$.

The equality $\bigcap_{k=1}^\infty J_{k,1}{(\Sigma,n)}
=J_{\infty,1}(\Sigma,n)$ is obvious.
\end{proof}

\begin{corollary}
\label{conj:Cinfinity}
Two $n$--string links in ${\Sigma\times[0,1]}$ are $C_\infty$--equivalent to each other
if and only if they are not distinguished by any finite type
invariants.
\end{corollary}

\begin{conjecture}
\label{conj:ckkVk}
Let $\Sigma$ be a connected oriented surface and let $n,k\ge 0$.  Then
two $n$--string links in ${\Sigma\times[0,1]}$ are ${C_{k+1}}$--equivalent if and only if
they are $V_k$--equivalent.
\end{conjecture}

If Conjecture~\ref{conj:ckkVk} holds, then Theorem~\ref{CinftyJinfty}
can be proved as a corollary to Conjecture~\ref{conj:ckkVk}.

\subsection{Vassiliev--Goussarov filtrations on string knots}

\begin{definition}
A {\em string knot} will mean a $1$--string link in the cylinder
$D^2\times$\break$[0,1]$.  
Let ${\calL(1)}=\calL_1(1)$ denote the commutative monoid
$\calL(D^2,1)=\calL_1(D^2,1)$ of string knots.
\end{definition}

There is a natural isomorphism between the monoid ${\calL(1)}$ and the
monoid of equivalence classes of knots in $S^3$ with multiplication
induced by the connected sum operation.  Therefore all results about
string knots in the rest of this section can be directly restated for
knots in $S^3$.

\begin{definition}
A $\z$--linear map $v\colon \z{\calL(1)}\to A$, where $A$ is an abelian
group, is {\em additive} (or {\em primitive}) if $v([1_1])=0$ and
$v({J_1(1)}{J_1(1)})=0$.
\end{definition}
Since the augmentation ideal ${J_1(1)}$ of $\z{\calL(1)}$ is spanned
by $\{[\gamma]-[{1_1}]\thinspace|\thinspace[\gamma] \in\calL(1)\}$,
the condition that $v({J_1(1)}{J_1(1)})=0$ is equivalent to
$v(([\gamma]-[{1_1}])([\gamma']- [{1_1}]))=0$, and hence to
$v([\gamma\gamma'])- v([\gamma])- v([\gamma'])+v([{1_1}])=0$ for any
two string knots $\gamma$ and $\gamma'$.  This is equivalent to
$v([\gamma\gamma'])=v([\gamma])+v([\gamma'])$ by the first condition.
Conversely, $v([\gamma\gamma'])=v([\gamma])+v([\gamma'])$ implies the
additivity of $v$.  In other words, $v$ is additive if and only if $v$
restricts to a homomorphism ${\calL(1)}\to A$ of commutative monoids.

Let ${\psi_k}\colon \z{\calL(1)}\to{\calL(1)}/{C_{k+1}}$ ($k\ge 0$) be
the homomorphism of abelian groups defined by
${\psi_k}([\gamma])=[\gamma]_{C_{k+1}}$ for each string knot $\gamma$,
where $[\gamma]_{C_{k+1}}$ denotes the ${C_{k+1}}$--equivalence class
of $\gamma$.  (${\psi_k}$ is a homomorphism of an additive group into
a multiplicative group.)

\begin{proposition}
\label{additivetypek}
For each $k\ge 0$, the homomorphism ${\psi_k}$ is an additive
invariant of type $k$.
\end{proposition}

\begin{proof}
For two string knots $\gamma$ and $\gamma'$, we have
\[
{\psi_k}([\gamma\gamma']-[\gamma]-[\gamma'])=[\gamma\gamma']_{C_{k+1}}[\gamma]_{C_{k+1}}^{-1}[\gamma']_
{C_{k+1}}^{-1}=[{1_1}]_{C_{k+1}}.
\]
  Hence $\psi_k$ is additive.

Now we prove that ${\psi_k}$ is of type $k$, ie,
${\psi_k}({J_{k+1}(1)})=\{1\}$.  By Proposition~\ref{newJ1J1}, we have
${J_{k+1}(1)}\subset J_{{k+1},1}(1)+J_1(1)J_1(1)$.  Clearly, $\psi_k$
vanishes on $J_{k+1,1}(1)$.  By the additivity of ${\psi_k}$ proved
above, ${\psi_k}$ vanishes also on ${J_1(1)}{J_1(1)}$.  Hence
${\psi_k}$ is of type $k$.
\end{proof}

\begin{theorem}
\label{UniversalAdditiveVassilievInvariant}
For $k\ge 1$, ${\psi_k}$ is {\em universal} in that for any additive
invariant $v\colon \z{\calL(1)}\to A$ of type $k$ with values in any abelian group
$A$, there is a unique homomorphism $\bar v\colon {\calL(1)}/{C_{k+1}}\to A$ such that
$v=\bar v {\psi_k}$.
\end{theorem}

\begin{proof}
Let $v\colon \z{\calL(1)}\to A$ be an additive invariant of type $k$
with $A$ an abelian group.  First we prove the uniqueness of $\bar v$.
Suppose that $\bar v,{\bar v}'\colon {\calL(1)}/{C_{k+1}}\to A$ are
two homomorphisms with $v=\bar v{\psi_k}={\bar v}'{\psi_k}$.  Then,
for each string knot $\gamma$, we have $\bar
v([\gamma]_{C_{k+1}})=\bar v{\psi_k}([\gamma])={\bar
v}'{\psi_k}([\gamma])={\bar v}'([\gamma]_{C_{k+1}})$.  Hence $\bar
v={\bar v}'$.  Next we prove the existence of $\bar v$.  By the
additivity of $v$, the restriction $v|_{\calL(1)}\colon {\calL(1)}\to
A$ is a homomorphism of monoids.  The restriction $v|_{\calL(1)}$
factors through $\psi_k|{\calL(1)}=\operatorname{proj}\colon
{\calL(1)}\to{\calL(1)}/{C_{k+1}}$, since if two string knots $\gamma$
and $\gamma'$ are ${C_{k+1}}$--equivalent, then we have $v([\gamma])
=v([\gamma'])$.  Hence there is a homomorphism $\bar v\colon
{\calL(1)}/{C_{k+1}}\to A$ such that
$\bar{v}({\psi_k}|_{\calL(1)})=v|_{\calL(1)}$, and hence such that
$\bar v{\psi_k}=v$.
\end{proof}

The following theorem gives a characterization of the information
carried by invariants of type $k$ in terms of ${C_{k+1}}$--equivalence.

\begin{theorem}
\label{JkCk}
If $k\ge 0$ and if $\gamma$ and $\gamma'$ are string knots, then the following
conditions are equivalent.
\begin{enumerate}
\item $\gamma$ and $\gamma'$ are ${C_{k+1}}$--equivalent.
\item $\gamma$ and $\gamma'$ are $P'_{k+1}$--equivalent.
\item $\gamma$ and $\gamma'$ are $V_k$--equivalent, ie, $\gamma$ and
$\gamma'$ are not distinguished by any invariant of type $k$ with
values in any abelian group.
\item $\gamma$ and $\gamma'$ are not distinguished by any additive
invariant of type $k$ with values in any abelian group.
\end{enumerate}	
Similar equivalence holds also for knots in $S^3$.
\end{theorem}

\begin{proof}
By Theorem~\ref{CkPnk}, 1 and 2 are equivalent.  By
Corollary~\ref{CkJk}, 1 implies 3.  It is clear that 3 implies 4.  By
Theorem~\ref{UniversalAdditiveVassilievInvariant}, 4 is equivalent to
${\psi_k}([\gamma])={\psi_k}([\gamma'])$.  This implies
$[\gamma]_{C_{k+1}}=[\gamma']_{C_{k+1}}$, and hence the condition 1.
\end{proof}

\begin{remark}
That 2 implies 3 is due to T~Stanford
\cite{Stanford:BraidCommutators}.  The above proof using claspers
provides another (very indirect) proof of this.

After a previous version of this paper
\cite{Habiro:ClaspersAndTheVassilievSkeinModules} was circulated,
T~Stanford proved that two knots in $S^3$ are $V_k$--equivalent if and
only if they are represented as two closed braids of the same number
of strands which differ only by an element of the $k+1$st lower
central series subgroup of the pure braid group
\cite{Stanford:VassilievInvariantsAndKnots}.  The equivalence of 2 and
3 in Theorem~\ref{JkCk} can be also derived from this result of
Stanford.  That 3 is equivalent to 4 is due to Stanford
\cite{Stanford:VassilievInvariantsAndKnots}.

The techniques used in \cite{Stanford:VassilievInvariantsAndKnots}
deeply involves commutator calculus on the pure braid groups and, at
first sight, they may look very different from the techniques used in
this paper (and in \cite{Habiro:ClaspersAndTheVassilievSkeinModules}).
However, they are related to each other in some deep sense.
Stanford's proof involves commutator calculus on pure braid groups,
while our proof implicitly involves commutator calculus on a Hopf
algebra in a category of $3$--dimensional cobordisms.  See
Section~\ref{ss:BraidedCommutatorCalculus} for more details.
\end{remark}

\begin{remark}
Since a rational invariant of type $k$ is a sum of an additive
invariant of type $k$ and a polynomial of invariants of degree $<k$
\cite{Kontsevich:VassilievsKnotInvariants}
\cite{Bar-Natan:OnTheVassiliev}, the conditions 3 and 4 above are
equivalent for rational finite type invariants.  This fact is noted in
\cite{Stanford:VassilievInvariantsAndKnots}.
\end{remark}

\begin{remark}
It is well known that there is an algorithm to determine whether or
not two given knots $\gamma$ and $\gamma'$ in $S^3$ are
$V_k$--equivalent for a given integer $k\ge 0$.  This algorithm also
works to determine whether or not two knots in $S^3$ are
${C_{k+1}}$--equivalent for $k\ge0$.
\end{remark}

\section{Examples and remarks}

In this section we give some examples of $C_k$--moves and also give
some remarks.

\subsection{Simple $C_k$--moves as band-sum operations}
As we have already seen, a simple $C_1$--move is equivalent to a
crossing change. It is also equivalent to band-summing a Hopf link
$L_2$, see Figure~\ref{exCk}a.  Hence any two knots in $S^3$ are
$C_1$--equivalent to each other.  On the other hand, any invariant of
knots in $S^3$ of type $0$ with values in any abelian group is a
constant function.

A simple $C_2$--move is equivalent to band-summing the Borromean rings
$L_3$, see Figure~\ref{exCk}b.  This operation has appeared in many places:
\cite{Murakami-Nakanishi:OnACertainMove},
\cite{Matveev:GeneralizedSurgeries},
\cite{Ohtsuki:FiniteTypeInvariants},
\cite{Garoufalidis-Ohtsuki:OnFiniteTypeIII}, etc.  H~Murakami and
Y~Nakanishi proved that any two knots in $S^3$ are related by a
sequence of operations of this kind, which they call
``$\Delta$--unknotting operations''
\cite{Murakami-Nakanishi:OnACertainMove}.  On the other hand, any
knot invariant of type $1$ with values in any
abelian group is again a constant function.

A simple $C_3$--move is equivalent to band-summing Milnor's link $L_4$
of $4$--component, see Figure~\ref{exCk}c.  As a corollary to
Theorem~\ref{JkCk}, we have the following result, which was originally
stated (in a slightly different form) and proved more directly in
\cite{Habiro:ClaspPassMoves}.

\begin{proposition}
\label{a2}
Two knots $\gamma $ and $\gamma '$ in $S^3$ are $C_3$--equivalent if
and only if $\gamma$ and $\gamma'$ has equal values of the Casson
invariant of knots, also known as the second coefficient in the
Alexander--Conway polynomial.  The group of $C_3$--equivalence classes
of knots in $S^3$ with multiplication induced by the connected sum
operation is isomorphic to $\z$.
\end{proposition}

\begin{proof}
This is clear from the fact that an invariant of type $2$ of knots in
$S^3$ is a linear combination of $1$ and the second coefficient of the
Alexander--Conway polynomial.
\end{proof}

More generally, a simple $C_k$--move ($k\ge 1$) is equivalent to
band-summing an iterated Bing double \cite{Cochran:Derivatives} of a
Hopf link with $k+1$ components.  The result of surgery on a simple
strict tree clasper $T$ of degree $k$ for a $(k+1)$--component unlink
$\gamma$ such that $\gamma$ bounds $k+1$ disjoint disks
$D_1,\dots,D_{k+1}$ in such a way that $D_i\cap T$ is an arc for
$i=1,\dots , k+1$ is an iterated Bing double of a Hopf link.  Iterated
Bing doubles are successfully used by T~Cochran
\cite{Cochran:Derivatives} to study the Milnor $\bar\mu$ invariants of
links.  It seems that claspers also work well in studying the Milnor
$\bar\mu$ invariants.  In the next subsection we give a few results
concerning the Milnor $\bar\mu$ invariants.

\begin{figure}
\cl{\includegraphics{exCk.eps}}
\nocolon\caption{}
\label{exCk}
\end{figure}

\subsection{$C_k$--equivalence and Milnor's $\bar\mu$ invariants}

For the definition of the Milnor $\mu$ and $\bar\mu$ invariants, see
\cite{Milnor:IsotopyOfLinks} or \cite{Cochran:Derivatives}.

\begin{theorem}
\label{mubar}
{\rm(1)}\qua  For $k,n\ge 1$, the Milnor $\mu$ invariants of length $k+1$ of
$n$--string links in $D^2\times I$ are invariants of
${C_{k+1}}$--equivalence.

{\rm(2)}\qua The Milnor $\bar\mu$ invariants of length $k+1$ of $n$--component
links in $S^3$ are invariants of ${C_{k+1}}$--equivalence.  (Recall that
each Milnor $\bar\mu$ invariant of length $k+1$ is only well-defined
modulo a certain integer determined by the Milnor $\bar\mu$ invariants
of length $\le k$.)
\end{theorem}

\begin{proof}
(1)\qua The Milnor $\mu$ invariants of length $k+1$ of string links are
invariants of type $k$ \cite{Bar-Natan:VassilievHomotopy}
\cite{Lin:PowerSeriesExpansions}, and hence invariants of
${C_{k+1}}$--equivalence by Corollary~\ref{CkJk}.

(2)\qua If two $n$--component links $\gamma$ and $\gamma'$ are
${C_{k+1}}$--equivalent, then they are equivalent to the closure of two
mutually $C_{k+1}$--equivalent $n$--string links $\gamma_1$ and
$\gamma_1'$ in $D^2\times[0,1]$.  By (1), $\gamma_1$ and $\gamma_1'$
have the same values of the Milnor $\mu$ invariants of length $\le
k+1$.  Hence $\gamma$ and $\gamma'$ have the same values of the Milnor
$\bar\mu$ invariants of length $\le k+1$.
\end{proof}

\begin{remark}
There is a more direct proof as follows.  We can prove that a
$C_k$--move on a link $\gamma$ preserves the $k$th nilpotent quotient
$(\pi_1E_\gamma)/(\pi_1E_\gamma)_{k+1}$ of the fundamental group of
the link exterior $E_\gamma$ of $\gamma$ in $M$ in a natural way.
(See also Section~\ref{ss:grope}.)  Theorem~\ref{mubar} follows
directly from this result.  We will give the details in a future
paper.
\end{remark}

By Theorem~\ref{JkCk}, the ${C_{k+1}}$--equivalence and the $V_k$--equivalence are
equal for knots in $S^3$.  For links in $S^3$ with more than $1$
component, we have the following.

\begin{figure}[ht!]
\cl{\includegraphics[width=4in]{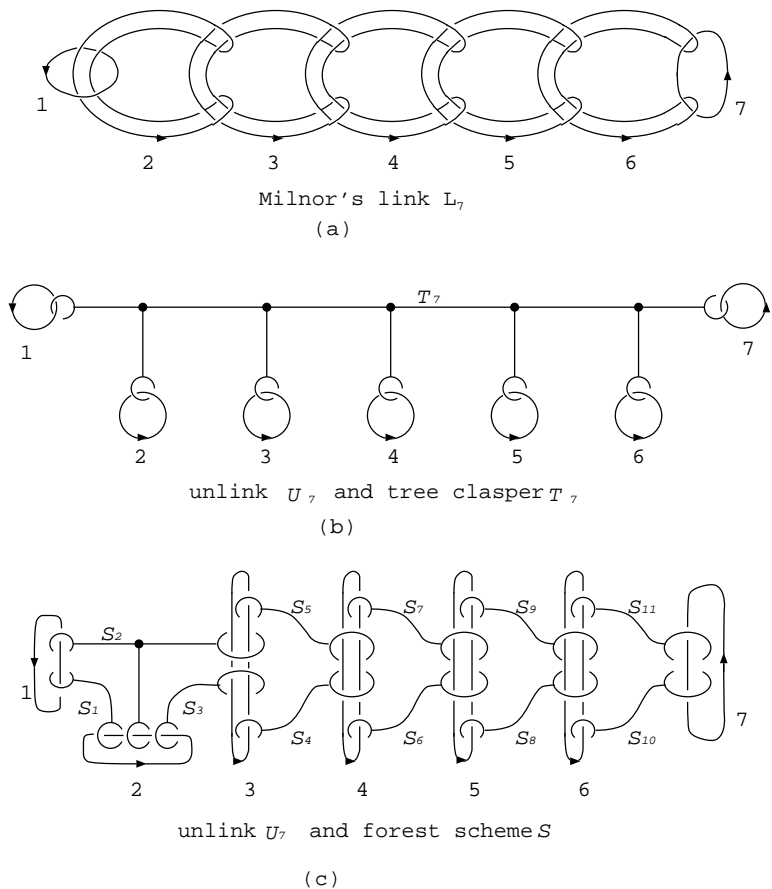}}
\caption{(a) Milnor's link $L_7$ of $7$--components ($k=6$).
(b)~Strict tree clasper $T_7$ for the unlink $U_7$.  (c) Strict forest
scheme $S= S_1\cup\dots\cup S_{11}$ for $U_7$.}
\label{Milnorlink1}
\end{figure}

\begin{proposition}
\label{Milnorlink}
For $k\ge 1$, let $U_{k+1}$ denote the $(k+1)$--component unlink and
let $L_{k+1}$ denote Milnor's link of $(k+1)$--components (which is a
($k+1$)--component iterated Bing double of a Hopf link), see
Figure~\ref{Milnorlink1}a.  Then we have the following.
\begin{enumerate}
\item $U_{k+1}$ and $L_{k+1}$ are ${C_k}$--equivalent but not
${C_{k+1}}$--equivalent. 
\item If $k=1$, then $U_2$ and $L_2$ are $V_0$--equivalent but not
$V_1$--equivalent.  If $k\ge 2$, then $U_{k+1}$ and $L_{k+1}$ are
$V_{2k-1}$--equivalent, but not $V_{2k}$--equivalent.
\end{enumerate}
\end{proposition}

\begin{proof}
We first prove 1.  That $U_{k+1}$ and $L_{k+1}$ are ${C_k}$--equivalent
follows from the fact that surgery on the simple strict tree clasper
$T_k$ of degree $k$ for $U_{k+1}$ as depicted in Figure~\ref{Milnorlink1}b
yields the link $L_{k+1}$.  That $U_{k+1}$ and $L_{k+1}$ are not
${C_{k+1}}$--equivalent follows from Theorem~\ref{mubar} and the fact that the link
$L_{k+1}$ has some non-vanishing Milnor $\bar\mu$ invariant of length
$k+1$ \cite{Milnor:LinkGroups}, but $U_{k+1}$ has vanishing
Milnor $\bar\mu$ invariants.

Now we prove 2.  If $k=1$, then $L_2$ is a Hopf link and the claim
clearly holds.  Assume that $k\ge 2$.  Let $S=\{S_1,\dots,S_{2k-1}\}$
be the forest scheme of degree $2k$ for $U_{k+1}$ as depicted in
Figure~\ref{Milnorlink1}c.  Then it is not difficult to prove that
$[U_{k+1};S_1,\dots ,S_{2k-1}]=[L_{k+1}]-[U_{k+1}]$.  (The proof goes
as follows: 
\[\begin{split}
 &[U_{k+1};S_1,\dots,S_{2k-1}] =-[U_{k+1};S_2,\dots ,
S_{2k-1}] =[U_{k+1};S_2,S_4,\dots,S_{2k-1}]\\
=&[U_{k+1}^{S_{2k-1}};S_2,S_4,\dots,S_{2k-2}]
=[U_{k+1}^{S_{2k-2}\cup{}S_{2k-1}};S_2,S_4,\dots,S_{2k-3}] \\
=&\dots=[U_{k+1}^{S_4\cup\cdots\cup S_{2k-1}};S_2] =[L_{k+1}]-[U_{k+1}] .
\end{split}\]
The details are left to the reader.)  Hence $L_{k+1}$ and $U_{k+1}$
are $V_{2k-1}$--equivalent.  That $L_{k+1}$ and $U_{k+1}$ are not
$V_{2k}$--equivalent can be verified, for example, by calculating the
linear combination of uni-trivalent graphs of degree $2k$
corresponding to the difference $L_{k+1}-U_{k+1}$ and taking the value
of it in, say, the $sl_2$--weight system (but not in the
Alexander--Conway weight system).
\end{proof}		

\begin{remark}
We can generalize a part of Proposition~\ref{Milnorlink} that
$L_{k+1}$ is both $C_k$--equivalent and $V_{2k-1}$--equivalent to $U_k$
for $k\ge 2$ as follows: If a $(k+1)$--component link $\gamma$ in $S^3$
is Brunnian (ie, every proper sublink of $\gamma $ is an unlink),
then $L$ is both $C_k$--equivalent and $V_{2k-1}$--equivalent to the
$(k+1)$--component unlink $U_{k+1}$.  We will prove this result in a
future paper.
\end{remark}

\section{Surveys on some other aspects of the calculus of claspers}

In this section we survey some applications of claspers to other field
of $3$--dimensional topology.  We will prove the results below in
forthcoming papers.

\subsection{Calculus of claspers and commutator calculus in braided category}
\label{ss:BraidedCommutatorCalculus}

\begin{figure}
\cl{\includegraphics[width=4in]{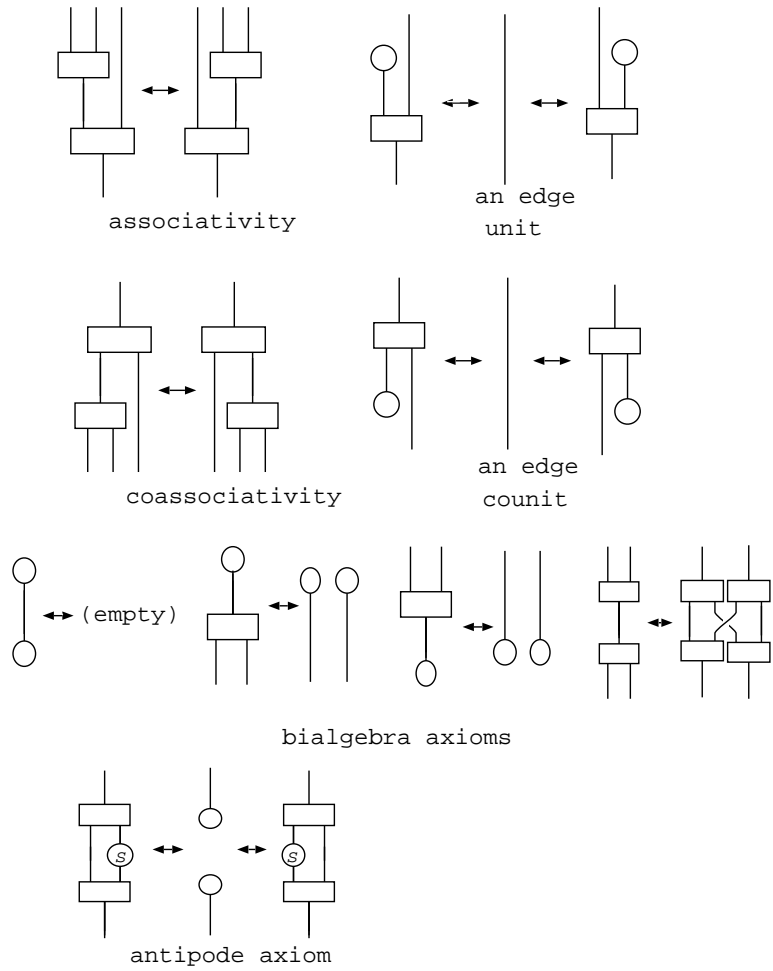}}
\caption{Claspers satisfy the axioms of Hopf algebra in a braided
category.}
\label{Hopf1}
\end{figure}

\begin{figure}
\cl{\includegraphics{PfAssoc.eps}}
\nocolon\caption{}
\label{PfAssoc}
\end{figure}
                  
The reader may have noticed that some of the moves introduced in
Proposition~\ref{moves} are similar to the axioms of a Hopf algebra in
a braided category.  To see this, we think of an edge as a Hopf
algebra, a box as a (co)multiplication, a trivial leaf as a (co)unit
and a positive half-twist as an antipode.  Then move~{3} corresponds
to the axiom of (co)unit, and move~{4} to that of antipode.  Other
axioms actually hold as illustrated in Figure~\ref{Hopf1}.  For the
proof of the ``associativity'' and the fourth of the ``bialgebra
axioms,'' see Figures~\ref{PfAssoc} and \ref{PfBialg}, respectively.
The proofs of the others are easy.  This Hopf algebra structure in
claspers is closely related to the Hopf algebra structure in
categories of cobordisms of surfaces with connected boundary by
L~Crane and D~Yetter \cite{Crane-Yetter:OnAlgebraicStructures} and
by T~Kerler \cite{Kerler:Genealogy}.

\begin{figure}
\cl{\includegraphics{PfBialg.eps}}
\nocolon\caption{}
\label{PfBialg}
\end{figure}

\begin{figure}
\cl{\includegraphics{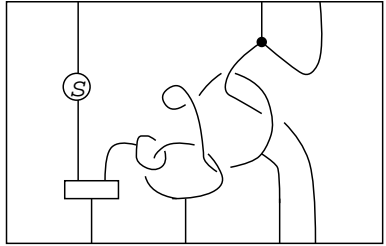}}
\caption{A clasper diagram representing a morphism from
$3$ to $4$ in the category ${\mathbf{Cl}}_0$.}
\label{ClasperDiagram}
\end{figure}

Let us give a rough definition of the braided category in which
claspers live.  A {\em clasper diagram} will mean a picture of a
clasper drawn in a square $[0,1]^2$ with some edges going out of the
top and the bottom edges of $[0,1]^2$, see for example
Figure~\ref{ClasperDiagram}.  Two clasper diagrams $D$ and $D'$ are
said to be equivalent if the numbers of edges of $D$ and $D'$ on the
top (resp.\ bottom) are equal and they represent two surfaces equipped
with decompositions in $[0,1]^3$ that are ambient isotopic to each
other relative to boundary of the cube $[0,1]^3$ (after a suitable
reparameterization near the top and the bottom squares).  Then the
category ${\mathbf{Cl}}_0$ of clasper diagrams is defined as follows.
The objects of ${\mathbf{Cl}}_0$ are nonnegative integers.  The
morphisms from $m$ to $n$ in ${\mathbf{Cl}}_0$ are equivalence classes
of clasper diagrams with $m$ edges on the top and $n$ edges on the
bottom.  The composition is induced by pasting two diagrams
vertically.  Identity $1_m\colon m\to m$ is the equivalence class of
the diagram consisting of $m$ vertical edges.  The tensor functor
$\otimes\colon {\mathbf{Cl}}_0\times{\mathbf{Cl}}_0\to{\mathbf{Cl}}_0$
is induced by addition of integers and placing two diagrams
horizontally.  The monoidal unit $I$ is $0$.  The braiding
$\Psi_{m,n}\colon m\otimes n\to n\otimes m$ is a positive crossing of
two parallel families of edges.

Let us also give a sketch of the definition of the category
${\mathbf{Cob}}$ of cobordisms of oriented connected surfaces with
connected boundary.  For a precise definition, see
\cite{Crane-Yetter:OnAlgebraicStructures} or \cite{Kerler:Genealogy}.
The objects in ${\mathbf{Cob}}$ are nonnegative integers.  For each
object $m$ in ${\mathbf{Cob}}$, we fix a surface $F_m$ of genus $m$
with one boundary component.  We assume that $F_0=[0,1]^2$, the
surface $F_1$ is a ``square with a handle,'' and $F_m$ with $m\ge2$
are obtained by pasting $m$ copies of $F_1$ side by side, see
Figure~\ref{F0F1F2}.  For $m\ge0$, the boundary of $F_m$ is
parameterized by $\partial([0,1]^2)$ in a natural way.  A cobordism
from $F_m$ to $F_n$ is a $3$--manifold with boundary parameterized by
the surface $(-F_m)\cup_{\partial([0,1]^2)\times\{0\}}
(\partial([0,1]^2)\times[0,1])\cup_{\partial([0,1]^2)\times\{1\}}F_n$, where
$-F_m$ is $F_m$ with orientation reversed.  The morphisms from $m$ to
$n$ are the diffeomorphism classes, respecting boundary
parameterizations, of cobordisms from $F_m$ to $F_n$.  The composition
in ${\mathbf{Cob}}$ is induced by ``pasting the bottom surface of one
cobordism with the top surface of another.''  The identity $1_m\colon
m\to m$ is the direct product $F_m\times[0,1]$ with the obvious
boundary parameterization.  The tensor functor is induced by addition
of integers and ``pasting two cobordisms side by side.''  The monoidal
unit in ${\mathbf{Cob}}$ is $0$.  (We identify the boundary connected
sum of $F_m$ and $F_n$ with $F_{m+n}$ via a certain predescribed
diffeomorphism.)  The braiding is obtained by ``letting two identity
cobordisms cross each other positively.''

\begin{figure}
\cl{\includegraphics{F0F1F2.eps}}
\nocolon\caption{}
\label{F0F1F2}
\end{figure}

Then the object $1$ in ${\mathbf{Cob}}$, which will be denoted by $H$,
has a Hopf algebra structure \cite{Crane-Yetter:OnAlgebraicStructures},
\cite{Kerler:Genealogy}.

\begin{figure}
\cl{\includegraphics{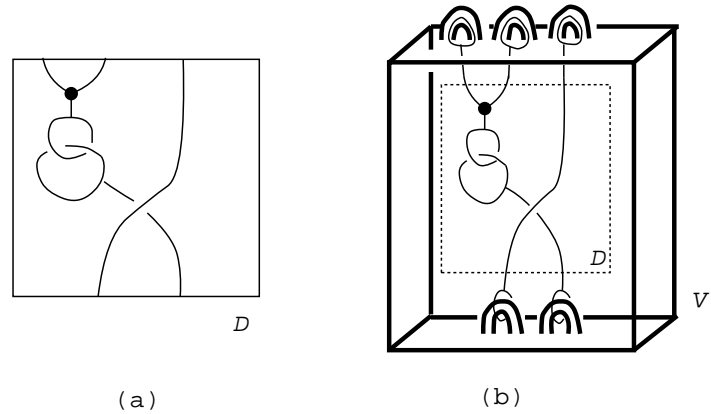}}
\caption{In (a) is a clasper diagram $D$ representing a
morphism $[D]$ from $2$ to $3$ in ${\mathbf{Cl}}_0$.  We embed it in a
``cube-with-handles-and-holes'' $V$ as depicted in (b) together with
some extra leaves running through the handles or linking with the
holes.  Let $G_D$ denote the clasper obtained in this way.  The image
of $[D]$ by $F$ is represented by the result of surgery from $V$ on
$G_D$.}
\label{CobordismDiagram}
\end{figure}

We define a functor $F\colon {\mathbf{Cl}}_0\to{\mathbf{Cob}}$
respecting the structure of braided strict monoidal categories.  On
the object level, $F$ maps a nonnegative integer $n$ into $n$.  On the
morphism level, $F$ maps a morphism in ${\mathbf{Cl}}_0$ into one in
${\mathbf{Cob}}$ as illustrated in Figure~\ref{CobordismDiagram}.
It is not difficult to see that $F$ is a functor and respects the
structure of braided strict monoidal category.

The relations among clasper diagrams depicted in Figure~\ref{Hopf1}
implies that there is a Hopf algebra structure on $H$ in
${\mathbf{Cob}}$.  We can check that this Hopf algebra structure is
essentially equivalent to that given in
\cite{Crane-Yetter:OnAlgebraicStructures} and \cite{Kerler:Genealogy}.
Thus clasper diagrams provides a new way to visualize the cobordisms
of surfaces.  This may be regarded as a variant of a similar
visualization of cobordisms using ``bridged links'' due to Kerler
\cite{Kerler:BridgedLinksAndTanglePresentations}.

Let ${\mathbf{Cl}}$ denote the coimage of the functor $F$, ie, the
category obtained from ${\mathbf{Cl}}$ by regarding each two morphisms
mapped by $F$ into equal morphism to be equal.  Of course,
${\mathbf{Cl}}$ is isomorphic to the image of $F$.  It is easy to
check that $F$ is surjective, and hence ${\mathbf{Cl}}$ is isomorphic
to ${\mathbf{Cob}}$.

\begin{figure}
\cl{\includegraphics{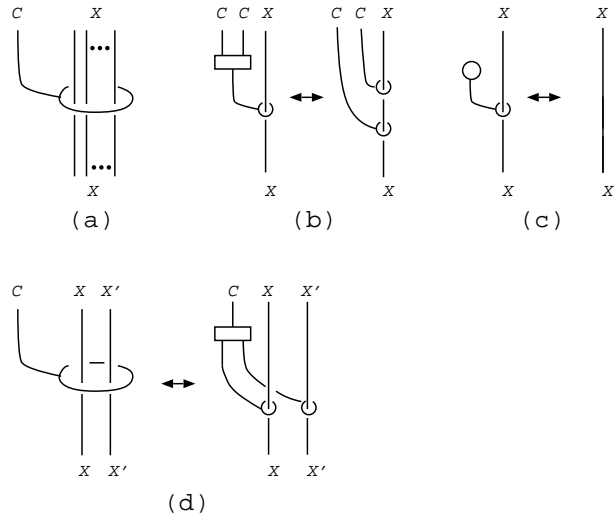}}
\caption{In (b), (c) and (d), the arcs labeled $X$ and $X'$
represents parallel families of edges and strings (but not leaves).
(a)  Action of $C$ on $X$.  (b) Associativity.  (c)
Unitality.  (d) Action on tensor product.}
\label{action}
\end{figure}

Now we give an interpretation of disk-leaves and leaves as {\em
actions} of the Hopf algebra $C$ on other objects.  For this, we
extend the notions of clasper diagrams and cobordisms to those involving
links and enlarge the categories ${\mathbf{Cl}}$ and ${\mathbf{Cob}}$
to ${\mathbf{Cl}}'$ and ${\mathbf{Cob}}'$, respectively.  Then we may
think of a leaf bounding an embedded disk as a left action of the Hopf
algebra on an object, see Figure~\ref{action}.  The ``associativity''
(b) is equivalent to move~{6} and the ``unitality'' (c) is a
consequence of move~{2}.  Figure~\ref{action}d is equivalent to
move~{8} and shows how the Hopf algebra $C$ acts on the tensor product
(ie, parallel) of two objects $X$ and $X'$.  Because of the
obvious self duality of the Hopf algebra $C$ in ${\mathbf{Cl}}$, we
may think of (disk-)leaves also as {\em coactions}.

Now we give an interpretation of nodes as {\em (co)commutators}.  See
Figure~\ref{8:commutator}.  We can transform a clasper diagram consisting of
a node on the left side to the clasper diagram on the right side.
Here the box with many input edges replaces as depicted in
Figure~\ref{manyinputbox}.  We explain how we can think of the right side as
a commutator.  Recall that one of the most typical examples of Hopf
algebras is the {\em group Hopf algebra} $kG$ of a group $G$ with $k$
a field, where the algebra structure is induced by the group
multiplication, the coalgebra structure is given by
$\Delta(g)=g\otimes g$ and $\epsilon(g)=1$ for $g\in G$, and the
antipode is given by $S(g)=g^{-1}$ for $g\in G$.  So, we try to input
two group elements $a$ and $b$ into the two top edges and see what we
obtain as the output from the bottom edge.  We think of the two upper
boxes as comultiplications, which duplicate $a$ and $b$.  We think of
the symbols `$S$' as antipodes, which invert the elements $a$ and $b$
in the middle.  The braiding permutes $a^{-1}$ and $b^{-1}$.  The
lower box acts as a multiplication map and multiplies $a$, $b^{-1}$,
$a^{-1}$ and $b$.  Hence we obtain a commutator $ab^{-1}a^{-1}b$ as
the output.  This explains why we think of the left side as a
commutator.  In the third in Figure~\ref{8:commutator} we consider the
fundamental group of the complement of two upper leaves and incident
half-edges in $[0,1]^3$, which is a free group of rank $2$ freely
generated by the meridians to the two leaves, $a$ and $b$.  Then the
element of this free group represented by a boundary component of the
lower leaf is again the commutator $ab^{-1}a^{-1}b$.

\begin{figure}
\cl{\includegraphics{8comm.eps}}
\nocolon\caption{}
\label{8:commutator}
\end{figure}

\begin{figure}
\cl{\includegraphics{manybox.eps}}
\nocolon\caption{}
\label{manyinputbox}
\end{figure}

In this group theoretic analogy, a tree clasper can be thought of as
an iterated commutator.  We can give group theoretic interpretation to
some of the results in the previous sections.  For example,
Proposition~\ref{trivialleaf} is similar to the fact that an iterated
commutators of group elements is $1$ if at least one on the elements
is $1$, Proposition~\ref{slidediskleaf} is similar to the fact that
``two iterated commutators of class $k$ and $k'$ commutes each other
up to an iterated commutator of class $k+k'$,'' and so on.  These
interpretations greatly help us understand the algebraic nature of
calculus of claspers and theory of finite type invariants.

However, this group theoretic analogy does not work very well in some
cases.  For example, let $\beta^*\colon C\to C\otimes C$ be the dual
(ie, rotation by $\pi$) of the commutator $\beta\colon C\otimes C\to
C$.  We call $\beta^*$ a {\em cocommutator}.  The cocommutator
$\beta^*$ replaces the dual of the last diagram in
Figure~\ref{8:commutator}.  It is easy to check that inputting any
group element $g$ on the top of this diagram yields $1\otimes 1$ as
output.  So, the group theoretic analogy or more generally, such an
analogy involving cocommutative Hopf algebras does not work well for
cocommutators.

Therefore to understand the algebraic nature of claspers more
accurately, we must seek such analogy for more general Hopf algebras
in braided categories.  This leads us to {\em commutator calculus of
Hopf algebras in braided categories}, or {\em braided commutator
calculus}, which may be regarded as a branch of ``braided
mathematics'' proposed by S. Majid.  Let us briefly explain
commutators and cocommutators appearing in this new commutator
calculus here.

Let ${\mathcal B}$ be a braided strict monoidal category and let
$H=(H,\mu,u,\Delta,\epsilon,S)$ be a Hopf algebra in ${\mathcal B}$.
Then we define the {\em commutator} $\beta\colon H\otimes H\to H$ via
Figure~\ref{8:commutator}.  ie, we set
\begin{equation}
\label{eq:beta}
\beta=\mu_4(H\otimes\Psi_{H,H}\otimes H)(H\otimes S\otimes S\otimes H) 
(\Delta\otimes\Delta),
\end{equation}
where $\mu_4\colon H^{\otimes4}\to H$ is the multiplication with four
inputs, and $\Psi_{H,H}$ is the braiding of $H$ and $H$.  Dually we
define the {\em cocommutator} $\gamma\colon H\to H\otimes H$ by
\begin{equation}
\label{eq:betastar}
\gamma=(\mu\otimes\mu)(H\otimes S\otimes S\otimes H)
(H\otimes\Psi_{H,H}\otimes H)\Delta_4,
\end{equation}
where $\Delta_4\colon H\to H^{\otimes4}$ is the comultiplication with
four outputs.  It seems that commutator calculus based on these
(co)commutators works well at least when $H$ is ``braided
cocommutative with respect to the adjoint action'' in S~Majid's sense
\cite{Majid:AlgebrasAndHopfAlgebras}.  This braided cocommutativity is
satisfied by the Hopf algebra $C$ in ${\mathbf{Cl}}$ and hence by $H$
in ${\mathbf{Cob}}$.  In this abstract setup, for example, variants of
some of the moves in Proposition~\ref{moves} holds, and zip
construction works.  Commutator calculus in braided category will
enable us to handle complicated lemmas on claspers purely
algebraically, and moreover help us formalize a large part of calculus
of claspers in the language of category theory.

\subsection{Graph claspers as topological realization of uni-trivalent graphs} 
The notion of tree claspers is generalized to that of graph claspers.
We may regard graph claspers as ``topological realizations'' of
uni-trivalent graphs that appear in theories of finite type invariant
of links and $3$--manifolds.

A {\em graph clasper} $G$ for a link $\gamma$ in $M$ is a clasper
consisting only of leaves, disk-leaves, nodes and edges.  $G$ is {\em
admissible} if each component of $G$ has at least one disk-leaf, and
is {\em strict} if, moreover, $G$ has no leaves.  $G$ is {\em simple}
if every disk-leaf of $G$ intersects the link with one point.  The
{\em degree} of connected strict graph clasper $G$ is half the number
of disk-leaves and nodes of $G$, and the degree of a general strict
graph clasper $G$ is the minimum of the degrees of components of $G$.

A {\em graph scheme} $S=\{S_1,\dots,S_l \}$ is a scheme consisting
of connected graph claspers $S_1,\dots,S_l $.  $S$ is strict
(resp.\ admissible, simple) if every element of $S$ is strict
(resp.\ admissible, simple).  The {\em degree} $\deg S$ of a strict
graph scheme $S$ is the sum of the degrees of its elements.

We can generalize a large part of definitions and results in previous
sections using graph claspers.  For example, we can prove that two
links related by a surgery on a strict graph clasper for a link
$\gamma$ of degree $k\ge1$ are ${C_k}$--equivalent.  So we may redefine
the notion of $C_k$--equivalence using strict graph claspers.  We can
also prove that, for a link ${\gamma_0}$ in $M$, the subgroup
$J_k({\gamma_0})$ of $\z\calL_1({\gamma_0})$ equals the subgroup
generated by the elements $[\gamma,S]$, where $\gamma$ is a link in
$M$ which is $C_1$--equivalent to ${\gamma_0}$ and $S$ is a strict
graph scheme for ${\gamma_0}$ of degree $k$.

We can generalize the definitions and results in Section~4 to simple
strict graph claspers.  For a link ${\gamma_0}$ in $M$, let
${\tilde{\mathcal G}}_k^h(\gamma_0)$ denote the free abelian group
defined similarly as $\tilde{\mathcal F}_k^h({\gamma_0})$ but we use
simple strict graph claspers instead of simple forest graph claspers.
Let $R_k$ denote the subgroup of ${\tilde{\mathcal
G}}_k^h({\gamma_0})$ generated by the elements depicted in
Figure~\ref{ASIHXSTU}.  They are called antisymmetry relations, IHX
relations and STU relations.\footnote{The sign of the last term in the
STU relation looks different from the usual one for a technical
reason.}  Here we allow only STU relations of a special kind which
involves only {\em connected} graph claspers.  We can prove that the
natural map $\nu_k\colon {\tilde{\mathcal
G}}_k^h({\gamma_0})\to{{\bar\calL}_k}({\gamma_0})$ which exists by an
analogue of Theorem~\ref{FactorThroughFhkg} factors through
${\tilde{\mathcal G}}_k^h({\gamma_0})/R_k$.

\begin{figure}
\cl{\includegraphics{ASIHXSTU.eps}}
\nocolon\caption{}
\label{ASIHXSTU}
\end{figure}

A {\em uni-trivalent graph} $D$ on a $1$--manifold $\alpha$ is an
abstract finite graph $D$ possibly with loop edges and multiple edges
such that every vertex of $D$ is of valence $1$ or $3$, to each
trivalent vertex of $D$ is equipped with a cyclic order on the three
incident edges, and to some of the univalent vertices of $D$ are
equipped with points on $\alpha$.  Here two distinct vertex must
corresponds to distinct points.  We call the univalent vertices of $D$
equipped with points in $\alpha$ the univalent vertices {\em on
$\alpha$}.  A uni-trivalent graph $D$ on a $1$--manifold $\alpha$ is
{\em strict} if every univalent vertex is on $\alpha$ and if each
connected component of $D$ have at least one univalent vertex.  The
{\em degree} of a strict uni-trivalent graph $D$ is half the number of
vertices of $D$.

In the following we restrict our attention to links in $S^3$ and
string links in $D^2\times[0,1]$ for simplicity.  Let ${\gamma_0}$ be
an unlink or a trivial string link.  We here refer to links of the
same pattern as ${\gamma_0}$ simply as ``links.''

For $k\ge0$, let ${\mathcal A}_k({\gamma_0})$ denote the abelian group
generated by strict uni-trivalent graphs of degree $k$ on
${\gamma_0}$, subject to the framing independence relations and the
(usual) STU relations (and hence subject to the antisymmetry and IHX
relations).  See \cite{Bar-Natan:OnTheVassiliev} for the definitions
of these relations.  We set
${\bar{J}_k}({\gamma_0})={J_k}({\gamma_0})/{J_{k+1}}({\gamma_0})$.
Let $\xi_k\colon {\mathcal A}_k({\gamma_0})\to{\bar{J}_k}({\gamma_0})$
denote a well-known surjective homomorphism which ``replaces chords
with double points''.\footnote{In a previous version, $\xi_k$ was
claimed to be an isomorphism, but it does not seem to be known whether
this is an isomorphism.  However, $\xi_k\otimes{\mathbf
Q}\colon{\mathcal A}_k({\gamma_0})\otimes{\mathbf
Q}\to{\bar{J}_k}({\gamma_0})\otimes{\mathbf Q}$ is injective and hence
an isomorphism by Kontsevich's theorem.}
Let $\iota_k\colon {\tilde{\mathcal G}}_k^h({\gamma_0})\to{\mathcal
A}_k({\gamma_0})$ denote the natural homomorphism which maps a class
of a connected simple strict graph clasper $G$ for ${\gamma_0}$ of
degree $k$ into the ``corresponding strict uni-trivalent graph'' of
$G$ with an appropriate sign.  See Figure~\ref{correspondingUTG}.  Let
$\chi_k\colon {{\bar\calL}_k}({\gamma_0})\to{\bar{J}_k}({\gamma_0})$
be the homomorphism defined by
$\chi_k([\gamma]_{C_{k+1}})=[\gamma-{\gamma_0}]_{{J_{k+1}}({\gamma_0})}$.
Then we can prove that the following diagram commutes (up to
sign).\footnote{In the case of links in $S^3$ with more than one
component, the map $\chi_k$ is not injective in general.
Conjecture~\ref{conj:ckkVk} for string links in $D^2\times[0,1]$ is
equivalent to that $\chi_k$ is injective for all $k,n\ge0$.  Hence,
for string knots and knots in $S^3$, $\chi_k$ is injective.}
\begin{equation}
\label{eq:diag}
\begin{CD}
{\tilde{\mathcal G}}_k^h({\gamma_0}) @>{\iota_k}>> {\mathcal
A}_k({\gamma_0}) \\ @V{\nu_k}VV @VV{\xi_k}V \\
{{\bar\calL}_k}({\gamma_0}) @>>{\chi_k}> {\bar{J}_k}({\gamma_0})
\end{CD}
\end{equation}

\begin{figure}
\cl{\includegraphics{corrUTG.eps}}
\nocolon\caption{}
\label{correspondingUTG}
\end{figure}

From these results, we may think of graph claspers as {\em topological
realizations of strict uni-trivalent graphs}.  In other words, any
primitive strict uni-trivalent graph, $D$, of degree $k$ on a
${\gamma_0}$ is ``realized'' by the knot obtained from the trivial
knot by surgery on the simple strict graph clasper $G_D$ such that the
``corresponding strict uni-trivalent graph'' of $G_D$ is $D$.  Related
realization results of uni-trivalent graphs are given by K\,Y~Ng
\cite{Ng:GroupsOfRibbonKnots} and by N~Habegger and G~Masbaum
\cite{Habegger-Masbaum:TheKontsevichIntegral}.  One of the advantages
of using graph claspers is that for any connected strict uni-trivalent
graph, $D$, we can immediately find a simple strict graph clasper
realizing $D$.

From the category-theoretical point of view described in 8.1, it is
important to note that {\em the Lie algebraic structures appearing in
theories of finite type invariants of links and $3$--manifolds
originate from the Hopf algebraic structure in the category
${\mathbf{Cl}}\cong{\mathbf{Cob}}$} (or in a suitably extended
category involving links).  This is just like that commutator calculus
in the associated graded Lie algebra of the lower central series of a
group can be explained in terms of commutator calculus in the group.


\subsection{New filtrations and equivalence relations on links based
on admissible graph claspers} 
\label{NewFiltration}
Using admissible graph claspers, we can define a new filtration on
links which is much coarser than the Vassiliev--Goussarov filtration
and from this filtration we can define a special class of finite type
invariants of links in $3$--manifolds.

For a connected graph clasper $G$, the {\em $A$--degree} of $G$ is the
number of disk-leaves and nodes of $G$, and the {\em $S$--degree} of
$G$ is ${\frac{1}{2}}(\mbox{$A$--deg } G - l(G))$, where $l(G)$ is the
number of leaves of $G$.  For a general graph clasper $G$, the
$A$--degree (resp.\ $S$--degree) of $G$ is the minimum of the $A$--degrees
(resp.\ $S$--degrees) of components of $G$.  Observe that the degree of
a strict graph clasper equals the $S$--degree.  We define the
$A$--degree (resp.\ $S$--degree) of a graph scheme $S$ to be the sum of
$A$--degrees (resp.\ $S$--degrees) of elements of $S$.

A uni-trivalent graph $D$ on a $1$--manifold $\alpha$ is said to be
{\em $H$--labeled}, where $H$ is an abelian group, if each univalent
vertex of $D$ that is not on $\alpha$ is labeled by an element of $H$.
An $H$--labeled uni-trivalent graph $D$ on $\alpha$ is {\em admissible}
if every component of $D$ has at least one univalent vertex on
$\alpha$.  For such $D$, the $A$--degree of $D$ is the sum of the
number of the trivalent vertices and the number of the univalent
vertices on $\alpha$, and the $S$--degree of $D$ is half the difference
of $A$--degree of $D$ and the number of univalent vertices of $D$ not
on $\alpha$.

Let $P=(\alpha,i)$ be a pattern on a $3$--manifold $M$.  Let
${\calL(P)}$ denote the set of equivalence classes of links in $M$ of
pattern $P$.  For $k\ge0$, Let $J^A_k(P)$ denote the subgroup of
$\z{\calL(P)}$ generated by all the elements $[\gamma,S]$, where
$\gamma$ is a link in $M$ of pattern $P$ and $S$ is an admissible
graph scheme for $\gamma$ of $A$--degree $k$.

For each $k,l\ge0$ with $0\le 2l\le k$, we define an abelian group
${\mathcal A}^A_{k,l}(P)$ to be generated by admissible
$H_1(M;\z)$--labeled uni-trivalent graphs of $A$--degree $k$ and of
$S$--degree $l$, on the $1$--manifold $\alpha$ whose univalent vertices
are labeled elements of the first homology group $H_1(M;\z)$ and to be
subject to the framing independence, antisymmetry, IHX, STU relations
and multilinearity of labels.

We set $\bar{J}^A_k(P)=J^A_k(P)/J^A_{k+1}(P)$.  For $l$ with $0\le
2l\le k$, let $\bar{J}^A_{k,l}(P)$ denote the subgroup of
$\bar{J}^A_k(P)$ generated by the elements $[\gamma,S] \mod
J^A_{k+1}(P)$, where $S$ is an admissible graph schemes for $\gamma$ of
$A$--degree $k$ and of $S$--degree $\ge l$.

\begin{figure}
\cl{\includegraphics{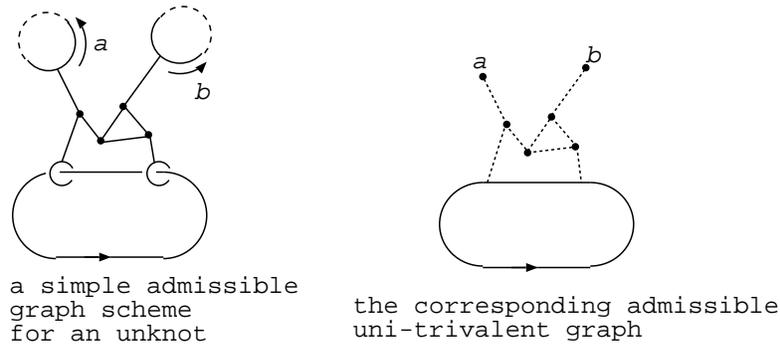}}
\caption{To the simple admissible graph scheme on the
left corresponds the admissible uni-trivalent graph on the right that
is labeled elements $a,b$ in $H_1(M;\z)$.}
\label{corresponding2}
\end{figure}

We can define a natural surjective homomorphism of ${\mathcal
A}^A_{k,l}(P)$ onto the graded quotient
$\bar{J}^A_{k,l}(P)/\bar{J}^A_{k,l+1}$.  By this homomorphism an
admissible uni-trivalent graph $D$ is mapped into the element
$[\gamma, S_D] \mod J^A_{k+1}(P)$, where $\gamma$ is a link of pattern
$P$ and $S_D$ is a simple admissible graph scheme whose
``corresponding admissible uni-trivalent graph'' is $D$, see
Figure~\ref{corresponding2}. Here the homology class of each leaf of
$S_D$ equals the label of the corresponding univalent vertex of $D$.

Since each ${\mathcal A}^A_{k,l}(P)$ is finitely generated,
$\bar{J}^A_{k,l}(P)/\bar{J}^A_{k,l+1}$ is finitely generated.  Hence
$\bar{J}^A_k(P)$ and $\z{\calL(P)}/J^A_{k+1}(P)$ are finitely
generated.

A homomorphism $f\colon \z{\calL(P)}\to X$, where $X$ is an abelian
group, is said to be {\em of $A$--type $k$} if $f$ vanishes on
$J^A_{k+1}(P)$.  We call the homomorphism $f\colon \z{\calL(P)}\to A$
of finite $A$--type a {\em special finite type invariant} since we can
prove that an invariant of $A$--type $2k$ is an invariant of type $k$.
For links in $S^3$ and string links in ${D^2\times[0,1]}$, we can
prove that the notions of $A$--type $2k$ and that of type $k$ are
equivalent.  We can prove that, for an integral homology $3$--sphere
$M$, there is an isomorphism between the group of invariants of type
$k$ for links in $M$ and that of $S^3$.  This implies that any finite
type invariant of links in $S^3$ canonically extends to a special
finite type invariant of links in integral homology $3$--spheres.  This
enables us to extend in a natural way the polynomial invariants such
as the Jones, HOMFLY and Kauffman polynomials to links in integral
homology spheres.

For $k\ge 1$, an {\em $A_k$--move} on a link is defined to be a surgery
on an admissible graph clasper of $A$--degree $k$.  It is clear that an
$A_{k+1}$--move preserves any invariant of $A$--type $k$.  The notion of
$A_k$--equivalence is defined in the obvious way.  The set of
$A_k$--equivalence classes of $n$--string links in
${\Sigma\times[0,1]}$, where $\Sigma$ is a connected oriented surface,
form a finitely generated nilpotent group (cf. Theorem~\ref{group}).
We can define the associated graded Lie algebra, say
$\bar\calL^A(\Sigma,n)$, where the $k$th homogeneous part is the
abelian group of $A_{k+1}$--equivalence classes of $n$--string links in
${\Sigma\times[0,1]}$ which are $A_k$--equivalent to the trivial
$n$--string link ${1_n}$.  This group is finitely generated.  This is a
(much more tractable) quotient of the Lie algebra
$\bigoplus_{k=1}^\infty{{\hat\calL}_k{(\Sigma,n)}}
/{\hat\calL}_{k+1}{(\Sigma,n)}$ defined in Section~5.  We can
completely determine the structure of the Lie algebra
$\bar\calL^A(\Sigma,n)\otimes\q$ using admissible
$H_1(\Sigma;\z)$--labeled uni-trivalent graphs, at least when $\Sigma$
is not closed.  In the proof we require the Le--Murakami--Ohtsuki
invariant \cite{Le-Murakami-Ohtsuki:OnAUniversal}.

To some extent the definition of the $A$--filtration resembles
M~Goussarov's filtration using ``interdependent modifications''
\cite{Goussarov:InterdependentModifications}.  If we use only
admissible graph claspers of a special kind such that all the leaves
are {\em null-homotopic} in the $3$--manifold, then we obtain a theory
equivalent to Goussarov's.

\subsection{Clasper surgeries and finite type invariants of $3$--manifolds}

Theories of clasper surgeries and finite type invariants of links in a
fixed $3$--manifold developed in previous sections are naturally
generalized to that of ($3$--manifold, link) pairs by allowing graph
claspers that are not necessarily tame.  These theories are very
closely related to known theories of finite type invariants and
surgery equivalence relations of $3$--manifolds
\cite{Ohtsuki:FiniteTypeInvariants}
\cite{Le-Murakami-Ohtsuki:OnAUniversal}
\cite{Garoufalidis-Levine:FiniteTypeBlinks}
\cite{Cochran-Melvin:FiniteTypeInvariants}
\cite{Cochran-Gerges-Orr:SurgeryEquivalenceRelations}.

After almost finished this paper, the author received a paper of
M. Goussarov \cite{Goussarov:NewTheoryOfInvariants}.  It seems that
some results in this section overlap that in
\cite{Goussarov:NewTheoryOfInvariants}.

%
\subsubsection{$A_k$--surgery equivalence relations}
For simplicity, we consider only compact connected closed
$3$--manifolds without links, though we can naturally generalize a
large part of the following definitions and results to $3$--manifolds
with boundaries and $3$--manifolds with links.

A graph clasper $G$ (for the empty link) in a $3$--manifold $M$ is {\em
allowable} if every component of $G$ is not a basic clasper.  Note
that every component of an allowable $G$ has no disk-leaf and has at
least one node. The $A$--degree of a connected graph clasper $G$ is
equal to the number of nodes in $G$, and the $S$--degree of $G$ is
equal to half the number of nodes minus half the number of leaves.
For a connected allowable graph clasper $G$, we have
${\text{$A$--deg}\thinspace} G\ge1$ and ${\text{$S$--deg}\thinspace}
G\ge -1$.  For $k\ge 1$, an {\em $A_k$--surgery} is defined to be a
surgery on a connected allowable graph clasper of $A$--degree $k$.  We
define the notion of {\em $A_k$--surgery equivalence} as the
equivalence relation on closed $3$--manifolds generated by
$A_k$--surgeries and orientation-preserving diffeomorphisms.

It turns out that two $3$--manifolds $M$ and $M'$ are $A_k$--surgery
equivalent if and only if there is a connected compact oriented
surface $F$ embedded in $M$ (which may be closed or not) and an
element $\alpha$ of the $k$th lower central series subgroup of the
Torelli group of $F$ such that the $3$--manifold obtained from $M$ by
cutting $M$ along $F$ and reglueing it using the self-diffeomorphism
of $F$ representing $\alpha$ is diffeomorphic to $M'$.  Such
modifications of $3$--manifolds by elements of the Torelli groups
appear in \cite{Morita:CassonsInvariant}
\cite{Garoufalidis-Levine:FiniteTypeBlinks} for integral homology
$3$--spheres.

A result of S\,V Matveev is restated that two closed $3$--manifolds
$M$ and $M'$ are $A_1$--surgery equivalent if and only if there is an
isomorphism of $H_1(M;\z)$ onto $H_1(M';\z)$ which preserves the
torsion linking pairing \cite{Matveev:GeneralizedSurgeries}.  We can
generalize this result to $3$--manifolds with boundaries.  An
$A_2$--surgery preserves the $\mu$--invariant of $\z_2$--homology
$3$--spheres.  The notion of $A_k$--surgery ($k\ge1$) works well also
for spin $3$--manifolds, and an $A_2$--surgery preserves the
$\mu$--invariant of any closed spin $3$--manifolds.  An $A_3$--surgery
preserves the Casson--Walker--Lescop invariant of closed oriented
$3$--manifolds.  Two integral homology $3$--spheres $M$ and $M'$ are
$A_2$-- (resp.\ $A_3$--, $A_4$--) surgery equivalent if and only if they
have equal values of the Rohlin (resp, Casson, Casson) invariant.  For
more informations on $A_k$--surgeries, see below, too.

%

\subsubsection{Definition of new filtrations on $3$--manifolds}
For a closed $3$--manifold $M$, let $\m(M)$ denote the free abelian
group generated by the orientation-preserving diffeomorphism classes
of $3$--manifolds which are $A_1$--equivalent to $M$.  In the following
we will construct a descending filtration
\begin{equation}
\label{8:eq:filtMM}
\m(M)=\m_1(M)\supset \m_2(M)\supset \cdots,
\end{equation}
which we call the {\em $A$--filtration}.  

A graph scheme $S$ in $M$ is said to be {\em allowable} if every
element of $S$ is allowable.  We define the $A$--degree
(resp.\ $S$--degree) of $S$ to be the sum of the $A$--degrees
(resp.\ $S$--degrees) of the elements of $S$.  For an allowable graph
scheme $S=\{S_1,\dots,S_m\}$ in $M$, we define an element $[M,S]$ of
$\m(M)$ by
\begin{displaymath}
[M,S]=\sum_{S'\subset S}(-1)^{|S'|}[M^{\cup S'}],
\end{displaymath}
where the sum is over all subset of $S$, $|S'|$ denotes the number of
elements in $S'$ and $[M^{\cup S'}]$ denotes the
orientation-preserving diffeomorphism class of the result $M^{\cup
S'}$ of surgery on the union $\cup S'$ in $M$.  Then, for each $k\ge
0$, we define $\m_k(M)$ as the subgroup of $\m(M)$ generated by the
elements $[M,S]$, where $S$ is an allowable graph scheme in $M$ of
$A$--degree $k$.

We can prove that the quotient group $\m(M)/\m_{k+1}(M)$ is finitely
generated by showing that there is a descending filtration
\begin{equation}
\label{eq:filtMk}
\bar{\m}_k(M)=\bar{\m}_{k,-1}\supset\bar{\m}_{k,0}\supset\bar{\m}_{k,1}
\supset\cdots\supset \bar{\m}_{k,[k/2]}\supset\{0\}
\end{equation}
on the group
$\bar{\m}_k(M){\overset{\text{def}}{=}}\m_k(M)/\m_{k+1}(M)$ such that
onto each graded quotient $\bar{\m}_{k,l}/\bar{\m}_{k,l+1}$ maps a
finitely generated abelian group ${\mathcal A}^\m_{k,l}(M)$ generated
by $H_1(M;\z)$--labeled uni-trivalent graphs of $A$--degree $k$ and of
$S$--degree $l$.

A homomorphism $f\colon \m(M)\to X$, where $X$ is an abelian group, is
{\em of $A$--type $k$} if $f$ vanishes on $\m_{k+1}(M)$.  Since
$\m(M)/\m_{k+1}(M)$ is a finitely generated abelian group, for a
commutative ring with unit, $R$, the $R$--valued invariants of $A$--type
$k$ form a finitely generated $R$--module.

Claspers enables us to prove realization theorems also for finite type
invariants of $3$--manifolds.  For example we can prove that for a
$3$--manifold $M$, any integral linear combination of connected
$H_1(M,\z)$--labeled uni-trivalent graphs with $k$ trivalent vertices
and with $k-2l$ univalent vertices can be ``realized'' by the
difference of $M$ and a $3$--manifold which is related to $M$ by an
$A_k$--surgery.

It is clear from the definition that if two $3$--manifolds $M$ and $M'$
are $A_{k+1}$--surgery equivalent, then the difference of $M$ and $M'$
lies in $\m_{k+1}(M)$, and hence they are not distinguished by any
invariant of $A$--type $k$.  The converse does not hold in general.  As
with the case of links, we may say that the notion of
$A_{k+1}$--surgery equivalence is more fundamental than the equivalence
relations determined by the $A$--filtration.  However, two integral
homology $3$--spheres are $A_{k+1}$--equivalent if and only if they are
not distinguished by any invariant of $A$--type $k$.  The proof of this
is very similar to that of Theorem~\ref{JkCk}.
Theorem~\ref{UniversalAdditiveVassilievInvariant} can be also
translated into integral homology spheres: We can define the {\em
universal additive $A$--type $k$ invariant} of integral homology
$3$--spheres.

%

\subsubsection{Comparison with other filtrations}
Here we compare the $A$--filtration \eqref{8:eq:filtMM} and other
filtrations in literature.  In
\cite{Garoufalidis-Levine:FiniteTypeBlinks}, S~Garoufalidis and
J~Levine defined a filtration on integral homology spheres using
framed links bounding surfaces, which they call ``blinks.''  This
filtration can be directly generalized to general $3$--manifolds and we
can prove that this filtration equals the $A$--filtration.  For
homology spheres, by a result of Garoufalidis and Levine, this
equality implies that the $A$--filtration is, after re-indexing and
tensoring $\z[\frac{1}{2}]$, equal to T~Ohtsuki's original filtration
using algebraically split framed links
\cite{Ohtsuki:FiniteTypeInvariants}.  Garoufalidis and Levine also
proved that there are no rational invariant of odd degree.  We can
generalize this to that for {\em closed} $3$--manifolds any rational
invariant of $A$--type $2k-1$ is of $A$--type $2k$.  (This cannot be
generalized for $3$--manifolds with boundaries.)

Now we compare the $A$--filtration with the Ohtsuki's filtration on
integral homology $3$--spheres and also with the generalization to more
general $3$--manifolds by T~Cochran and P~Melvin
\cite{Cochran-Melvin:FiniteTypeInvariants}.  Here we call these
filtrations the Ohtsuki--Cochran--Melvin filtrations.  It turns out that
the $3k$th subgroup of the Ohtsuki--Cochran--Melvin filtration is
contained in $\m_{k}(M)$, hence an invariant of $A$--type $k$ is an
invariant of Ohtsuki--Cochran--Melvin type $3k$.  A
$\z[{\frac{1}{2}}]$--module valued invariant of $A$--type $2k$ is an
invariant of Ohtsuki--Cochran--Melvin type $3k$.  Hence the
$A$--filtration is coarser than the Ohtsuki--Cochran--Melvin filtration.
In some respects, the $A$--filtration is easier to handle than the
Ohtsuki--Cochran--Melvin filtration.  Using graph schemes, we can also
re-define the Ohtsuki--Cochran--Melvin filtration.  We can define this
filtration like the $A$--filtration, but, instead of the notion of
$A$--degree, we use that of {\em $E$--degree}, which is defined to be
the number of edges either connecting two nodes or connecting a node
with an unknotted leaf with $-1$ framing not linking with other leaves
nor edges.  This definition enables us to study the
Ohtsuki--Cochran--Melvin filtration using claspers.

Now we compare the notion of $A_k$--surgery equivalence with the notion
of $k$--surgery equivalence introduced by T. D. Cochran, A. Gerges and
K. Orr \cite{Cochran-Gerges-Orr:SurgeryEquivalenceRelations}.  Recall
that two $3$--manifolds $M$ and $M'$ are {\em $k$--surgery equivalent}
to each other if they are related by a finite sequence of Dehn
surgeries on $\pm1$--framed knots whose homotopy classes lie in the
$k$th lower central series subgroups of the fundamental groups of the
$3$--manifolds.  It is easy to see that $2$--surgery equivalence implies
$A_1$--surgery equivalence.  For each $k\ge2$, $A_{2k-2}$--surgery
equivalence implies $k$--surgery equivalence.  However, it is clear
that every integral homology sphere is $k$--surgery equivalent to $S^3$
for all $k\ge2$, while the $A_{2k}$--equivalence becomes strictly finer
for integral homology spheres as $k$ increases.

%

\subsubsection{Examples of invariants of finite $A$--type}
There are many nontrivial invariants of finite $A$--type.  First of
all, we can prove that, for $k\ge0$, the Le--Murakami--Ohtsuki invariant
$\Omega_k$ of closed $3$--manifolds
\cite{Le-Murakami-Ohtsuki:OnAUniversal} is of $A$--type $2k$ (and hence
of Ohtsuki--Cochran--Melvin type $3k$, since any rational invariants of
Ohtsuki--Cochran--Melvin type $3k$ are of $A$--type $2k$).

We can generalize a result of T\,Q\,T Le
\cite{Le:AnInvariantOfIntegral} to rational homology $3$--spheres:
$\Omega_k$ is {\em the universal} rational-valued invariant of
rational homology $3$--spheres of $A$--type $2k$.

S~Garoufalidis and N~Habegger
\cite{Garoufalidis-Habegger:TheAlexanderPolynomial} proved that the
coefficient $C_{2k}$ of $z^{2k}$ in the Conway polynomial of a closed
$3$--manifold with first homology group isomorphic to $\z$ factors
through $\Omega_k$.  Hence $C_{2k}$ is an invariant of $A$--type $2k$.
Recall that $C_{2k}$ is an invariant of Ohtsuki--Cochran--Melvin type
$2k$ \cite{Cochran-Melvin:FiniteTypeInvariants}.

N~Habegger proved that the Le--Murakami--Ohtsuki invariant vanishes for
closed $3$--manifolds with first Betti number $\ge4$
\cite{Habegger:AComputationOfTheUniversal}.  It turns out that for
$3$--manifolds that are $A_1$--equivalent to a fixed closed $3$--manifold
with first Betti number $3k\ge0$, the $\q$--vector space of rational
invariants of $A$--type $2k$ of such manifolds is isomorphic to the
$GL(3k;\z)$--invariant subspace of
${\operatorname{Sym}}^{2k}(\wedge^3V)$, where $V$ is $\q^{3k}$ with
the canonical action of $GL(3k;\z)$.  This invariant subspace is
non-zero, and hence there are nontrivial rational invariants (and
hence integral invariants) of $A$--type $2k$ of closed $3$--manifolds of
first Betti number $3k$ for every $k\ge0$.  These invariants are
homogeneous polynomial of order $2k$ of triple cup products
$\alpha\cup\beta\cup\gamma\in H^3(M;\z)\cong\z$ of
$\alpha,\beta,\gamma\in H^1(M;\z)$ evaluated at the fundamental class
of $M$.  Hence they are of Ohtsuki--Cochran--Melvin type $0$.  For
closed $3$--manifolds with first Betti number $b$, there are no
non-constant rational invariant of $A$--type $k< 2b/3$.

Theory of finite $A$--type invariants suggests that there should be a
refinement of the Le--Murakami--Ohtsuki invariant which does not vanish
for $3$--manifolds with high first Betti numbers and which is universal
among the rational valued finite $A$--type invariants.


\subsection{Groups of homology cobordisms of surfaces}

In Section 5, we proved that for a connected oriented surface
$\Sigma$, the set of $C_k$--equivalence classes of $n$--string links in
${\Sigma\times[0,1]}$ forms a group.  This group plays a fundamental
role in studying the $C_k$--equivalence relations and finite type
invariants of links.  For $A_k$--equivalence relations and finite type
invariants of $3$--manifolds, the {\em group of $A_k$--equivalence
classes of homology cobordisms of a surface} plays a similar role.
This group will serve as a new tool in studying the mapping class
groups of surfaces.

Let $\Sigma$ be a connected compact oriented surface of genus $g\ge0$
possibly with some boundary components.  We set $H=H_1(M;\z)$.

A {\em homology cobordism} $C=(C,\phi)$ of $\Sigma$ is a pair of a
$3$--manifold $C$ and an orientation-preserving diffeomorphism
$\phi\colon \partial (\Sigma\times
[0,1]){\overset{\cong}{\longrightarrow}}\partial C$ such that both the
two inclusions $\phi|_{\Sigma\times[0,1]}\colon
\Sigma\times\{i\}\hookrightarrow C$ for $i=0,1$ induce isomorphisms on
the first homology groups with integral coefficients.  Two homology
cobordisms $(C,\phi)$ and $(C',\phi')$ are said to be {\em equivalent}
if there is an orientation-preserving diffeomorphism $\Phi\colon
C{\overset{\cong}{\longrightarrow}} C'$ such that
$\phi'=(\Phi|_{\partial C})\phi$.  For two homology cobordisms
$C_1=(C_1,\phi_1)$ and $C_2=(C_2,\phi_2)$, the {\em composition}
$C_1C_2=(C_1,\phi_1)(C_2,\phi_2)$ is defined by ``pasting the bottom
of $C_1$ and the top of $C_2$.''  The set of equivalence classes of
homology cobordisms of $\Sigma$, $\C(\Sigma)$ forms a monoid with
multiplication induced from the composition operation defined above,
and with unit the equivalence class of the {\em trivial} homology
cobordism $1_\Sigma=({\Sigma\times[0,1]},\operatorname{id}_{\partial
({\Sigma\times[0,1]})})$.

A homology cobordism $C$ is {\em homologically trivial} if, for the
two embeddings 
\begin{displaymath}
i_\epsilon\colon
\Sigma\overset{\cong}{\longrightarrow}\Sigma\times\{\epsilon\}
\hookrightarrow C,\quad(\epsilon=0,1),
\end{displaymath}
the composition $(i_1)_*^{-1}(i_0)_*\colon H\to H$ of the induced
isomorphisms is the identity.  Let $\C_1(\Sigma)$ denote the submonoid
of $\C(\Sigma)$ consisting of the 
equivalence classes of homologically trivial cobordisms of $\Sigma$.

For each $k\ge 1$, we define the notion of $A_k$--equivalence of
homology cobordisms in the obvious way.  For $k\ge1$, let
$\C_k(\Sigma)$ denote the submonoid of $\C(\Sigma)$ consisting of the
equivalence classes of homology cobordisms that are $A_k$--equivalent
to the trivial cobordism $1_\Sigma$.  This defines a descending
filtration on $\C_1(\Sigma)$,
\begin{equation}
\label{8:eq:c1c2}
\C_1(\Sigma)\supset \C_2(\Sigma)\supset\cdots.
\end{equation}

We can prove that the two definitions of $\C_1(\Sigma)$ are
equivalent, ie, a homology cobordism of $\Sigma $ is homologically
trivial if and only if it is $A_1$--equivalent to $1_\Sigma$.

\begin{figure}[ht!]
\cl{\includegraphics{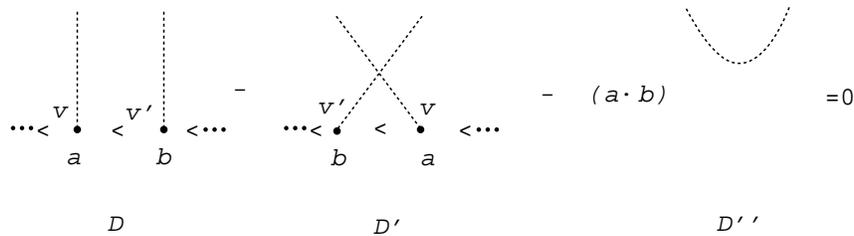}}
\caption{Let $D$ be a uni-trivalent graph and let $v<v'$ be two
consecutive univalent vertices in $D$ labeled $a,b\in H_1(\Sigma;\z)$.
Let $D'$ be the uni-trivalent graph obtained from $D$ by exchanging
the order of $v$ and $v'$.  Let $D''$ denote the uni-trivalent graph
obtained from $D$ by connecting two vertices $v$ and $v'$.  Then the
``STU--like relation'' states that $D-D'-(a\cdot b)D''=0$, where
$a\cdot b\in\z$ denote the intersection number of $a$ and $b$.  In
this figure the univalent vertices are placed according to the total
order.}
\label{STUlike}
\end{figure}

\begin{figure}[ht!]
\cl{\includegraphics{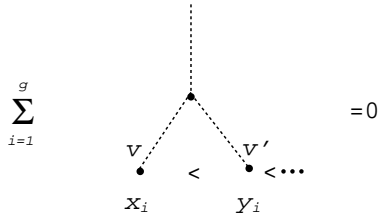}}
\caption{This relation states that $\sum_{i=1}^{g}D_{x_i,y_i}=0$,
where the elements $x_1,y_1,\ldots,x_g,y_g$ form a symplectic basis of
$H_1(\Sigma;\z)$, and $D_{x_i, y_i}$ is a uni-trivalent graph with the 
smallest two univalent vertices $v$ and $v'$ adjacent to the same
trivalent vertex such that $v$ and $v'$ are labeled $x_i$ and $y_i$,
respectively.   This relation does not depend on the choice
of the symplectic basis.}
\label{omega}
\end{figure}

Now we consider the descending filtration of quotient monoids by the
$A_{k+1}$--equivalence relation
\begin{equation}
\label{eq:filtcsigmaak}
\C(\Sigma)/A_{k+1}\supset \C_1(\Sigma)/A_{k+1}\supset\dots\supset
\C_k(\Sigma)/A_{k+1}.
\end{equation} 
These monoids are {\em finitely generated groups}, and moreover
$\C_i(\Sigma)/A_{k+1}$ is nilpotent for $i=1,\dots, k$.  Especially,
$\bar\C_k(\Sigma){\overset{\text{def}}{=}}\C_k(\Sigma)/A_{k+1}$ is an
abelian group.  We define, when $\Sigma$ is not closed and $k\ge2$, a
finitely generated abelian group ${\mathcal A}_k(\Sigma)$ generated by
allowable $H$--labeled uni-trivalent graphs of $A$--degree $k$ on the
empty $1$--manifold equipped with a total order on the set of univalent
vertices.  Here an $H$--labeled uni-trivalent graph $D$ is {\em
allowable} if each components of $D$ has at least one trivalent
vertex.  These uni-trivalent graphs are subject to the antisymmetry
relations, the IHX relations, the ``STU--like relations'' and the
multilinearity of labels.  Here the ``STU--like relation'' is depicted
in Figure~\ref{STUlike}.  When $\Sigma$ is closed and $k\ge2$, we
define $A_k(\Sigma)$ to be the quotient of
$A_k(\Sigma\setminus\operatorname{int}D^2)$ by the relation depicted
in Figure~\ref{omega}.  When $\Sigma$ is not closed and $k=1$, we set
${\mathcal A}_1(\Sigma)=\wedge^3H\oplus \wedge^2H_2\oplus H_2\oplus
\z_2$, where we set $H_2=H_1(\Sigma;\z_2)\cong H\otimes\z_2$.  When
$k=1$ and $\Sigma$ is closed, we set ${\mathcal
A}_1(\Sigma)=\wedge^3H/(\omega\wedge H)\oplus
\wedge^2H_2/(\omega_2)\oplus H_2\oplus \z_2$, where
$\omega=\sum_{i=1}^g x_i\wedge y_i \in \wedge^2H$ for a symplectic
basis $x_1,y_1,\dots,x_g,y_g\in H$, and $\omega_2$ is the
$\operatorname{mod} 2$ reduction of $\omega$.

There is a natural surjective homomorphism of
$A_k(\Sigma)$ onto $\bar\C_k(\Sigma)$.  We conjecture that this is an
isomorphism.  This conjecture holds when $k=1$.  We can
also prove this conjecture over $\q$ for $k\ge1$ with $\Sigma$
non-closed, 
using the Le--Murakami--Ohtsuki invariant.

We can naturally define a graded Lie algebra structure on the graded
abelian group $\bar\C(\Sigma){\overset{\text{def}}{=}}
\bigoplus_{k=1}^\infty\bar\C_k(\Sigma)$.  When $\Sigma$ is not closed,
we can give a presentation of the Lie algebra
$\bar\C(\Sigma)\otimes\q$ in terms of uni-trivalent graphs.  (Again, the proof requires the
Le--Murakami--Ohtsuki invariant.)

Groups and Lie algebras of homology cobordisms of surfaces will serve
as {\em new tools in studying the mapping class groups of surfaces}.
This is because we can think of a self-diffeomorphism of a surface
$\Sigma$ as a homology cobordism of $\Sigma$ via the mapping cylinder
construction.  The filtration \eqref{8:eq:c1c2} restricts to a
filtration on the Torelli group of $\Sigma $, which is coarser than or
equal to the lower central series of the Torelli group\footnote{At low
genus we can prove that they are different, but at high genus it is
open if they are different or not.  We conjecture that they are stably
equal.} and is finer than the filtration given by considering the
action of the Torelli group on the fundamental group $\pi_1\Sigma$
\cite{Johnson:AnAbelianQuotient}\cite{Morita:AbelianQuotients}.  We
can naturally extend the Johnson homomorphisms to homologically
trivial cobordisms and describe it in terms of tree claspers and
tree-like uni-trivalent graphs.  It is extremely important to clarify
the relationships between the presentation of the Lie algebra
$\bar\C(\Sigma)\otimes\q$ in terms of uni-trivalent graphs and
R~Hain's presentation of the associated graded of the lower central
series of the Torelli group \cite{Hain:Infinitesimal}.

\subsection{Claspers and gropes}
\label{ss:grope}

Some authors use {\em gropes} to study links and $3$--manifolds
\cite{Cochran-Gerges-Orr:SurgeryEquivalenceRelations}
\cite{Krushkal:AdditivityProperties}.  We explain here some
relationships between claspers and gropes {\em embedded} in
$3$--manifolds.

For the definitions of {\em gropes} and {\em capped gropes}, see
\cite{Freedman-Teichner:4ManifoldTopologyII}.  We define a {\em
(capped) $k$--grope} $X$ for a link $\gamma$ in $M$ to be a (capped)
grope $X$ of class $k$ embedded in $M$ intersecting $\gamma$ only by
some transverse double points in the caps of $X$.  (In the non-capped
case, $X$ and $\gamma$ are disjoint.)

Two links $\gamma$ and $\gamma'$ in $M$ are said to be related by a
{\em (capped) $k$--groping} if there is a (capped) $k$--grope $X$ for
$\gamma$ and a band $B$ connecting a component of $\gamma$ and the
bottom $b$ of $X$ in such a way that $B\cap X=\partial B\cap b$ and
$B\cap\gamma=\partial B\cap \gamma$, and if the band sum of $\gamma$
and $b$ along the band $B$ is equivalent to $\gamma'$.

We can prove that two links in $M$ are related by a sequence of capped
$k$--gropings (resp.\ $k$--gropings) if and only if they are
$C_k$--equivalent (resp.\ $A_k$--equivalent).  As corollaries to this, we
can prove that an $A_k$--move on a link in $M$ preserves the homotopy
classes of the components of a link up to the $k$th lower central
series subgroup of $\pi_1M$, and that the $k$th nilpotent quotient
(ie, the quotient by the $k+1$st lower central series subgroup) of
the fundamental group of the link exterior is an invariant of
$A_k$--equivalence classes of links (and hence of $C_k$--equivalence
classes).  From this we can also prove that an $A_k$--surgery on a
$3$--manifold preserves the $k$th nilpotent quotient of the fundamental
group of $3$--manifolds.

Recall that for a knot $\gamma$ in a $3$--manifold $M$, the homotopy
class of $\gamma$ lies in the $k$th lower central series subgroup of
$\pi_1M$ if and only if there is map $f$ of a grope $X$ of class $k$
into $M$ such that the bottom of $X$ is mapped diffeomorphically onto
$\gamma$.  This condition is much weaker than that $\gamma$ bounds an
embedded $k$--grope in $M$.  In some sense, embedded gropes, and hence
tree and graph claspers, may be thought of as a kind of ``{geometric
commutator}'' in a $3$--manifold.  Gropes thus provide us  another way of
thinking of calculus of claspers as a commutator calculus of a new
kind.

\newpage


\end{document}